\documentclass[11pt]{amsart}
\usepackage[a4paper,margin=2cm]{geometry}
\usepackage[comma, round]{natbib}

\usepackage{float}
\usepackage{amssymb}
\usepackage{amsthm}
\usepackage{amsmath}
\usepackage{stmaryrd}
\usepackage{titletoc}
\usepackage{mathrsfs}
\usepackage[colorlinks, allcolors=blue, linktocpage]{hyperref}
\usepackage{amsfonts}
\usepackage{graphicx}
\usepackage[linesnumbered, ruled, vlined]{algorithm2e}
\usepackage{esint}
\usepackage{bbm}
\usepackage{bm}
\usepackage{mathtools}
\usepackage[shortlabels, inline]{enumitem}

\setlist[enumerate,1]{label=\arabic*.}
\usepackage{todonotes}
\usepackage{color}

\numberwithin{equation}{section}
\numberwithin{figure}{section}

\theoremstyle{plain}
\newtheorem{thm}{Theorem}[section]
\newtheorem{cor}{Corollary}[section]
\newtheorem{lem}{Lemma}[section]
\newtheorem{prop}{Proposition}[section]

\theoremstyle{definition}
\newtheorem{defn}{Definition}[section]
\newtheorem{asmp}{Assumption}[section]
\newtheorem{rmk}{Remark}[section]
\newtheorem{exam}{Example}[section]

\newcommand{\cA}{\mathcal{A}}
\newcommand{\cB}{\mathcal{B}}

\newcommand{\cD}{\mathcal{D}}
\newcommand{\cE}{\mathcal{E}}
\newcommand{\cF}{\mathcal{F}}
\newcommand{\cG}{\mathcal{G}}
\newcommand{\cH}{\mathcal{H}}

\newcommand{\cK}{\mathcal{K}}

\newcommand{\cM}{\mathcal{M}}
\newcommand{\cN}{\mathcal{N}}

\newcommand{\cP}{\mathcal{P}}

\newcommand{\cS}{\mathcal{S}}
\newcommand{\cT}{\mathcal{T}}

\newcommand{\cW}{\mathcal{W}}
\newcommand{\cX}{\mathcal{X}}

\newcommand{\cZ}{\mathcal{Z}}

\newcommand{\bbE}{\mathbb{E}}

\newcommand{\bbP}{\mathbb{P}} 
\newcommand{\bbQ}{\mathbb{Q}} 
\newcommand{\bbR}{\mathbb{R}} 
\newcommand{\bbS}{\mathbb{S}} 
\newcommand{\bbT}{\mathbb{T}}

\newcommand{\bbW}{\mathbb{W}}

\newcommand{\scrA}{\mathscr{A}}

\newcommand{\scrH}{\mathscr{H}}

\newcommand{\scrK}{\mathscr{K}} 
 
\newcommand{\scrM}{\mathscr{M}}

\newcommand{\scrZ}{\mathscr{Z}}

\newcommand{\frD}{\mathfrak{D}}

\newcommand{\frF}{\mathfrak{F}}

\newcommand{\Id}{\operatorname{Id}}

\newcommand{\Law}{\operatorname{Law}}

\newcommand{\var}{\mathrm{var}}

\newcommand{\tr}{\operatorname{tr}}

\newcommand{\op}{\mathrm{op}}

\def\lb{\mathopen{}\mathclose\bgroup\left}
\def\rb{\aftergroup\egroup\right}

\newcommand{\identity}{\text{Id}}

\newcommand{\tb}{\tilde{b}}

\newcommand{\tsigma}{\tilde{\sigma}}
\newcommand{\tZ}{\tilde{Z}}

\newcommand{\ocH}{\overline{\scrH}}

\newcommand{\optproj}[2]{\prescript{o}{}{#2}^{#1}}
\newcommand{\jointoptproj}[2]{\prescript{\mathbf{o}}{}{#2}^{#1}}

\newcommand{\Leb}{\text{Leb}}
\newcommand{\onevar}{1\text{-var}}

\usepackage[]{todonotes}

\title{Adapted transport inequalities in continuous time}

\author{Erhan Bayraktar}
\address{Department of Mathematics, University of Michigan, Ann Arbor, MI 48109, USA}
\email{erhan@umich.edu}

\author{Fang Rui Lim}
\address{Department of Mathematics, University of Michigan, Ann Arbor, MI 48109, USA}
\email{ruilim@umich.edu}

\begin{document}

\bibliographystyle{abbrvnat}



\begin{abstract}
    Given a strong solution $X$ to a stochastic differential equation (SDE), we establish sufficient conditions under which $\Law(X)$ satisfies an adapted quadratic transport cost inequality (AQTCI) with respect to the semimartingale adapted Wasserstein distance.
    When the drift and diffusion coefficients are Markovian, we further show that the AQTCI tensorizes in a dimension-free manner. As a corollary, we show that various stochastic control problems, e.g. optimal stopping and mean-variance portfolio selection, are stable under bounded, measurable drift perturbations. Additionally, we derive concentration inequalities (with respect to certain dynamic risk measures) for stochastic integrals $\int H_s dX_s$ and deduce appropriate bounds on their tail behavior. Furthermore, using the functional It\^o formula, we show that the AQTCI for $\Law(X)$ allows one to deduce classical transport inequalities for $\Law(F(X))$ even when $F$ is not Lipschitz on path space. Finally, we give sufficient conditions for $\Law(X)$ to satisfy AQTCIs  when $X$ is a solution to a scalar fractional SDE or stochastic Volterra equation.
\end{abstract}

\maketitle

\tableofcontents
\newpage

\section{Introduction}

Our starting point is the following representation of relative entropy as a transport problem established in \cite{lassalle2018causaltp}: \begin{equation}
    \label{eqn:lassalletptinequality}
    2 \cH(\bbP| \bbW^n) = \inf_{\pi \in \Pi_{bc}(\bbW^n,\bbP)} \bbE^\pi c_{\mathfrak{h}}.
\end{equation} Here, $\bbW^n$ is the law of an $n$-dimensional standard Brownian motion, $\bbP$ is any probability measure on $\cW^n$, the space of $\bbR^n$-valued continuous functions on $[0,1]$ starting at $0$, $\Pi_{bc}(\bbW^n, \bbP)$ is the set of all \emph{bicausal} couplings between $\bbW^n$ and $\bbP$, $\cH$ denotes the relative entropy, i.e. \[
\cH(\bbP | \bbW^n) := \begin{cases}
    \bbE^\bbP \log \frac{d\bbP}{d\bbW^n} & \text{ if } \bbP \ll \bbW^n, \\
    +\infty & \text{ otherwise,}
\end{cases} 
\] and $c_{\mathfrak{h}} : \cW^n \times \cW^n \rightarrow \bbR_{\geq 0} \cup \{+\infty\}$ denotes the Cameron-Martin cost given by \[
c_{\mathfrak{h}}(\omega, \omega') := |\omega - \omega'|^2_{\mathfrak{h}}, \text{ where } |\omega|^2_{\mathfrak{h}} := \begin{cases}
    \int_0^1 |\dot{\omega}_s|^2 ds & \text{ if } \omega \text{ is absolutely continuous,} \\
    +\infty & \text{ otherwise.}
\end{cases}
\] Since the uniform norm $|\cdot|_\infty$ is trivially bounded by $|\cdot|_{\mathfrak{h}}$ on $\cW^n$, we immediately obtain from \eqref{eqn:lassalletptinequality} an adapted version of the quadratic transport cost inequality: \begin{equation}
    \label{eqn:AQTCIBM}
    \cA\cW^2_{|\cdot|_\infty,2}(\bbP, \bbW^n) := \inf_{\pi \in \Pi_{bc}(\bbW^n,\bbP)} \bbE^\pi |\cdot - \cdot|^2_\infty \leq 2 \cH(\bbP| \bbW^n),
\end{equation} where $\cA\cW_{|\cdot|_\infty,2}$ is the so-called adapted Wasserstein distance.
In this article, we seek to generalize this beyond $\bbW^n$ and consider its implications.


Classical transport inequalities for a probability measure $\bbQ$ on $\cW^n$ typically take the form: \begin{equation}
    \label{eqn:classicaltptinequality}
    \cW_{d,p}(\bbP, \bbQ) \leq \sqrt{C \cH(\bbP|\bbQ)} \text{ for all probability measures } \bbP \text{ on } \cW^n, 
\end{equation} where $p \geq 1$ and $\cW_{d,p}$ is the $p^{th}$-Wasserstein distance built from some metric $d$ on $\cW^n$. An example of such a $\bbQ$ satisfying \eqref{eqn:classicaltptinequality}  when $d = |\cdot|_\infty$ and $p = 2$ is $\mu^{\sigma,b}$, the law of the stochastic differential equation (SDE): \[
dX_t = b(t,X) dt + \sigma(t,X) dB_t, \quad X_0 = 0, \quad \Law(X) =: \mu^{\sigma,b}.
\] See \cite{pal2012concentration,ustunel2012transportationci} for the case when $b,\sigma$ are Lipschitz and $\sigma$ is bounded or \cite{bahlalimouchtabihtangpi2021quadraticti,palsarantsev2019transportationci} for more general conditions on $b,\sigma$. These inequalities and their many fascinating implications, such as concentration-of-measure phenomena, have been studied since the 1990s and continue to be a fruitful area of research to this day; we refer readers to the wonderful survey by \cite{gozlanleonard2010transport} and \cite{ledoux2001conc}. However, it was observed starting from the works of \cite{aldous1981weakcvg,hooverkeisler1984adapted,hellwig1996seq} that the usual weak convergence of probability measures, metrized by the Wasserstein distance, is agnostic to the filtrations generated by stochastic processes, thereby limiting its usefulness in applications to stochastic analysis and mathematical finance. Hence, one might wish to consider a stronger distance that incorporates this temporal information, namely the adapted Wasserstein distance. This distance is induced by causal transports, which were introduced by \cite{lassalle2018causaltp,backhoffbeiglbocklinzalashko2017causal} and subsequently studied by many others, e.g. \cite{backhoffbartlbeiglbockeder2020all,bartlbeiglbockpammer2026,acciaio2020causalot,bayraktarhan2025fittedvalue,backhoffbayraktarekrenzitridis2026analytical} to name a few.
It is then natural to consider the following variant of the transport inequality: \begin{equation}
    \label{eqn:adaptedtptinequality}
     \cA\cW_{d,p}(\bbP, \bbQ) \leq \sqrt{C \cH(\bbP|\bbQ)} \text{ for all probability measures } \bbP \text{ on } \cW^n.
\end{equation} Here, $\cA\cW_{d,p}$ denotes the adapted $p^{th}$-Wasserstein distance, whose precise definition is deferred to Section \ref{subsection:adaptedwassersteindistance}. Since $\cW_{d,p} \leq \cA\cW_{d,p}$, we trivially have $\eqref{eqn:adaptedtptinequality} \implies \eqref{eqn:classicaltptinequality}$. 

Going a step further, owing to the semimartingale nature of $\mu^{\sigma,b}$, we can even consider an adapted distance tailored specifically to semimartingales developed in \cite{backhoff2020adaptedwd}: \begin{equation}
    \label{eqn:adaptedsmdistance}
    \cA\cW_{sm,p}(\bbP,\bbQ) := \inf_{\Law(X,Y) \in \Pi_{bc}(\bbP,\bbQ)} \bbE (|V^\bbP - V^\bbQ|_{1-\var}^p + \tr([M^\bbQ - M^\bbP]_1)^{p/2})^{1/p}.
\end{equation} Here, $X,Y$ are continuous $\bbR^n$-valued semimartingales w.r.t. their own filtrations $\cF^X,\cF^Y$ on some abstract probability space with finite variation parts $V^\bbP,V^\bbQ$, martingale parts $M^\bbP, M^\bbQ$ and laws $\bbP,\bbQ$ respectively. Crucially, this distance only makes sense due to the (bi-)causality constraint on the couplings, yielding that $X$ and $Y$ are actually semimartingales w.r.t. the filtration \emph{jointly} generated by $X$ and $Y$, $\cF^{X,Y}$. In any case, we can consider the following adapted transport inequality: \begin{equation}
    \label{eqn:adaptedtptinequalitysemimart}
     \cA\cW_{sm,p}(\bbP, \bbQ) \leq \sqrt{C \cH(\bbP|\bbQ)} \text{ for all probability measures } \bbP \text{ on } \cW^n.
\end{equation} By the BDG inequality, we have $\cA\cW_{|\cdot|_\infty,p} \leq c_p \cA\cW_{sm,p}$ for some constant $c_p > 0$ and consequently, $\eqref{eqn:adaptedtptinequalitysemimart}\implies \eqref{eqn:adaptedtptinequality} \implies \eqref{eqn:classicaltptinequality}$. We will demonstrate in the sequel that $\mu^{\sigma,b}$ satisfies \eqref{eqn:adaptedtptinequalitysemimart}, the strongest of the aforementioned transport inequalities, and investigate its consequences.

In discrete time path space, sufficient conditions for the law of a discrete-time stochastic process to satisfy \eqref{eqn:adaptedtptinequality} were given in \cite{backhoffbeiglbocklinzalashko2017causal}. Surprisingly, \cite{park2026t1} showed that $\eqref{eqn:classicaltptinequality} \iff \eqref{eqn:adaptedtptinequality}$ for $p=1$. To the best of the authors' knowledge, with the exception of \cite{follmer2022optimal} for a related transport inequality using the notion of specific relative entropy,  adapted transport inequalities in continuous time have been unexplored thus far, which is the focus of the present article.




\subsection{Contributions} 

\subsubsection*{Adapted transport inequalities for diffusions} We provide sufficient conditions in Theorem \ref{thm:adaptedQTCIsmcost} under which $\mu^{\sigma,b}$ satisfies \eqref{eqn:adaptedtptinequalitysemimart} for $p = 2$, hereafter referred to as the semimartingale adapted quadratic transport cost inequality (AQTCI). In order to do so, we first characterize bicausal couplings between probability measures under a finite relative entropy condition in Theorem \ref{thm:bicausalcplsbetweensemimartfiniterelent} by extending the martingale representation result of \cite{ustunel2019mrp} in Appendix \ref{appendix:mrpfiniterelent}. 

It is known that many quantities from stochastic analysis and control that are not stable w.r.t. Wasserstein distances are nevertheless stable w.r.t. the adapted Wasserstein distances. For example, consider the case of utility maximization under model approximation studied by \cite{bayraktardolinskyguo2020continuity} in frictionless markets and by \cite{bayraktardolinskyidolinsky2020extended} in markets with proportional transaction costs. The former established continuity under weak convergence only after assuming additional structural and integrability conditions, while the latter uses the extended weak convergence \cite{aldous1981weakcvg}, a notion intimately connected with adapted transports \cite{backhoffbartlbeiglbockeder2020all}. We collect in Example \ref{exam:awlipschitzbd} an (incomplete) list of such quantities in the literature and supplement two additional quantities: optional projections (Appendix \ref{appendix:optprojstab}) and the value function of mean-variance portfolio selection problems (Appendix \ref{appendix:mvpselect}). Following a calculation analogous to \citeauthor{follmer1986timereversal}'s (\citeyear{follmer1986timereversal}) seminal work establishing the connection between relative entropy and the $L^2$-integral of the drift term of semimartingales, we show that the semimartingale AQTCI for $\mu^{\sigma,b}$ implies the stability of these quantities w.r.t. bounded, measurable perturbations of the drift coefficient $b$ (Corollary \ref{cor:stabwrtuniformdrift}).

\subsubsection*{Concentration phenomena and related results}
We juxtapose the concentration phenomena implied by classical transport inequalities with the corresponding results for our adapted variant:
\begin{enumerate}
\item (Tail behaviors) The results of \cite{marton1996bounding} and \cite{bobkovgotze1999exponentiali} imply that \eqref{eqn:classicaltptinequality}, with $p = 1$, for probability measure $\bbQ$ gives a quantitative bound on the deviation of any $1$-Lipschitz function $F : \cW^n \rightarrow \bbR$ from its mean under $\bbQ$, i.e. \begin{equation}
    \label{eqn:lipdeviation}
    \bbQ(|F - \bbE^\bbQ F| \geq r) \leq 2 \exp(- \frac{r^2}{C}) \text{ for } r > 0.
    \end{equation} In other words, the tails of centered Lipschitz functions exhibit sub-Gaussian behavior. We show that if $\bbQ$ is a semimartingale measure satisfying the semimartingale AQTCI, then one can derive similar tail bounds for stochastic integrals $\int_0^\cdot H_s dX_s$ w.r.t. any semimartingale $X$ with law $\bbQ$. Here, $H = (H_t)_{t \in [0,1]}$ may be bounded (Proposition \ref{prop:tailbehaviorforstointsemimartbddintegrand}) or square integrable but H\"older continuous (Proposition \ref{prop:tailbehaviorforstointmartsqintintegrand}). In particular, if $X$ is a martingale, we recover via optimal transport methods the well-known result in stochastic calculus that martingales $M = (M_t)_{t \in [0,1]}$ with bounded quadratic variation exhibit sub-Gaussian tail behavior.
    \item (Tensorization) One can also deduce that \eqref{eqn:classicaltptinequality} with $p = 2$ implies a dimension-free tensorization, i.e. for all $N \geq 1$, \[
     \cW_{d_{\cW^n}^{\oplus N}, 2}(\bbP, \bbQ^{\otimes N}) \leq \sqrt{C \cH(\bbP|\bbQ^{\otimes N})} \text{ for all probability measures } \bbP \text{ on } {\cW^n}^N,
    \] where $\bbQ^{\otimes N}$ is the product measure of $\bbQ$ on ${\cW^n}^N$, $d_{{\cW^n}}^{\oplus N}(x,x') := (\sum_{i = 1}^N d^2_{\cX}(x_i,x'_i))^{1/2}$ and crucially, the constant $C$ is \emph{independent} of $N$. Such a property, initially established by \cite{talagrand1996transportation} for Gaussian distributions on $\bbR^N$, implies a dimension-free concentration of measure phenomenon for $\bbQ^{\otimes N}$. In fact, \cite{gozlan2009characterization} demonstrates that \eqref{eqn:classicaltptinequality} with $p = 2$ is equivalent to dimension-free concentration of measure. In Corollary \ref{cor:adaptedQTCItensorization}, we show that the semimartingale AQTCI for probability measure $\mu^{\sigma,b}$ can be tensorized in a dimension-free manner when $b,\sigma$ are Markovian. This, in turn, implies a quantitative propagation of chaos for mean field systems as studied in \cite{lacker2023hierarchies} under this adapted Wasserstein metric for semimartingales.
    \item (Transferring transport inequalities) A particularly convenient way to establish \eqref{eqn:classicaltptinequality} is via \cite[Lemma $2.1$]{djelloutguillinwu2004transportationci}, which allows one to ``transfer" transport inequalities via Lipschitz maps: if $\bbQ$ satisfies \eqref{eqn:classicaltptinequality}, then so does $F_{\#}\bbQ$ for $1$-Lipschitz $F$. This was used in \cite{bartltangpi2020functional} very effectively to establish classical transport inequalities for forward and backward diffusions. We show that if semimartingale measure $\bbQ$ satisfies the semimartingale AQTCI, then one can ``transfer" this to a classical QTCI for its quadratic variation (Corollary \ref{cor:QTCIforQV}) or along non-anticipative functionals $F$ that are horizontally and vertically differentiable in the sense of \cite{Dupire2019} (Proposition \ref{prop:W1TCIfornonanticipativefnls}). We emphasize here that neither the map sending a semimartingale to its quadratic variation nor $F$, under our proposed regularity assumptions, is necessarily Lipschitz.
    \item (Liquidity risk profiles) As an application with a more mathematical finance flavor, we consider liquidity risk profiles of a random variable $L$ representing a financial loss w.r.t. a convex risk measure $\cE$. Introduced by \cite{lacker2018liquidity}, the liquidity risk profile is the curve $(\cE(\lambda L))_{\lambda \geq 0}$ which models how the risk of $L$ scales with its size $\lambda$. For  certain dynamic risk measures $\cE$, \cite{tangpi2019concentration} provided a quadratic upper bound for the liquidity risk profile of losses $L$ given by Lipschitz functionals of Brownian motion or its increments. As a consequence of the semimartingale AQTCI, we derive quadratic upper bounds for financial losses $L$ represented by stochastic integrals instead (Proposition \ref{prop:convexriskmeasureconcinequality}).
\end{enumerate}  

\subsubsection*{Transport inequalities beyond diffusions on path space} Finally, classical QTCIs, i.e. \eqref{eqn:classicaltptinequality} with $p = 2$ and appropriately chosen $d$, are known to hold for laws of stochastic processes that are neither Markovian nor semimartingale; see \cite{fan2015stochastic} for certain stochastic Volterra equations (SVEs) and \cite{saussereau2012transportation} for differential equations driven by fractional Brownian motion (fBM) or \cite{riedel2017transportation} for more general Gaussian rough paths. We show that AQTCIs in the form of \eqref{eqn:adaptedtptinequality} can also be established for these types of processes. In order to do so, we first give an abstract generalization of \eqref{eqn:lassalletptinequality} in Proposition \ref{prop:relentequality} for probability measures $\bbQ$ that can be transformed into $\bbW^n$ in a bicausal, bijective manner. This bicausal, bijective transformation is analogous to the transfer principle used in \cite{jianglim2025transferprinciple} to compute adapted Wasserstein distances between Gaussian processes in continuous time. By first identifying appropriate conditions for SDEs driven by fractional Brownian motion or SVEs so that their laws satisfy the assumptions of Proposition \ref{prop:relentequality}, we prove the AQTCI w.r.t. the uniform or $L^2$-metric on $\cW^n$ for these laws in Theorems \ref{thm:adaptedQTCIfSDEs} and \ref{thm:adaptedQTCIforSVEs}.  

\begin{figure}[t!]
    \centering 
\includegraphics{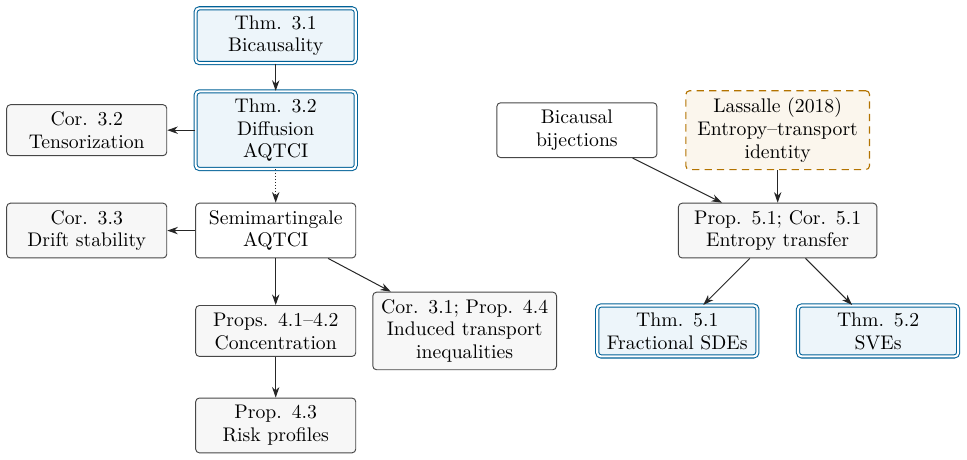}
\caption{Dependence of the main results. Double borders mark theorems; dashed borders mark an external input. The dotted arrow supplies a diffusion instance of the general semimartingale AQTCI.}
\label{fig:dependencies}
\end{figure}

\subsubsection*{Outline} 

Section \ref{section:preliminaries} introduces notations and the adapted Wasserstein distance. Section \ref{section:AQTCIfordiffusions} discusses the main result establishing the semimartingale AQTCI for laws of SDEs and the tensorization and stability corollaries specific to this diffusion setup. Section \ref{section:concentration} discusses general concentration of measure phenomena as a result of the semimartingale AQTCI, namely tail behaviors, transferring transport inequalities and liquidity risk profiles. Section \ref{section:AQTCIbeyonddiffusions} establishes the AQTCI for the laws of some non-Markovian, non-semimartingale scalar stochastic equations. We collect in Appendix \ref{section:appendix} the supporting results on martingale representation and bicausal, a.s. bijective maps. It also proves the stability of optional projections and mean-variance portfolio selection.
Figure \ref{fig:dependencies} shows how the results depend on one another.

\section{Preliminaries}
\label{section:preliminaries}

\subsection{Notation}

Given matrix $A \in \bbR^{n \times n'}$, $|A|$ always denotes the Frobenius norm, so that if $n' = 1$, then this is simply the Euclidean norm of vector $A$. When necessary, we denote by $|A|_{\op}$ the operator norm. 

On any abstract probability space $(\Omega, \cF, \bbP)$ with some random variable $M$ taking values in some measurable space $(\Omega', \cF')$, the probability measure $M_{\#}\bbP$ on $\Omega'$ denotes the pushforward of $\bbP$ under $M.$ We define $\cM_b([0,1], \bbR^n)$ to be the space of bounded, measurable $\bbR^n$-valued functions and $L^p(\bbP)$ to be the standard $L^p$ space of $\cF$-measurable real-valued functions $f : \Omega \rightarrow \bbR$ with norm $\Vert \cdot \Vert_{L^p(\bbP)}$. On occasion, we also consider spaces $L^p(\bbP; \scrZ)$, where $\scrZ$ is a Banach space with norm $\Vert \cdot \Vert_{\scrZ}$ and Borel $\sigma$-algebra $\cZ$, to be the space of $\cF\backslash \cZ$-measurable functions $f : \Omega \rightarrow \scrZ$ such that $\bbE^\bbP(\Vert f \Vert_{\scrZ}^p) < \infty$. $\Leb$ will always denote the Lebesgue measure on $[0,1]$. 

Throughout, $\cW^n$ denotes the space of $\bbR^n$-valued continuous functions on $[0,1]$ starting from $0$ endowed with the topology induced by the supremum norm $|\cdot|_\infty$ and associated Borel $\sigma$-algebra $\cF^X_1$. $X = (X_t)_{t \in [0,1]}$ will always denote the canonical process on $\cW^n$ and its \emph{right-continuous} canonical filtration $\cF^X = (\cF^X_t)_{t \in [0,1]}$. On the product space $\cW^n \times \cW^n$, without risk of confusion, the first marginal process is still denoted $X : \cW^n \times \cW^n \times [0,1] \rightarrow \bbR^n, X_t(\omega, \omega') := \omega(t)$, while $Y = (Y_t)_{t \in [0,1]}$ denotes the second marginal, i.e. $Y_t(\omega, \omega') := \omega'(t)$. Furthermore, $\cP(\cW^n)$ denotes the space of probability measures on $\cW^n$ with $\cP(\cW^n \times \cW^n)$ similarly defined.

We say that $\bbP \in \cP(\cW^n)$ is a semimartingale probability measure if, 
under $\bbP$, $X$ is an $\bbR^n$-valued semimartingale with respect to the canonical filtration $\cF^X$ with characteristics absolutely continuous w.r.t. the Lebesgue measure. In particular, $X$ has decomposition \[
X = M^\bbP + V^\bbP,
\] where $M^\bbP,V^\bbP : \cW^n \rightarrow \cW^n$ are such that $M^\bbP$ is a continuous local martingale w.r.t. $\cF^X$ and $V^\bbP$ a finite-variation $\cF^X$-predictable continuous process satisfying $[M^\bbP], |V^\bbP|_{\onevar}\ll \Leb$. Here, $|V^\bbP|_{1\text{-}\var}$ denotes the $1$-variation: \[
|V^\bbP|_{1\text{-}\var} := \sup_{\substack{0= t_0 < t_1 < \dots < t_n = 1 \\
n \geq 1}}  \sum_{i = 0}^{n-1} |V^\bbP_{t_{i + 1}} - V^\bbP_{t_i}|
\] and $[M^\bbP] =([M^\bbP]_t)_{t \in [0,1]}$ the quadratic variation of a martingale, \[
[M^\bbP]_t \in \bbS^n_{+} \text{ such that } [M^\bbP]^{i,j}_t := [(M^\bbP)^i,(M^\bbP)^j]_t \text{ for } 1 \leq i,j \leq n.
\] We write $\cS\cM(\cW^n)$ as the subset of all semimartingale probability measures and further write for $p \geq 1$, $ \cS\cM_p(\cW^n)$ as the set of all $\bbP \in \cS\cM(\cW^n)$ such that  \[
\bbE^\bbP|V^\bbP|_{1\text{-}\var}^p + \tr([M^\bbP]_1)^{p/2} < \infty.
\] In particular, for $ p = 2$, since $\bbE^{\bbP}\tr([M^{\bbP}]_1)<\infty$, $M^\bbP$ is actually a square-integrable martingale. In the case that $V^\bbP \equiv 0$, i.e. $X$ is a local martingale, we say $\bbP$ is an $\bbR^n$-martingale measure and denote the set of all such martingale measures as $\cM(\cW^n)$ with $\cM_2(\cW^n) := \cM(\cW^n) \cap \cS \cM_2(\cW^n)$.

Fix $\texttt{Prog}_n$ as the $\cF^X$-progressively measurable $\sigma$-algebra on $[0,1] \times \cW^n$, \[
    \scrA^{n,n'} := \{ a : [0,1] \times \cW^n \rightarrow \bbR^{n \times n'} : a \text{ is } \texttt{Prog}_n\text{-measurable} \}
    \] and \[
    \scrM^{n,n'} := \{ a \in \scrA^{n,n'} : a(t,\omega) = h(t,\omega_t) \text{ for some } h : [0,1] \times \bbR^n \rightarrow \bbR^{n \times n'} \}.
    \] Finally, we define $\ocH_k := \{ H \in \scrA^{n,1} : \sup_{t, \omega} |H(t,\omega)| \leq k \}$ for $k \geq 0$.

\subsection{Adapted Wasserstein distances}
\label{subsection:adaptedwassersteindistance}

Let $\bbP, \bbQ \in \cP(\cW^n)$ and $\Pi(\bbP,\bbQ)$ denote the set of all couplings between $\bbP$ and $\bbQ$, i.e. $\pi \in \Pi(\bbP, \bbQ) \iff X_{\#}\pi = \bbP, Y_{\#}\pi = \bbQ$. 

Central to the notion of adapted Wasserstein distances is the following definition of causal couplings introduced in \cite{lassalle2018causaltp,backhoffbeiglbocklinzalashko2017causal}:
\begin{defn}
    \label{defn:causalcplingsandmaps}
    \mbox{}
  \begin{enumerate}
    \item   $\pi \in \Pi(\bbP, \bbQ)$ is said to be causal if \[
    \theta_\pi(B) \text{ is } \prescript{\bbP}{}{\cF}^X_t\text{-measurable for all } t \in [0,1] \text{ and } B \in \cF^X_t,
    \] where $\theta_\pi :\cW^n \rightarrow \cP(\cW^n)$ denotes the disintegration of $\pi$ w.r.t. $\bbP$ and $\prescript{\bbP}{}{\cF}^X_t$ the completion of $\cF^X_t$ w.r.t. $\bbP$. $\pi$ is said to be \emph{bicausal} if $\pi$ and $((\omega, \omega') \mapsto (\omega',\omega))_{\#} \pi$ are both causal. We denote by $\Pi_c(\bbP, \bbQ)$ and $\Pi_{bc}(\bbP, \bbQ)$ the set of all causal and bicausal couplings respectively.
    \item A measurable map $T : \cW^n \rightarrow \cW^n$ is said to be causal or bicausal if the induced coupling $\pi_T := (\Id,T)_{\#}\bbP \in \Pi_c(\bbP, T_{\#}\bbP)$ or $\in \Pi_{bc}(\bbP,T_{\#}\bbP)$ respectively. The set $\cT_{bc}(\bbP,\bbQ)$ denotes the set of all $\pi \in \Pi_{bc}(\bbP,\bbQ)$ such that $\pi = \pi_T$ for some $T$. Abusing notation, we write $T \in \cT_{bc}(\bbP,\bbQ)$ if $\pi_T \in \cT_{bc}(\bbP,\bbQ)$.
    \item Given $T \in \cT_{bc}(\bbP,\bbQ)$, we further say that $T$ is an a.s. bijective bicausal map if there exists $S : \cW^n \rightarrow \cW^n$ such that \[
    S \circ T = \identity \quad \bbP\text{-a.s. and } T \circ S = \identity \quad 
    \bbQ\text{-a.s..}
    \] The set of all a.s. bijective bicausal maps between $\bbP$ and $\bbQ$ is denoted $\cT_{bcb}(\bbP,\bbQ)$.
  \end{enumerate}
\end{defn}

Given a metric $d$ on $\cW^n$, not necessarily the uniform one $|\cdot|_\infty$, the Wasserstein distance and its adapted variant are defined by \[
\cW_{d,p}^p(\bbP, \bbQ) := \inf_{\pi \in \Pi(\bbP,\bbQ)} \bbE^\pi d^p \text{ and } \cA\cW_{d,p}^p(\bbP, \bbQ) := \inf_{\pi \in \Pi_{bc}(\bbP,\bbQ)} \bbE^\pi d^p
\] for $p \geq 1$. In the sequel, we will focus on $p = 2$, in which case we simply write $\cW_{d}(\bbP, \bbQ) = \cW_{d,2}(\bbP, \bbQ)$ and $\cA \cW_d(\bbP,\bbQ) = \cA \cW_{d,2}(\bbP,\bbQ)$. 

\begin{defn}
    \label{def:AQTCIdef}
    A probability measure $\bbQ \in \cP(\cW^n)$ is said to satisfy an adapted quadratic transport cost inequality (AQTCI) w.r.t. distance $d$ with constant $C > 0$ if for all $\bbP \in \cP(\cW^n)$, we have \[
    \cA\cW_d^2(\bbP,\bbQ) \leq C \cH(\bbP|\bbQ),
    \] where $\cH(\bbP|\bbQ)$ denotes the relative entropy of $\bbP$ w.r.t. $\bbQ$: \[
    \cH(\bbP|\bbQ) := \begin{cases}
    \bbE^\bbP \log \frac{d\bbP}{d\bbQ} &\text{ if } \bbP \ll \bbQ, \\
    +\infty &\text{ otherwise. } 
\end{cases}
\]
\end{defn}

\section{Adapted QTCI for diffusions}
\label{section:AQTCIfordiffusions}

In this section, we take advantage of the fact that solutions to SDEs are semimartingales and demonstrate that the AQTCI holds with respect to a semimartingale ``distance" considered in \cite{backhoff2020adaptedwd}.  Let $\bbP_i \in \cS\cM_p(\cW^n), i = 1,2$ for $p \geq 1$ and $\pi \in \Pi_{bc}(\bbP_1,\bbP_2)$. Suppose the semimartingale decomposition of canonical process $X$ under $\bbP_i$ is given by $X = M^{\bbP_i} + V^{\bbP_i}$. Then, using \cite[Remark $2.3$]{acciaio2020causalot} or \cite[Theorem $2.2$]{cont2024causaltransportpathspace}, bicausality of $\pi$ implies that under $\pi$, the marginal processes $X$ and $Y$ are $\cF^{X,Y}$-semimartingales with semimartingale decompositions: \[
X = M^{\bbP_1}(X) + V^{\bbP_1}(X) \text{ and } Y = M^{\bbP_2}(Y) + V^{\bbP_2}(Y).
\] From now on, we simply write $X = M^{\bbP_1} + V^{\bbP_1}$ and $Y = M^{\bbP_2} + V^{\bbP_2}$ and define for $p \geq 1$ \[
\cA\cW^p_{sm,p}(\bbP_1,\bbP_2) := \inf\{\Phi_{1,p}(\pi) : \pi \in \Pi_{bc}(\bbP_1,\bbP_2)\},
\] where, for $t \in [0,1]$, $\Phi_{t,p}(\pi) := \bbE^\pi |(V^{\bbP_1})_{\cdot \wedge t} - (V^{\bbP_2})_{\cdot \wedge t}|_{1\text{-}\var}^p + \tr([M^{\bbP_1} - M^{\bbP_2}]_t)^{p/2}$. In the case $p = 2$, we simply write $\Phi_t := \Phi_{t,2}$ and $\cA\cW_{sm} = \cA \cW_{sm,2}$. Note that if $\bbP_i \in \cS\cM_2(\cW^n), i = 1,2$, then by the BDG inequality, \begin{multline}
    \label{eqn:upperboundinginftynorm}
    \bbE^\pi |X - Y|_\infty^2 \leq 2 \bbE^\pi (|V^{\bbP_1} - V^{\bbP_2}|^2_\infty + |M^{\bbP_1} - M^{\bbP_2}|^2_\infty) \\
    \leq 8\bbE^\pi (|V^{\bbP_1} - V^{\bbP_2}|_{1\text{-}\var}^2 + \tr([M^{\bbP_1} - M^{\bbP_2}]_1)) = 8 \Phi_1(\pi),
\end{multline} from which we see that \[
\cA\cW_{|\cdot|_\infty}(\bbP_1,\bbP_2) \leq \sqrt{8} \cA \cW_{sm}(\bbP_1,\bbP_2).
\] If $\bbP_i \notin \cS\cM_2(\cW^n)$ for $i = 1$ or $2$, then we simply set $\cA\cW_{sm}(\bbP_1,\bbP_2) = \infty$.

In order to prove the semimartingale AQTCI for laws of diffusions, we first elaborate on the bicausality of couplings between diffusions under a finite relative entropy condition.

\subsection{Bicausality under finite relative entropy}

Suppose that $b \in \scrA^{n,1}$ and $\sigma \in \scrA^{n,n}$ are such that the SDE \begin{equation}
    \label{eqn:sde}
    dX_t = b(t,X) dt + \sigma(t,X) dB_t, \quad X_0 = 0,
\end{equation} has a unique (in law) weak solution, i.e. there exists a unique $\mu^{\sigma,b} \in \cS\cM(\cW^n)$ such that under $\mu^{\sigma,b}$, the canonical process $X$ is an $\cF^X$-semimartingale with decomposition $X = V^{\mu^{\sigma,b}} + M^{\mu^{\sigma,b}}$, where $V^{ \mu^{\sigma,b}}$ (resp. $M^{\mu^{\sigma,b}}$) is a finite-variation process (resp. continuous local martingale) satisfying \[
V^{\mu^{\sigma,b}}_t = \int_0^t b(s,X) ds \text{ and } [M^{\mu^{\sigma,b}}]_t = \int_0^t \sigma(s,X) \sigma^*(s,X) ds.
\] It is shown in \cite{ustunel2019mrp} that weak uniqueness is sufficient for $X$ to satisfy a martingale representation property (MRP), which was key in characterizing bicausal couplings between laws of SDEs in \cite{cont2024causaltransportpathspace}. As an extension, we demonstrate in Appendix \ref{appendix:mrpfiniterelent} that for \emph{any} $\bbP \in \cP(\cW^n)$ such that $\cH(\bbP | \mu^{\sigma,b}) < \infty$, it holds that $\bbP \in \cS\cM(\cW^n)$ and under $\bbP$, $X$ will also satisfy an MRP. This leads to analogous characterizations of bicausality between laws of semimartingales with finite relative entropy w.r.t. laws of SDEs with weak uniqueness (cf. \cite[Theorem $3.2$]{cont2024causaltransportpathspace}).

\begin{thm}
    \label{thm:bicausalcplsbetweensemimartfiniterelent}
    Let $b_i \in \scrA^{n,1}$ and $\sigma_i \in \scrA^{n,n}$ be such that the corresponding SDE has a unique weak solution, denoted $\mu^{\sigma_i,b_i}$, and $\bbP_i \in \cP(\cW^n)$ satisfy $\cH(\bbP_i | \mu^{\sigma_i,b_i}) < \infty$ for $i =1,2$. Then, for any $\pi \in \Pi(\bbP_1, \bbP_2)$, $\pi $ is bicausal if and only if there exists an abstract probability space $(\Omega, \cG, \tilde{\bbP})$ and processes $Z^i = (Z^i_t)_{t \in [0,1]}, i = 1,2$ such that \begin{enumerate}
        \item $\pi$ is the joint law of $(Z^1, Z^2)$: $\pi = (Z^1,Z^2)_{\#}\tilde{\bbP},$ 
        \item under $\tilde{\bbP}$, $Z^i$ is an $\cF^{Z^1,Z^2}$-semimartingale whose decomposition w.r.t. $\cF^{Z^1,Z^2}$ agrees with its decomposition w.r.t. $\cF^{Z^i}$ for $i = 1,2$.
    \end{enumerate} 
\end{thm}

\begin{rmk}
    \label{rmk:semimartdecompsame}
    In $2.$, it is essential that the semimartingale decomposition of $Z^i$ w.r.t. the larger filtration $\cF^{Z^1, Z^2}$ remains the same as the decomposition w.r.t. $\cF^{Z^i}$. As an example, take $Z^1 = B$ to be a Brownian motion which by \cite[Chapter $VI$, Theorem $3$]{protter2005stochasticint} remains a semimartingale after the filtration $\cF^B$ is enlarged with the terminal value. In particular, $Z^1$ is a $(\cG_t)_{t\in [0,1]}$-semimartingale for $\cG_t := \sigma(\cF^B_t, B_1)$ with decomposition \[
    B_t = \beta_t + \int_0^t \frac{B_1 - B_s}{1 - s} ds,
    \] where $\beta = (\beta_t)_{t \in [0,1]}$ is a $\cG$-Brownian motion. Furthermore, set $Z^2 =(Z^2_t)_{t \in [0,1]}, Z^2_t := B_t - 2tB_1$. By checking the covariance, one can easily conclude that $Z^2$ is $\cF^{Z^2}$-Brownian motion. We compute \[
    \bbE(Z^2_t - Z^2_s | \cG_s) = \bbE (B_t - B_s|\cG_s) - 2(t - s)B_1 = \frac{t - s}{1 - s}(B_1 - B_s) - 2(t - s) B_1,
    \] from which we conclude $Z^2$ has decomposition (w.r.t. $\cG_t$) \[
    Z^2_t = \beta_t + \int_0^t \frac{B_1 - B_s}{1 - s} - 2B_1 ds.
    \] Finally, noting that $\cF^{Z^1,Z^2} = \cG$, we obtain that even though $Z^i$ are $\cF^{Z^1,Z^2}$-semimartingales, their joint law $\Law(Z^1,Z^2) \in \Pi(\bbW^1,\bbW^1)$ cannot be bicausal since they do not remain Brownian motions under their joint filtration by \cite[Theorem $3.1$]{cont2024causaltransportpathspace}.
\end{rmk}

\begin{proof}
    The proof follows that of \cite[Theorem $3.2$]{cont2024causaltransportpathspace} very closely and we only provide a sketch.

    To see necessity, we fix $\pi \in \Pi_{bc}(\bbP_1, \bbP_2)$ and note that by the equivalence between (bi-)causality and the so-called H-hypothesis (see \cite[Remark $2.3$]{acciaio2020causalot} or \cite[Theorem $2.2$]{cont2024causaltransportpathspace}), we can choose $(\Omega, \cG, \tilde{\bbP}) = (\cW^n \times \cW^n, \cF^X_1 \otimes \cF_1^X, \pi)$ and $Z^1_t(\omega, \omega') := \omega_t, Z^2_t(\omega, \omega') := \omega'_t$. 

    To see sufficiency, $1.$ and $2.$ above imply that the marginal processes $X$ and $Y$ on $\cW^n \times \cW^n$ are $\cF^{X,Y}$-semimartingales under $\pi$. Furthermore, since the martingale part of $X$ w.r.t. $\cF^X$ and $\cF^{X,Y}$ remains the same by $2.$, the stochastic integral representation of square-integrable $\cF^X$-martingales from Proposition \ref{prop:mrpunderfiniterelent} implies that every square-integrable $\cF^X$-martingale remains a martingale in the joint filtration. The equivalence between causality and H-hypothesis then gives $\pi \in \Pi_c(\bbP_1, \bbP_2)$. A completely symmetric argument gives $\pi \in \Pi_{bc}(\bbP_1, \bbP_2)$.
\end{proof}

\subsection{Adapted QTCI}

With the preceding theorem, we can now prove the main result of this article.

\begin{thm}
    \label{thm:adaptedQTCIsmcost}
    Let $b \in \scrA^{n,1}, \sigma \in \scrA^{n,n}$ be such that $\int_0^1 |b(t,0)| dt < \infty$, $\Vert \sigma \Vert_{\op,\infty} := \sup_{t,\omega}  |\sigma(t,\omega)|_{\op} \in (0,\infty)$ and there exists some $K > 0$ satisfying \begin{equation}
        \label{eqn:lipschitzsdecoeff}
         \sup_{t \in [0,1]}\max\{|b(t,\omega) - b(t, \omega')|, |\sigma(t,\omega) - \sigma(t,\omega')|\} \leq K |\omega - \omega'|_\infty \text{ for all } \omega, \omega' \in \cW^n.
    \end{equation} Then, $\mu^{\sigma,b} \in \cS\cM_2(\cW^n)$ and $\mu^{\sigma,b}$ satisfies the semimartingale AQTCI with constant $4\Vert\sigma \Vert_{\op,\infty}^2 e^{32K^2}$: for all $\bbP \in \cP(\cW^n)$, we have \begin{equation}
        \label{eqn:adaptedQTCIsmcost}
        \cA\cW_{sm}^2(\bbP, \mu^{\sigma,b}) \leq 4\Vert\sigma \Vert_{\op,\infty}^2 e^{32K^2} \cH(\bbP|\mu^{\sigma,b}).
    \end{equation} Furthermore, whenever $\cH(\bbP |\mu^{\sigma,b}) < \infty$, it holds that $\bbP \in \cS\cM_2(\cW^n).$
\end{thm}

\begin{rmk}
    \label{rmk:unambiguouslawsde}
    By the assumption on $b,\sigma$, we have the existence of a pathwise unique, strong solution and in particular, a unique weak solution. 
\end{rmk}

\begin{rmk}
    \label{rmk:sigmabounded}
    We note that Lipschitz $\sigma$ alone is insufficient. For example, consider the geometric Brownian motion \[
    X_t = \exp(B_t) - 1, \quad t \in [0,1],
    \] which is the solution to the SDE \[
    dX_t = \frac{1}{2} (X_t + 1) dt + (X_t + 1) dB_t, \quad X_0 = 0.
    \] If $\Law(X)$ satisfies \eqref{eqn:adaptedQTCIsmcost}, then it would also satisfy the classical $\cW_1$-transport inequality: \begin{multline*}
            \cW_{|\cdot|_\infty,1}(\bbP, \Law(X)) \leq \cW_{|\cdot|_\infty,2}(\bbP, \Law(X)) \leq \cA\cW_{|\cdot|_\infty,2}(\bbP, \Law(X)) \\
            \leq 2\sqrt{2}\cA\cW_{sm}(\bbP, \Law(X)) \leq \sqrt{8 C \cH(\bbP| \Law(X))} \text{ for any } \bbP \in \cP(\cW^1).
    \end{multline*} By \cite[Theorem $2.3$]{djelloutguillinwu2004transportationci}, this implies that $\Law(X)$ has finite quadratic exponential moments, but this is not possible since for any $\alpha > 0$, \[
     \infty = \exp(-2\alpha) \bbE \exp(\alpha \exp(2B_1)) \leq \bbE \exp(2\alpha |X|^2_\infty).
    \]
\end{rmk}

\begin{proof}

Standard Gronwall arguments and the Lipschitz assumptions on $\sigma,b$ yield that $\mu^{\sigma,b} \in \cS\cM_2(\cW^n)$, i.e. $\bbE^{\mu^{\sigma,b}}(|V^{\mu^{\sigma,b}}|^2_{1\text{-}\var}+ \tr([M^{\mu^{\sigma,b}}]_1)) < \infty$. 

We first show that if $\cH(\bbP | \mu^{\sigma,b}) < \infty$, then $\bbP \in \cS\cM_2(\cW^n)$. Girsanov's theorem established in \cite{leonard2012girsanov} and recalled in Lemma \ref{thm:cleonardgirsanov} gives the existence of an $\bbR^n$-valued $\cF^X$-adapted process $v : [0,1] \times \cW^n \rightarrow \bbR^n$ defined $\bbP$-a.s. such that \[
\bbE^\bbP \int_0^1 | \sigma^*(s, X) v(s,X)|^2 ds = 2 \cH(\bbP | \mu^{\sigma,b})
\] and under $\bbP$, $X$ is an $\cF^X$-semimartingale with decomposition \begin{equation}
    \label{eqn:AQTCIsmcostsmdecomp}
    X_t = V^\bbP_t + M^\bbP_t,
\end{equation} where $V^\bbP$ is a finite-variation process given by \[
V^\bbP_t := \int_0^t b(s,X) + \sigma \sigma^*(s,X) v(s,X) ds
\] and $M^\bbP$ is a continuous local martingale with quadratic variation $[M^\bbP]_t = \int_0^t \sigma\sigma^*(s,X) ds.$ To verify that $\bbP = \Law(X) \in \cS\cM_2(\cW^n)$, that is, \[
\Lambda := \bbE^{\bbP} (|V^\bbP|^2_{1\text{-}\var} + \tr([M^\bbP]_1)) < \infty,
\] it is sufficient to show, using the boundedness of $\sigma$ and linear growth of $b$, that $\bbE^{\bbP} |X|_\infty^2 < \infty$. By the Gronwall's, H\"older's and BDG inequalities, one can deduce that \[
\bbE^\bbP |X|_\infty^2 \leq 3 e^{2K} ((\int_0^1 |b(s,0)| ds)^2 + 4n \Vert \sigma \Vert^2_{\op,\infty} (\cH(\bbP |\mu^{\sigma,b}) + 1)).
\] It is clear the RHS is finite since $\cH(\bbP|\mu^{\sigma,b}) <\infty$ by assumption.

To prove \eqref{eqn:adaptedQTCIsmcost}, it suffices to consider the case that  $\cH(\bbP | \mu^{\sigma,b}) < \infty$. In particular, by \cite[Theorem $7.1'$]{ikedawatanabe1989stochasticde}, one can extend the filtered probability space $(\Omega, \cF^X_1, \cF^X, \bbP)$ to $(\Omega, \cG_1, \cG = (\cG_t)_{t \in [0,1]}, \tilde{\bbP})$ such that under $\tilde{\bbP}$, \begin{enumerate}
    \item $X$ is a $\cG$-semimartingale with the same decomposition as in \eqref{eqn:AQTCIsmcostsmdecomp}, i.e. $M^\bbP$ is a $\cG$-martingale, and
    \item there exists a $\cG$-Brownian motion $B$ such that \[
    M^\bbP = \int_0^\cdot \sigma(s,X) dB_s, \quad \tilde{\bbP}\text{-a.s..} 
    \]
\end{enumerate} Under the assumption of \eqref{eqn:lipschitzsdecoeff} on $b,\sigma$, there exists a pathwise unique strong solution to \eqref{eqn:sde}, denoted $Z$, driven by $B$. In particular, $Z$ is a $\cG$-semimartingale with decomposition \[
Z_t = \int_0^t b(s,Z) ds + \int_0^t \sigma(s,Z) dB_s. 
\] Choosing $\bbP_1 = \bbP$ and $\bbP_2 = \mu^{\sigma,b}$ in Theorem \ref{thm:bicausalcplsbetweensemimartfiniterelent} yields that $\Law(X,Z) \in \Pi_{bc}(\bbP, \mu^{\sigma,b})$. Hence, to show \eqref{eqn:adaptedQTCIsmcost}, we will show that \begin{multline*}
    \Phi_1(\Law(X,Z)) = \bbE^{\tilde{\bbP}} ((\int_0^1 |b(s,X) + \sigma\sigma^*(s,X)v(s,X) - b(s,Z)|ds)^2 + \int_0^1 |\sigma(s,X) - \sigma(s,Z)|^2 ds) \\
    \leq 2 \Vert \sigma \Vert^2_{\op,\infty} e^{32K^2} \bbE^\bbP \int_0^1 | \sigma^*(s, X) v(s,X)|^2 ds = 4\Vert \sigma \Vert^2_{\op,\infty} e^{32K^2} \cH(\bbP| \mu^{\sigma,b}).
\end{multline*}

Using It\^o's isometry, H\"older's inequality, the Lipschitz continuity of $\sigma,b$ and Lemma \ref{thm:cleonardgirsanov}, we estimate \begin{align*}
     \Phi_t(\Law(X,Z)) = &\bbE^{\tilde{\bbP}} ((\int_0^t |b(s,X) + \sigma\sigma^*(s,X)v(s,X) - b(s,Z)|ds)^2 + \int_0^t|\sigma(s,X) - \sigma(s, Z)|^2 ds)\\
     \leq &2 (\int_0^t \bbE^{\tilde{\bbP}} |b(s,X) - b(s,Z)|^2 + |\sigma(s,X) - \sigma(s,Z)|^2 ds + \Vert \sigma\Vert_{\op, \infty}^2 \bbE^{\tilde{\bbP}} \int_0^1 |\sigma^*(s,X)v(s,X)|^2 ds) \\
     \leq &4K^2 \int_0^t \bbE^{\tilde{\bbP}} |X_{\cdot \wedge s} - Z_{\cdot \wedge s}|^2_\infty ds + 4\Vert \sigma \Vert_{\op, \infty}^2 \cH(\bbP|\mu^{\sigma,b}).
 \end{align*} Furthermore, a calculation analogous to that in \eqref{eqn:upperboundinginftynorm} gives \[
 \bbE^{\tilde{\bbP}} |X_{\cdot \wedge s} - Z_{\cdot \wedge s}|^2_\infty \leq 8 \Phi_s(\Law(X,Z))
 \] so that \[
 \Phi_t(\Law(X,Z)) \leq 32K^2 \int_0^t \Phi_s(\Law(X,Z)) ds + 4\Vert \sigma \Vert_{\op, \infty}^2 \cH(\bbP|\mu^{\sigma,b})
 \] and Gronwall's inequality yields \[
  \Phi_1(\Law(X,Z)) \leq 4\Vert \sigma \Vert_{\op, \infty}^2 e^{32K^2} \cH(\bbP|\mu^{\sigma,b}) \implies \cA\cW_{sm}^2(\bbP, \mu^{\sigma,b}) \leq 4\Vert\sigma \Vert_{\op, \infty}^2 e^{32K^2} \cH(\bbP|\mu^{\sigma,b})
 \] as desired. 
\end{proof}

\begin{rmk}
    \label{rmk:uniformcostAQTCI}
    In the proof of Theorem \ref{thm:adaptedQTCIsmcost}, the use of Girsanov's theorem to construct the coupling $\Law(X,Z) \in \Pi(\bbP, \mu^{\sigma,b})$ was used in previous works (see \cite{ustunel2012transportationci,pal2012concentration,bahlalimouchtabihtangpi2021quadraticti}) to establish classical quadratic transport inequalities w.r.t. $d = | \cdot |_\infty$ or $\sqrt{\int_0^1 | \cdot |^2 ds}$ for the law of diffusions with various conditions on the drift and diffusion coefficients. For example, in \cite[Proposition $2$]{bahlalimouchtabihtangpi2021quadraticti}, the drift and diffusion coefficients are taken to be Sobolev differentiable instead of the usual Lipschitz assumptions. Since Theorem \ref{thm:bicausalcplsbetweensemimartfiniterelent} gives that $\Law(X,Z)$ is actually bicausal, the classical transport inequalities established in the aforementioned works can be improved to their adapted variants without much effort.
\end{rmk}

As a first corollary, we obtain a classical transport inequality for the quadratic variation of $\bbQ$:

\begin{cor}
\label{cor:QTCIforQV}
    Suppose $\bbQ \in \cS\cM_2(\cW^n)$ satisfies the semimartingale AQTCI with constant $C > 0$ and define $S : \cW^n \rightarrow \bbR_{\geq 0}$, $S(X) := \sqrt{\tr[X]_1}$. Then, $S_{\#}\bbQ \in \cP(\bbR)$ satisfies a classical QTCI w.r.t. $d = |\cdot|$ with constant $C$: for all $\mu \in \cP(\bbR)$, \[
    \cW_{|\cdot|}^2(\mu, S_{\#}\bbQ) \leq C \cH(\mu|S_{\#}\bbQ).
    \] In particular, for all $\epsilon > 0$ sufficiently small and $1 \leq p < \infty$, \[
    \bbE^\bbQ \exp(\epsilon \tr[X]_1 ) + \Vert \sqrt{\tr[X]_1} \Vert_{L^p(\bbQ)} = \bbE^\bbQ \exp(\epsilon \tr[M^\bbQ]_1 )  + \Vert \sqrt{\tr[M^\bbQ]_1} \Vert_{L^p(\bbQ)} < \infty.
    \]
\end{cor}

\begin{proof}
    We first note that for any $\bbP \ll \bbQ$ such that $\cH(\bbP|\bbQ) < \infty$, we have $\bbP \in \cS\cM(\cW^n)$ and for any $\pi \in \Pi_{bc}(\bbP,\bbQ)$, \[
    \sqrt{\tr[M^\bbQ]_1} = \sqrt{\tr[M^\bbQ - M^\bbP + M^\bbP]_1} \leq \sqrt{\tr[M^\bbP]_1} + \sqrt{\tr[M^\bbP - M^\bbQ]_1} \quad \pi\text{-a.s.}
    \] so that symmetrizing gives \[
    |\sqrt{\tr[M^\bbQ]_1} - \sqrt{\tr[M^\bbP]_1}|^2 \leq \tr[M^\bbP - M^\bbQ]_1
    \] and we obtain for any $\pi \in \Pi_{bc}(\bbP,\bbQ)$, \[
    \cW_{|\cdot|,2}^2(S_{\#}\bbP, S_{\#}\bbQ) \leq \bbE^\pi |\sqrt{\tr[M^\bbQ]_1} - \sqrt{\tr[M^\bbP]_1}|^2 \leq \Phi_1(\pi).
    \] Taking infimum over $\pi$ and using the AQTCI gives \[
    \cW_{|\cdot|,2}(S_{\#}\bbP, S_{\#}\bbQ) \leq \cA \cW_{sm}(\bbP,\bbQ) \leq  \sqrt{C \cH(\bbP | \bbQ)}.
    \] Finally, since $(2.1)$ in \cite[Lemma $2.1$]{djelloutguillinwu2004transportationci} gives \[
    \cH(\mu| S_{\#}\bbQ) = \inf \{ \cH(\bbP| \bbQ) : S_{\#}\bbP = \mu \} \text{ for any } \mu \in \cP(\bbR),
    \] we obtain the desired transport inequality. Furthermore, applying \cite[Proposition $6.1$]{gozlanleonard2010transport} with $\alpha(t) := \frac{t}{C}$ and $\theta(t) := t^2$ gives for all $\epsilon > 0$ sufficiently small that $\bbE^\bbQ \exp(\epsilon\tr[X]_1 ) < \infty$. In turn, since $\sup_{x \geq 0} x^p e^{- \epsilon x^2} < \infty $, it is straightforward to conclude $\Vert \sqrt{\tr[M^\bbQ]_1} \Vert_{L^p(\bbQ)} < \infty$.
\end{proof}





\subsubsection*{Tensorization}

It is interesting to note that the constant in \eqref{eqn:adaptedQTCIsmcost} seemingly only depends on the Lipschitz constant $K$ and the bound on $\sigma$, $\Vert \sigma \Vert_{\op,\infty}$, and not on the dimension $n$. This suggests that some sort of tensorization property is possible, which is the content of the next corollary:

\begin{cor}
    \label{cor:adaptedQTCItensorization}
    Let $b \in \scrM^{n,1}, \sigma \in \scrM^{n,n}$ be such that \[
b(t,\omega) = f(t,\omega_t) \text{ and } \sigma(t, \omega) = g(t,\omega_t)
\] for some $f : [0,1] \times \bbR^n \rightarrow \bbR^n$ satisfying $\int_0^1 |f(t,0)| dt < \infty$ and non-zero bounded $g : [0,1] \times \bbR^n \rightarrow \bbR^{n \times n}$ such that \[
\sup_{t \in [0,1]} \max\{|f(t,x) - f(t,y)|, |g(t,x) - g(t,y)|\} \leq K | x- y| \text{ for all } x,y \in \bbR^n
\] for some $K > 0$. Then, for any $N \geq 1$ and $\bbP \in \cP(\cW^{n N})$, \[
\cA\cW_{sm}^2(\bbP, (\mu^{\sigma,b})^{\otimes N}) \leq 4 \Vert \sigma \Vert_{\op, \infty}^2 e^{32 K^2} \cH(\bbP|(\mu^{\sigma,b})^{\otimes N}).
\]
\end{cor}

\begin{rmk}
    \label{rmk:ctsfnspaceprdt}
    Here, $(\mu^{\sigma,b})^{\otimes N}$ is a probability measure on $(\cW^n)^N$, but since $(\cW^n)^N \cong \cW^{nN}$, we identify $(\mu^{\sigma,b})^{\otimes N}$ as a probability measure on $\cW^{nN}$. 
\end{rmk}

\begin{proof}
    We prove the case $N=2$, noting that general $N$ follows very similarly. We define the drift and diffusion coefficients $\tb \in \cM^{2n,1},\tsigma \in \cM^{2n, 2n}$ as \[
    \tb(t,\omega) := \begin{pmatrix}
        f(t, \omega^{1:n}_t) \\
        f(t, \omega^{n+1:2n}_t) 
    \end{pmatrix} \text{ and } \tsigma(t,\omega) := \begin{pmatrix}
        g(t, \omega^{1:n}_t) & 0 \\
        0 & g(t, \omega^{n+1:2n}_t) 
    \end{pmatrix},
    \]  where $\omega^{1:n}_t$ and $\omega^{n+1:2n}_t$ denote the first and last $n$ coordinates of $\omega_t \in \bbR^{2n}$. It is clear that $\Vert \tsigma \Vert_{\op,\infty} = \Vert \sigma \Vert_{\op, \infty}$. Furthermore, one can easily verify that \eqref{eqn:lipschitzsdecoeff} holds for $\tb, \tsigma$ for the \emph{same} Lipschitz constant $K$. Hence, noting that $\mu^{\tsigma,\tb} = (\mu^{\sigma,b})^{\otimes 2} $, Theorem \ref{thm:adaptedQTCIsmcost} yields the desired inequality.
\end{proof}

An immediate consequence of Corollary \ref{cor:adaptedQTCItensorization} is a quantitative propagation of chaos result 
    under the adapted Wasserstein metric for mean field diffusions as studied in \cite{lacker2023hierarchies}. Indeed, fix $b_0 \in \scrA^{n,1}$ and progressively measurable $b : [0,1] \times \cW^n \times \cW^n \rightarrow \bbR^n$ satisfying \[
    |b_0(t,\omega) - b_0(t, \omega')| + |b(t, \omega, \overline{\omega}) - b(t, \omega', \overline{\omega}')| \leq K (|\omega - \omega'|_\infty+ |\overline{\omega} - \overline{\omega}'|_\infty)
    \] for all $t \in [0,1], \omega, \omega', \overline{\omega}, \overline{\omega}' \in \cW^n$ and \[
    \int_0^1 |b_0(t,0) + b(t,0,0)| dt < \infty.
    \] Consider the SDE system \begin{equation}
    \label{eqn:interactingparticle}
        dX^i_t = (b_0(t,X^i) + \frac{1}{N-1} \sum_{j \neq i} b(t,X^i,X^j)) dt + dW^i_t, \quad i = 1,...,N,
    \end{equation} where $W^1,\dots, W^N$ are independent $n$-dimensional Brownian motions and the associated McKean-Vlasov equation \begin{equation}
        \label{eqn:mckeanvlasov}
        dX_t = (b_0(t,X) + \int_{\cW^n} b(t,X,x) \bbQ(dx)) dt + dW_t, \quad \bbQ = \Law(X),
    \end{equation} where $W$ is an $n$-dimensional Brownian motion. Denoting the laws of \eqref{eqn:interactingparticle} and \eqref{eqn:mckeanvlasov} as $\bbP^{(N)} \in \cP(\cW^{nN})$ and $\bbQ$ respectively, \cite{lacker2023hierarchies} establishes a quantitative local propagation of chaos for the relative entropy: there exist constants $C,\gamma > 0$ such that for all $N$ sufficiently large and $k \in [1,N]$,  \[
    \cH(\bbP^{(N,k)}| \bbQ^{\otimes k}) \leq 2 C \frac{k^2}{N^2} + C \exp(-2N(e^{-\gamma} - \frac{k}{N})^2_+), 
    \] where $\bbP^{(N,k)}$ is the marginal law of the first $k$ coordinates of $\bbP^{(N)}$. Using a classical transport inequality that tensorizes, the above entropy bound implies propagation in Wasserstein distance \cite[Corollary $2.7$]{lacker2023hierarchies}, i.e. there exist constants $C',\gamma' > 0$ such that \[
    \cW_{|\cdot|_\infty,2}^2(\bbP^{(N,k)}, \bbQ^{\otimes k}) \leq 2 C' \frac{k^2}{N^2} + C' \exp(-2N(e^{-\gamma'} - \frac{k}{N})^2_+).
    \] Suppose now $b_0(t,\omega) = f_0(t,\omega_t)$ and $b(t,\omega, \overline{\omega}) = f(t,\omega_t,\overline{\omega}_t)$ for some Lipschitz $f_0 : [0,1] \times \bbR^n \rightarrow \bbR^n, f : [0,1] \times \bbR^n \times \bbR^n \rightarrow \bbR^n$, uniformly in $t \in [0,1]$. By replacing the use of the classical transport inequality in the final step of the proof of \cite[Corollary $2.7$]{lacker2023hierarchies} with Corollary \ref{cor:adaptedQTCItensorization}, we obtain propagation under semimartingale adapted Wasserstein distance, that is, the existence of constants $C'',\gamma'' > 0$ such that \[
    \cA\cW_{sm}^2(\bbP^{(N,k)}, \bbQ^{\otimes k}) \leq 2 C'' \frac{k^2}{N^2} + C'' \exp(-2N(e^{-\gamma''} - \frac{k}{N})^2_+).
    \]
    

\subsection{Stability with respect to bounded drift perturbation}
It is well known that several quantities arising from stochastic analysis/control are well-behaved with respect to the adapted Wasserstein distance. We collect here some examples:
\begin{exam}
    \label{exam:awlipschitzbd}
    Let $p \geq 1$ and $\bbP, \bbQ \in \cP_{d,p}(\cW^n) := \{ \bbP \in \cP(\cW^n) : \bbE^\bbP d(X,0)^p < \infty \}$ for some distance $d$ on $\cW^n$.
    \begin{enumerate}
	    \item(Optimal stopping problems) Let $F : [0,1] \times \cW^n \rightarrow \bbR_{\geq 0}$ be $\texttt{Prog}_n$-measurable and \[
        v^F(\bbP) := \inf_{\tau \; \cF^X\text{-stopping time}} \bbE^\bbP F(\tau,X).
        \] Then, \cite[Theorem $8.1$]{hitzrobinson2026bicausalot} shows that if $F$ satisfies \[
        |F(t,\omega) - F(t, \omega')| \leq K d(\omega, \omega'),
        \] where $d(\omega,\omega') := (\int_0^1 |\omega_s - \omega'_s|^p ds)^{1/p}$ and $K > 0$ is some constant independent of $t$, we have \[
        |v^F(\bbP) -v^F(\bbQ)| \leq C_F \cA\cW_{d,p}(\bbP,\bbQ),
        \] for some constant $C_F > 0$. 
        \item(Optional projections) Let $p=2,$ $F : \cW^n \rightarrow \cW^n$ and $h : \cM_b([0,1],\bbR^n) \rightarrow \bbR$ both be Lipschitz w.r.t. $|\cdot|_\infty$. Then, interpreting $F$ as a process, i.e. $F = (F_t)_{t \in [0,1]}$ where $F_t(\omega) := F(\omega)(t)$, we have \begin{equation}
            \label{eqn:optprojstab}
            | \bbE^\bbP h((\bbE^\bbP(F_t|\cF^X_t))_{t \in [0,1]}) - \bbE^\bbQ h((\bbE^\bbQ(F_t|\cF^X_t))_{t \in [0,1]}) | \leq C_{F,h} \cA\cW_{|\cdot|_\infty,2}(\bbP, \bbQ).
        \end{equation} for some constant $C_{F,h} > 0$. Here, the somewhat informal $(\bbE^\bbP(F_t|\cF^X_t))_{t \in [0,1]}$ can be made rigorous via the notion of optional projections. We detail this along with the proof of \eqref{eqn:optprojstab} in Appendix \ref{appendix:optprojstab}.
    \end{enumerate}
    Suppose further that $\bbP,\bbQ \in \cS\cM_{p}(\cW^n)$.
    \begin{enumerate}[resume]
        \item(Stochastic integrals)  
        For any $p \geq 2$ and $H \in \ocH_k$ such that $H_t : \cW^n \rightarrow \bbR^n$ is Lipschitz for all $t \in [0,1]$ with Lipschitz constant independent of $t$, one can easily obtain using \cite[Lemma $3.9$]{backhoff2020adaptedwd} the estimate \[
        |\Vert (H \cdot X)_1 \Vert_{L^{p/2}(\bbP)} - \Vert (H \cdot X)_1 \Vert_{L^{p/2}(\bbQ)} | \leq C_{k,p,H} R_{\bbP,\bbQ}(\cA\cW_{sm,p/2}(\bbP,\bbQ) + \cA\cW_{sm,p}(\bbP,\bbQ))
        \] for some constant $C_{k,p,H} > 0$ and $R_{\bbP, \bbQ} := 1 + \cA\cW_{sm,p}(\bbP, \delta_0) + \cA\cW_{sm,p}(\bbQ,\delta_0)$. 
        \item(AVaR{$_\alpha$}) The average value at risk at level $\alpha \in (0,1)$ under $\bbP$ of some random variable $Z : \cW^n \rightarrow \bbR$ is defined as $\text{AVaR}_\alpha^\bbP(Z) := \inf_{m \in \bbR} \bbE^\bbP ((Z-m)_+/\alpha + m)$. Then, for $p=2$, any $F : \cW^n \rightarrow \bbR$ Lipschitz w.r.t. $|\cdot|_\infty$ and $H \in \ocH_k$ as in $3.$, \cite[Theorem $1.6$]{backhoff2020adaptedwd} gives \begin{multline*}
            |\text{AVaR}_\alpha^\bbP(F(X) - (H \cdot X)_1) - \text{AVaR}_\alpha^\bbQ(F(X) - (H \cdot X)_1)| \\
        \leq C_{\alpha,F,k,H} R_{\bbP,\bbQ} (\cA\cW_{sm,1}(\bbP, \bbQ) + \cA\cW_{sm}(\bbP,\bbQ))
        \end{multline*} for some constant $C_{\alpha,F,k,H} > 0$. Similar results with regard to utility maximization or indifference pricing or law-invariant risk measures can also be established (see \cite{backhoff2020adaptedwd}).
        
        \item(Mean-variance portfolio selection)
        For $p = 2$ and $\gamma \geq 0$, the mean-variance portfolio selection problem is defined as \[
        MV(\bbP) := \sup_{H \in \ocH_k} \bbE^\bbP (H\cdot X)_1 - \gamma \text{Var}^\bbP((H\cdot X)_1),
        \] where $\text{Var}^\bbP(Z) := \bbE^\bbP (Z - \bbE^\bbP Z)^2.$ Then, we have \[
        |MV(\bbP) - MV(\bbQ)| \leq C_{k,\gamma}R_{\bbP,\bbQ} \cA\cW_{sm}(\bbP,\bbQ),
        \] for some constant $C_{k,\gamma} > 0$. We note that a similar result was shown in \cite{aghapour2026dynamicdatagenerationdynamic} for the discrete-time case and refer readers to Appendix \ref{appendix:mvpselect} for the proof in our continuous-time setting.
    \end{enumerate}
\end{exam}

As a consequence of the AQTCI for diffusions, the above quantities can be shown to be stable w.r.t. bounded drift perturbations. For example,
\begin{cor}
    \label{cor:stabwrtuniformdrift}
    Suppose $V : \cS\cM_2(\cW^n) \rightarrow \bbR$ satisfies \[
    |V(\bbP) - V(\bbQ)| \leq \tilde{K} \cA\cW_{sm}(\bbP,\bbQ)
    \] for some $\tilde{K} > 0$. Further suppose that $b \in \scrA^{n,1}, \sigma \in \scrA^{n,n}$ are such that the SDE \eqref{eqn:sde} has a unique weak solution, $\mu^{\sigma,b} \in \cS\cM_2(\cW^n)$, satisfying the semimartingale AQTCI with constant $C > 0$. Then, for any bounded $v \in \scrA^{n,1}$, the SDE \begin{equation}
        \label{eqn:perturbedsde}
        dX_t = \tb(t,X) dt + \sigma(t,X) dB_t, X_0 = 0, \quad \tb(t,X) := b(t,X) + \sigma(t,X) v(t,X)
    \end{equation} has a unique weak solution, $\mu^{\sigma,\tb} \in \cS\cM_2(\cW^n)$ and \[
    |V(\mu^{\sigma,b}) - V(\mu^{\sigma,\tb})| \leq \tilde{K}\sqrt{\frac{C}{2}} \Vert v \Vert_\infty,
    \] where $\Vert v \Vert_\infty := \sup_{t,\omega} |v(t,\omega)|$.
\end{cor}
In the case where $V$ is only locally Lipschitz as in Example \ref{exam:awlipschitzbd} $5.$, one can still derive a similar bound in an analogous fashion to the proof below. 

\begin{proof}
    The assertion that \eqref{eqn:perturbedsde} has a unique weak solution follows from the fact that \eqref{eqn:sde} is assumed to have a unique weak solution (see \cite[Chapter $4.4$]{ikedawatanabe1989stochasticde}).
    By the semimartingale AQTCI, the claim that $\mu^{\sigma,\tb} \in \cS \cM_2(\cW^n)$ and the desired inequality follow immediately if we show that \[
    \cH(\mu^{\sigma,\tb} | \mu^{\sigma,b}) = \frac{1}{2} \bbE^{\mu^{\sigma,\tb}} \int_0^1 | \sigma^\dag \sigma(t,X) v(t,X)|^2 dt \leq \frac{1}{2} \Vert v \Vert_\infty^2,
    \] where $\sigma^\dag$ denotes the Moore-Penrose pseudoinverse of matrix $\sigma$. Here, the inequality is a consequence of the fact that $\sigma^\dag \sigma$ is an orthogonal projection and hence, a contraction. To do this, we note that on the filtered probability space $(\cW^n, \cF^X_1, \cF^X, \mu^{\sigma,\tb})$, the canonical process $X$ has decomposition \[
    X_t = \int_0^t b(s,X) + \sigma(s,X)v(s,X) ds + M^{\mu^{\sigma,\tb}}_t,
    \] for some $\cF^X$-continuous local martingale $M^{\mu^{\sigma,\tb}}$, whose quadratic variation is given by $[M^{\mu^{\sigma,\tb}}]_\cdot = \int_0^\cdot \sigma\sigma^*(s,X) ds$. Define $\alpha(s,X) := \sigma^\dag(s,X)^* v(s,X)$ and the stochastic exponential \[
    Z_t := \exp(-\int_0^t \alpha(s,X) dM^{\mu^{\sigma,\tb}}_s - \frac{1}{2} \int_0^t \alpha(s,X)^* \sigma(s,X) \sigma(s,X)^* \alpha(s,X) ds) > 0.
    \] For the rest of the proof, we temporarily omit the dependence on $(s,X)$ to simplify notation. Using the fact that the orthogonal projection $\sigma^\dag \sigma$ is symmetric and $\sigma^* \alpha = \sigma^* (\sigma^\dag)^* v = (\sigma^\dag \sigma)^* v = \sigma^\dag \sigma v$, we conclude that \[
    [-\alpha \cdot M^{\mu^{\sigma, \tb}}]_t = \int_0^t |\sigma^\dag \sigma v|^2 ds \leq \Vert v \Vert_\infty^2 <\infty.
    \] Hence, Novikov's condition yields that $Z := (Z_t)_{t \in [0,1]}$ is a true martingale and $\bbE^{\mu^{\sigma,\tb}} Z_1 = 1$. Define $d\bbP := Z_1 d\mu^{\sigma,\tb}$. By Girsanov's theorem \citep[Theorem $5.22$]{legall2016brownianmotion}, $M^\bbP := M^{\mu^{\sigma,\tb}} + \int_0^\cdot \sigma \sigma^* \alpha ds$ is a $\bbP$-continuous local martingale under $\bbP$. By properties of the pseudoinverse (see \cite[Section $5.5.4$]{golubvanloan1996matrix}), we obtain \[
    \sigma v = \sigma \sigma^\dag \sigma v = \sigma (\sigma^\dag \sigma)^* v = \sigma \sigma^* (\sigma^\dag)^* v = \sigma \sigma^* \alpha 
    \] so that, under $\bbP$, $X$ has decomposition \[
    X_t = \int_0^t b(s,X) ds + M^{\bbP}_t, X_0 = 0
    \] and $[M^\bbP] = [M^{\mu^{\sigma,b}}] = \int_0^\cdot \sigma\sigma^* ds$. In particular, up to extending the probability space, there exists a Brownian motion $B$ such that $M^{\bbP} = \int_0^t \sigma(s,X) dB_s$. By weak uniqueness, we obtain that $\bbP = \Law(X) = \mu^{\sigma,b}$. Hence, using the fact that $Z_1 > 0$, we obtain that $\frac{d\mu^{\sigma,\tb}}{d\mu^{\sigma,b}} = (\frac{d\bbP}{d\mu^{\sigma,\tb}})^{-1} = Z_1^{-1}$ and one can verify in a straightforward manner that \[
    \cH(\mu^{\sigma,\tb}| \mu^{\sigma,b}) = -\bbE^{\mu^{\sigma,\tb}} \log Z_1 = \frac{1}{2}\bbE^{\mu^{\sigma,\tb}} \int_0^1 | \sigma^\dag \sigma(t,X) v(t,X)|^2 dt
    \] as required.
\end{proof}

\section{Concentration}
\label{section:concentration}

We now study concentration phenomena implied by the semimartingale AQTCI. The proofs of our concentration results mirror the classical method very strongly. For example, \eqref{eqn:lipdeviation} follows by first establishing a bound of the form \[
| \bbE^\bbP F - \bbE^\bbQ F | \leq \cW_{d,p}(\bbP,\bbQ)
\] for appropriate $F$, then applying the classical QTCI together with Donsker-Varadhan duality for entropy. In particular, any $1$-Lipschitz $F$ would be sufficient. Here, we replace $F$ by stochastic integrals, which are well controlled by the semimartingale adapted Wasserstein distance, a property already exploited in \cite{backhoff2020adaptedwd}.

\subsection{Stochastic integrals}

\begin{prop}
    \label{prop:tailbehaviorforstointsemimartbddintegrand}
    Suppose $\bbQ \in \cS\cM_2(\cW^n)$ satisfies the semimartingale AQTCI with constant $C > 0$. Then, for any $H \in \scrA^{n,1}$ such that $\Vert H \Vert_\infty := \Vert \Vert H \Vert_{L^\infty(\Leb)} \Vert_{L^\infty(\bbQ)} \in (0,\infty)$, we have for any $\lambda \geq 0$, \begin{equation}
        \label{eqn:laplaceboundforstointbddintegrand}
        \bbE^\bbQ \exp(\lambda ((H \cdot X)_1 - \bbE^\bbQ (H \cdot X)_1)) \leq \exp(2 \lambda \Vert H \Vert_\infty \bbE^\bbQ |V^\bbQ|_{\onevar} + \frac{1}{4} \lambda^2 \Vert H \Vert^2_\infty C).
    \end{equation}  In particular, for any $r \geq 0$, \[
    \bbQ(|(H \cdot X)_1 - \bbE^\bbQ (H \cdot X)_1| \geq r) \leq 2 \exp(- \frac{\max\{r - 2\Vert H \Vert_\infty \bbE^\bbQ | V^\bbQ|_{\onevar},0\}^2}{\Vert H \Vert^2_\infty C}).
    \]
\end{prop}

\begin{rmk}
    \label{rmk:martingalemeastailbehavior}
    Suppose $\bbQ \in \cM_2(\cW^n)$ is a martingale measure satisfying the semimartingale AQTCI. Then, $V^\bbQ = 0$ and $\bbE^\bbQ ( H \cdot X)_1 = 0$, from which we conclude by Doob's maximal inequality and convexity of $x \mapsto \exp(\lambda x)$ that \[
    \bbQ( \sup_{t \in [0,1]} (H \cdot X)_t \geq r) \leq \exp(-\lambda r) \bbE^\bbQ (\exp(\lambda (H \cdot X)_1))
         \leq \exp(-\lambda r + \frac{1}{4} \lambda^2 \Vert H \Vert^2_\infty C).
    \] Optimizing over $\lambda$ and replacing $H$ by $-H$, we obtain that \[
    \bbQ( |(H\cdot X)|_\infty \geq r) \leq 2 \exp(- \frac{r^2}{\Vert H \Vert^2_\infty C}).
    \] For example, if $\bbQ = \bbW^n$, then we recover the sub-Gaussian tail behavior for (Brownian) martingales whose quadratic variation has a uniformly bounded derivative (see \cite[Chapter $IV$, Exercise $3.16$]{revuzyor1999ctsmart}).
\end{rmk}


\begin{proof}
    Since $H$ is bounded, then $\bbE^\bbQ (H \cdot M^\bbQ)_1 = 0$ and $\bbE^\bbQ(H \cdot X)_1 = \bbE^\bbQ ( H \cdot V^\bbQ)_1$. Furthermore, fix any $\bbP \ll \bbQ$ such that $\cH(\bbP | \bbQ) < \infty$. Since the AQTCI, together with $\bbQ \in \cS \cM_2(\cW^n)$, implies $\bbP$ also lies in $\cS\cM_2(\cW^n)$, we similarly have $\bbE^\bbP(H \cdot X)_1 = \bbE^\bbP ( H \cdot V^\bbP)_1$. For any $\pi \in \Pi_{bc}(\bbP, \bbQ)$, we estimate \begin{align*}
        |\bbE^\bbP (H \cdot X)_1 - \bbE^\bbQ (H \cdot X)_1| = &|\bbE^\pi (H(X) \cdot V^\bbP)_1 - \bbE^\bbQ (H(Y) \cdot V^\bbQ)_1| \\
        \leq & \bbE^\pi |((H(X) - H(Y)) \cdot V^\bbQ)_1| + \bbE^\pi | (H(X) \cdot (V^\bbP - V^\bbQ))_1 | \\
        \leq & 2 \Vert H \Vert_\infty \bbE^\bbQ | V^\bbQ|_{\onevar} + \Vert H \Vert_\infty \bbE^\pi |V^\bbP - V^\bbQ|_{\onevar}. 
    \end{align*} Taking infimum over $\pi$ and using the AQTCI, we obtain \begin{equation}
        \label{eqn:stointAQTCIbound}
        |\bbE^\bbP (H \cdot X)_1 - \bbE^\bbQ (H \cdot X)_1| \leq \Vert H \Vert_\infty (2 \bbE^\bbQ | V^\bbQ|_{\onevar} +  \sqrt{C \cH(\bbP| \bbQ)}).
    \end{equation} For any $\lambda > 0 $, by the Donsker-Varadhan duality for entropy, we have \begin{align*}
            \log \bbE^\bbQ \exp(\lambda ((H\cdot X)_1 - \bbE^\bbQ (H\cdot X)_1)) \leq &\sup_{\bbP \ll \bbQ} \lambda(\bbE^\bbP( H \cdot X)_1 - \bbE^\bbQ( H \cdot X)_1) - \cH(\bbP|\bbQ) \\
            \leq &\sup_{\bbP \ll \bbQ} \lambda \Vert H \Vert_\infty (2 \bbE^\bbQ | V^\bbQ|_{\onevar} + \sqrt{C \cH(\bbP| \bbQ)}) - \cH(\bbP |\bbQ) \\
            \leq &\sup_{\zeta \geq 0}  \lambda \Vert H \Vert_\infty (2 \bbE^\bbQ | V^\bbQ|_{\onevar} + \sqrt{C \zeta}) - \zeta \\
            \leq & 2 \lambda \Vert H \Vert_\infty \bbE^\bbQ |V^\bbQ|_{\onevar} + \frac{1}{4} \lambda^2 \Vert H \Vert_\infty^2 C.
    \end{align*} In particular, for any $\lambda, r > 0$, Markov's inequality yields \begin{align*}
        \bbQ((H\cdot X)_1 - \bbE^\bbQ (H \cdot X)_1 \geq r) \leq & \bbQ(\exp(\lambda ((H\cdot X)_1 - \bbE^\bbQ (H \cdot X)_1))\geq \exp (\lambda r))\\
        \leq & \exp(-\lambda r) \bbE^\bbQ \exp(\lambda ((H\cdot X)_1 - \bbE^\bbQ (H\cdot X)_1)) \\
        \leq & \exp(- \lambda r + 2 \lambda \Vert H \Vert_\infty \bbE^\bbQ |V^\bbQ|_{\onevar} + \frac{1}{4} \lambda^2 \Vert H \Vert_\infty^2 C).
    \end{align*} Optimizing over $\lambda > 0$, we obtain \[
    \bbQ((H\cdot X)_1 - \bbE^\bbQ (H \cdot X)_1 \geq r) \leq \exp( - \frac{\max\{r - 2 \Vert H \Vert_\infty \bbE^\bbQ |V^\bbQ|_{\onevar},0\}^2}{\Vert H \Vert_\infty^2 C}).
    \] Up to replacing $H$ by $-H$, we obtain the desired.
\end{proof}

The boundedness assumption on $H$ can be replaced by a uniform H\"older condition in the path variable.

\begin{prop}
    \label{prop:tailbehaviorforstointmartsqintintegrand}
    Suppose $\bbQ \in \cS\cM_2(\cW^n)$ satisfies the semimartingale AQTCI with constant $C > 0$ and $ H \in \scrA^{n,1}$ satisfies \[
    \llbracket H \rrbracket_\alpha := \sup_{t \in [0,1]} |H_t(0)| + \sup_{\omega \neq \omega'} \sup_{t \in [0,1]} \frac{|H_t(\omega) - H_t(\omega')|}{| \omega - \omega'|^\alpha_\infty} \in (0,\infty) \text{ for some } \alpha \in (0,1].
    \] Then, \begin{enumerate}
        \item for $\alpha \in (0,1)$ and  any $\lambda > 0$, we have \begin{equation}
            \label{eqn:laplaceboundforstointl2integrand}
            \bbE^\bbQ \exp(\lambda((H\cdot X)_1 - \bbE^\bbQ( H \cdot X)_1 )) \leq \exp ( \sup_{\zeta > 0 } \lambda g_\alpha (\zeta) - \zeta),
        \end{equation} where $g_\alpha(\zeta) := A_\alpha \zeta^{\alpha/2} + B_\alpha \zeta^{(\alpha + 1)/2} +D \zeta^{1/2}$ for constants \[
         A_\alpha := \llbracket H \rrbracket_\alpha C^{\alpha/2}  \Vert |V^\bbQ|_{\onevar} \Vert_{L^2(\bbQ)} 5^{\alpha/2}, B_\alpha :=\llbracket H \rrbracket_\alpha C^{(\alpha+1)/2} 5^{\alpha/2},D:= \Vert H \Vert_{2, \infty} \sqrt{C}
         \] and $\Vert H \Vert_{2,\infty} := \Vert \Vert H \Vert_{L^\infty(\Leb)} \Vert_{L^2(\bbQ)} $. In particular, \[
         \bbQ(|(H\cdot X)_1 - \bbE^\bbQ (H \cdot X)_1| \geq r) \leq 2\exp(-g^{-1}_\alpha (r)).
         \]
        \item for $\alpha = 1$ and $0 < \lambda < 1/B_1$, we have \begin{equation}
            \label{eqn:laplaceboundsforstointlipint}
            \bbE^\bbQ \exp(\lambda ((H \cdot X)_1 - \bbE^\bbQ ( H \cdot X)_1)) \leq \exp(- \frac{1}{4} \frac{\lambda^2 (A_1 + D)^2}{\lambda B_1 - 1}).
        \end{equation} In particular, \[
        \bbQ(|(H\cdot X)_1 - \bbE^\bbQ (H \cdot X)_1| \geq r) \leq 2\exp(- \frac{4r^2}{(A_1 + D + \sqrt{(A_1 + D)^2 + 4 B_1 r})^2}).
        \]
    \end{enumerate} 
\end{prop}

\begin{rmk}
    \label{rmk:tailbehaviorforholderintegrand}
    For $\alpha \in (0,1)$ and $r > 0$, choosing $\zeta_r := \min \{ (\frac{r}{3 A_\alpha})^{2/\alpha}, (\frac{r}{3 B_\alpha})^{2/(1 + \alpha)}, (\frac{r}{3D})^2 \}$ gives \[
    g_\alpha(\zeta_r) \leq r \implies \zeta_r \leq g_\alpha^{-1}(r)
    \] and hence, \[
    \bbQ(|(H \cdot X)_1 - \bbE^\bbQ ( H \cdot X)_1| \geq r) \leq 2\exp(- \zeta_r).
    \] For $r$ large, we obtain that \[
    \bbQ(|(H \cdot X)_1 - \bbE^\bbQ ( H \cdot X)_1| \geq r) \leq 2 \exp(- (\frac{r}{3 B_\alpha})^{2/(1 + \alpha)}),
    \] which indicates a sub-Weibull tail behavior since $1 < \frac{2}{1 + \alpha} < 2$ (see \cite{vladimirova2020subweibull}). In general, one cannot improve the exponent of $r$ since one can take, for example, $\bbQ = \bbW^1$ and $H_t(\omega) := \mathbbm{1}_{[1/2,1)}(t) |\omega_{1/2}|^\alpha$ so that \[
    (H \cdot X)_1 = |X_{1/2}|^\alpha (X_1 - X_{1/2})
    \] and up to replacing $r$ by $r^{(1 + \alpha)/2}$, we obtain \[
    \bbQ(|X_{1/2}|^\alpha |X_1 - X_{1/2}| \geq r^{(1 + \alpha)/2}) \leq 2 \exp(- (\frac{1}{3 B_\alpha})^{2/(1 + \alpha)} r).
    \] However, by independence of increments and the fact that $\Law(X_{1/2}) = \Law(X_1 - X_{1/2}) = \cN(0,1/2)$, we have \[
    \frac{2}{\pi} \frac{2r}{(1 + 2 r)^2} e^{-2r} \leq \bbQ(| \sqrt{2 }X_{1/2}| \geq \sqrt{2r})^2 \leq \bbQ(|X_{1/2}|^\alpha |X_1 - X_{1/2}| \geq r^{(1 + \alpha)/2}),
    \] where the lower bound follows from integration by parts of the Gaussian density. Hence, \[
    \frac{2}{\pi} \frac{2r^{\frac{2}{1 + \alpha}}}{(1 + 2 r^{2/(1 + \alpha)})^2} e^{-2r^{\frac{2}{1 + \alpha}}} \leq \bbQ(|(H \cdot X)_1| \geq r).
    \] This indicates that the $\frac{2}{1 + \alpha}$ exponent of $r$ governing the tail behavior cannot be improved in general. For $\alpha = 1$, for large $r$, \[
    \frac{4r^2}{(A_1 + D + \sqrt{(A_1 + D)^2 + 4 B_1 r})^2} \approx A' r
    \] for some constant $A'$, indicating that $(H \cdot X)_1$ has sub-exponential tail behavior. The sharpness of this tail behavior is clear if we take $H_t(\omega) := \omega_t$ and $\bbQ = \bbW^1$ since we know that $(H \cdot X)_1 = \frac{1}{2}(X_1^2 - 1)$. $X_1^2$ follows a chi-squared distribution, which is well-known to have sub-exponential tails. 
\end{rmk}

\begin{proof}
    As a preliminary estimate, we have \[
    \sup_{t \in [0,1]} |H_t(\omega)| \leq \sup_{t \in [0,1]} |H_t(\omega) - H_t(0)| + \sup_{t \in [0,1]} |H_t(0)| \leq \llbracket H \rrbracket_\alpha (1 + |\omega|^\alpha_\infty)
    \] so that \[
    \Vert H \Vert_{2,\infty} \leq \llbracket H \rrbracket_\alpha ( 1 + \Vert |X|_\infty \Vert_{L^{2}(\bbQ)}^\alpha) < \infty.
    \] In particular, for any $\bbP \ll \bbQ$ such that $\cH(\bbP|\bbQ) < \infty$, under the AQTCI, we also have $\Vert \Vert H \Vert_{L^\infty(\Leb)} \Vert_{L^2(\bbP)} < \infty$. Hence, as in Proposition \ref{prop:tailbehaviorforstointsemimartbddintegrand}, for any $\pi \in \Pi_{bc}(\bbP, \bbQ)$, we obtain that \begin{align*}
         &|\bbE^\bbP (H \cdot X)_1 - \bbE^\bbQ (H \cdot X)_1| \\
         = &|\bbE^\pi (H(X) \cdot V^\bbP)_1 - \bbE^\bbQ (H(Y) \cdot V^\bbQ)_1| \\
        \leq & \bbE^\pi |((H(X) - H(Y)) \cdot V^\bbQ)_1| + \bbE^\pi | (H(X) \cdot (V^\bbP - V^\bbQ))_1 | \\
        \leq & \bbE^\pi \Vert H(X) - H(Y) \Vert_{L^\infty(\Leb)} |V^\bbQ|_{\onevar} + \bbE^\pi \Vert H(X) \Vert_{L^\infty(\Leb)}|V^\bbP - V^\bbQ|_{\onevar} \\
        \leq &\llbracket H \rrbracket_\alpha \bbE^\pi |X - Y|_\infty^\alpha |V^\bbQ|_{\onevar} + \bbE^\pi (\Vert H(X) - H(Y) \Vert_{L^\infty(\Leb)} + \Vert H \Vert_{L^\infty(\Leb)}) |V^\bbP - V^\bbQ|_{\onevar} \\
        \leq & \llbracket H \rrbracket_\alpha \Vert |X-Y|_\infty \Vert_{L^{2}(\pi)}^\alpha (\Vert |V^\bbQ|_{\onevar} \Vert_{L^2(\bbQ)} + \Vert |V^\bbQ - V^\bbP|_{\onevar} \Vert_{L^2(\pi)}) + \Vert H \Vert_{2, \infty} \Vert |V^\bbQ - V^\bbP|_{\onevar} \Vert_{L^2(\pi)},
    \end{align*} where we have used H\"older's inequality for $p = 2/\alpha \geq 2$ with H\"older conjugate $p^* \leq 2$ and monotonicity of $L^p$-norms. Furthermore, the BDG inequality yields \[
    \Vert |X - Y|_\infty \Vert_{L^2(\pi)} \leq \Vert |V^\bbP - V^\bbQ|_{\onevar} \Vert_{L^2(\pi)} + 2\sqrt{ \Vert \tr[M^\bbP - M^\bbQ]_1 \Vert_{L^1(\pi)}}.
    \] Hence, taking infimum over $\pi \in \Pi_{bc}(\bbP, \bbQ)$ and applying the AQTCI yields \begin{align*}
        &|\bbE^\bbP (H \cdot X)_1 - \bbE^\bbQ (H \cdot X)_1| \\
        \leq &\llbracket H \rrbracket_\alpha 5^{\alpha/2} \cA\cW_{sm}(\bbP, \bbQ)^\alpha (\Vert |V^\bbQ|_{\onevar} \Vert_{L^2(\bbQ)} + \cA\cW_{sm}(\bbP, \bbQ)) + \Vert H \Vert_{2, \infty} \cA\cW_{sm}(\bbP, \bbQ) \\
        \leq &\llbracket H \rrbracket_\alpha 5^{\alpha/2} C^{\alpha/2} \sqrt{\cH(\bbP|\bbQ)}^\alpha (\Vert |V^\bbQ|_{\onevar} \Vert_{L^2(\bbQ)} + \sqrt{C}\sqrt{\cH(\bbP|\bbQ)}) + \Vert H \Vert_{2, \infty} \sqrt{C}\sqrt{\cH(\bbP|\bbQ)}.
    \end{align*} The Donsker-Varadhan duality for entropy yields for any $\lambda > 0$ \begin{align*}
        & \log \bbE^\bbQ \exp( \lambda ((H\cdot X)_1 - \bbE^\bbQ (H \cdot X)_1)) \\
        \leq &\sup_{\bbP \ll \bbQ} \lambda (\bbE^\bbP (H \cdot X)_1 - \bbE^\bbQ (H \cdot X)_1) - \cH(\bbP |\bbQ) \\
        \leq &\sup_{\bbP \ll \bbQ} \lambda(\llbracket H \rrbracket_\alpha 5^{\alpha/2} C^{\alpha/2} \sqrt{\cH(\bbP|\bbQ)}^\alpha (\Vert |V^\bbQ|_{\onevar} \Vert_{L^2(\bbQ)} + \sqrt{C}\sqrt{\cH(\bbP|\bbQ)}) + \Vert H \Vert_{2, \infty} \sqrt{C}\sqrt{\cH(\bbP|\bbQ)}) - \cH(\bbP | \bbQ) \\
        \leq & \sup_{\zeta \geq 0} \lambda g_\alpha(\zeta) - \zeta.
    \end{align*}
    
    
    In particular, for any $\lambda, r \geq 0$, Markov's inequality yields \begin{align*}
        \bbQ((H\cdot X)_1 - \bbE^\bbQ (H \cdot X)_1 \geq r) \leq & \bbQ(\exp(\lambda ((H\cdot X)_1 - \bbE^\bbQ (H \cdot X)_1))\geq \exp (\lambda r))\\
        \leq & \exp(-\lambda r) \bbE^\bbQ \exp(\lambda ((H\cdot X)_1 - \bbE^\bbQ (H\cdot X)_1)) \\
        \leq & \exp(- \lambda r + \sup_{\zeta > 0} \lambda g_\alpha(\zeta) - \zeta).
        \end{align*} Taking infimum over $\lambda$ and using Lemma \ref{lem:infsupprob} when $\alpha \in (0,1)$ yields \[
        \bbQ((H\cdot X)_1 - \bbE^\bbQ (H \cdot X)_1 \geq r) \leq \exp(\inf_{\lambda > 0} \sup_{\zeta > 0} \lambda (g_\alpha(\zeta) - r) - \zeta ) = \exp(-g^{-1}_\alpha(r)).
        \] Up to replacing $H$ by $-H$, we obtain the desired. 
        
        When $\alpha = 1$ and $0 < \lambda < \frac{1}{B_1}$, then \[
        \log \bbE^\bbQ \exp( \lambda((H \cdot X)_1 - \bbE^\bbQ(H \cdot X)_1)) \leq \sup_{\zeta \geq 0} \lambda g_1(\zeta) - \zeta = \frac{\lambda^2 (A_1 + D)^2}{4(1 - \lambda B_1)}
        \] with accompanying tail bound \[
        \bbQ(|(H \cdot X)_1 - \bbE^\bbQ(H \cdot X)_1| \geq r) \leq 2\exp( \inf_{0 < \lambda < 1/B_1}\frac{\lambda^2 (A_1 + D)^2}{4(1 - \lambda B_1)} - \lambda r ).
        \] Solving the optimization on the RHS via standard calculus yields the desired inequality.
\end{proof}

\subsection{Dynamic risk measures}

We can also obtain concentration inequalities 
analogous to those in Propositions \ref{prop:tailbehaviorforstointsemimartbddintegrand} and \ref{prop:tailbehaviorforstointmartsqintintegrand} with expectation replaced by certain risk measures.

Given $b \in \scrA^{n,1},\sigma \in \scrA^{n,n}$ such that $\sigma$ is invertible, we assume that \eqref{eqn:sde} has a pathwise unique strong solution and set $\bbQ := \mu^{\sigma,b}$ to be the law of this SDE. Under such an assumption, the canonical process $X$ can always be ``inverted" to obtain a $\cF^X$-Brownian motion $B$ such that $X$ is a solution to \eqref{eqn:sde} driven by $B$ and $\cF^X = \cF^B$ (see Proposition \ref{prop:sdebcbmap}).

Let $g : [0,1] \times \cW^n \times \bbR^n \rightarrow \bbR$ be $\texttt{Prog}_n \otimes \cB(\bbR^n)$-measurable satisfying the following assumption: 
\begin{asmp}
\label{asmp:bsdegengasmp}
\mbox{}
\begin{enumerate}
    \item $\bbE^\bbQ \int_0^1 |g(t,0)|^2 dt < \infty$,
    \item there exists a constant $K \geq 0$ such that $\bbQ \otimes \Leb$-a.s., for all $z,z' \in \bbR^n$, \[
    |g(t,z) - g(t,z')| \leq K|z - z'|,
    \]
    \item $z \mapsto g(t,z)$ is convex, i.e. for all $z, z' \in\bbR^n$, $\lambda \in (0,1)$, \[
    g(t,\lambda z + (1-\lambda)z') \leq \lambda g(t,z) + (1-\lambda) g(t,z'),
    \]
    \item there exist constants $a,b > 0$ such that $\bbQ \otimes \Leb$-a.s., for all $z \in \bbR^n$, $g(t,z)\leq a + b |z|^2$.  
\end{enumerate}
\end{asmp}
 It is well-known that under $1.$ and $2.$ above, for any $\xi \in L^2(\bbQ)$, there exists a unique solution to the backward stochastic differential equation (BSDE) associated with $(g, \xi)$, i.e. there exists a pair of $\cF^X$-progressively measurable processes $(Y,Z) = (Y_t,Z_t)_{t \in [0,1]}$ with values in $\bbR \times \bbR^{1 \times n}$ such that \[
\begin{cases}
    \bbE^\bbQ \sup_{t \in [0,1]} |Y_t|^2 < \infty, \bbE^\bbQ  \int_0^1 |Z_t|^2 dt <\infty \\
    Y_t = \xi + \int_t^1 g(s,Z_s) ds - \int_t^1 Z_s dB_s, \quad t \in [0,1], \quad \bbQ\text{-a.s..}
\end{cases}
\] We refer readers to \cite{zhang2017backward} for more background on BSDEs. Let $\cE : L^2(\bbQ) \rightarrow \bbR,  \xi \mapsto \cE(\xi) := Y_0$. Then, with the addition of $(3)$, \cite{barrieuelkaroui2008pricing} demonstrate that $\cE$ is a convex risk measure with dual representation \begin{equation}
    \label{eqn:dualrepconvexriskmeasure}
    \cE(\xi) = \sup_{h \in \ocH_K} \bbE^{\bbP_h}(\xi) - \bbE^{\bbP_h}(\int_0^1 g^*(s,h_s) ds),
\end{equation} where $g^*$ is the Legendre-Fenchel transform of $g$, \[
g^*(t,\alpha) = \sup_{z \in \bbR^n} \{\alpha \cdot z - g(t,z) \}
\] and $\bbP_h \ll \bbQ$ defined as \[
\frac{d\bbP_h}{d\bbQ} := \exp\left ( \int_0^1 h_s dB_s - \frac{1}{2} \int_0^1 |h_s|^2 ds \right ).
\] In fact, one can compute that \begin{equation}
    \label{eqn:convexriskrelentvalue}
    \cH(\bbP_h | \bbQ) = \frac{1}{2}\bbE^{\bbP_h} \int_0^1 |h_s|^2 ds \leq \frac12 K^2.
\end{equation}

\begin{prop}
    \label{prop:convexriskmeasureconcinequality}
    Suppose $\bbQ$ as above lies in $\cS \cM_2(\cW^n)$ and satisfies the semimartingale AQTCI with constant $C > 0$. Then, for any bounded $H \in \scrA^{n,1}$, we have $(H \cdot X)_1 \in L^2(\bbQ)$ and for any $\lambda \geq 0$, \[
    \cE(\lambda (H \cdot X)_1) \leq \lambda \bbE^\bbQ (H\cdot X)_1 + a + \lambda 2 \Vert H \Vert_\infty \bbE^\bbQ |V^\bbQ|_{\onevar}+ \frac{\lambda^2 \Vert H \Vert_\infty^2 C b}{2}. 
    \]
\end{prop}

\begin{proof}

$(H \cdot X)_1 \in L^2(\bbQ)$ follows from the fact that $\bbQ \in \cS \cM_2(\cW^n)$ and \[
\bbE^\bbQ(H\cdot X)_1^2 \leq 2\Vert H \Vert_{\infty}^2 \bbE^\bbQ( |V^\bbQ|_{\onevar}^2 + \tr[M^\bbQ]_1).
\] 

To see the second assertion, we first note that Assumption \ref{asmp:bsdegengasmp} $4.$ yields a quadratic lower bound on $g^*$: \begin{equation}
    \label{eqn:quadlowerboundonconjofg}
        g^*(t,\alpha)  \geq \sup_{z \in \bbR^n} \{ \alpha \cdot z - a -b |z|^2\} = \frac{|\alpha|^2}{4b} - a.
\end{equation}
Computations similar to those leading to \eqref{eqn:stointAQTCIbound} in Proposition \ref{prop:tailbehaviorforstointsemimartbddintegrand} yield for any $h \in \ocH_K$, \[
|\bbE^{\bbP_h} (H \cdot X)_1 - \bbE^\bbQ (H \cdot X)_1| \leq \Vert H \Vert_\infty(2 \bbE^\bbQ |V^\bbQ|_{\onevar} + \sqrt{C\cH(\bbP_h| \bbQ)}).
\]
    By \eqref{eqn:dualrepconvexriskmeasure}, \eqref{eqn:quadlowerboundonconjofg} and \eqref{eqn:convexriskrelentvalue}, \begin{align*}
            &\cE(\lambda (H \cdot X)_1) \\
         \leq &\sup_{h \in \ocH_K} \bbE^{\bbQ} [\lambda (H \cdot X)_1] +\lambda 2 \Vert H \Vert_\infty \bbE^\bbQ |V^\bbQ|_{1\text{-}\var} + \lambda \Vert H \Vert_\infty \sqrt{C\cH(\bbP_h|\bbQ)} - \bbE^{\bbP_h}[\int_0^1 g^*(s,h_s) ds] \\
        \leq &\lambda \bbE^\bbQ [(H \cdot X)_1]+ \lambda 2 \Vert H \Vert_\infty \bbE^\bbQ |V^\bbQ|_{1\text{-}\var} + \sup_{h \in \ocH_K} \lambda \Vert H \Vert_\infty \sqrt{C\cH(\bbP_h|\bbQ)} + a - \frac{1}{4b} \bbE^{\bbP_h}[\int_0^1 |h_s|^2 ds ] \\
        = &\lambda \bbE^\bbQ [(H \cdot X)_1] +a + \lambda 2 \Vert H \Vert_\infty \bbE^\bbQ |V^\bbQ|_{1\text{-}\var} + \sup_{h \in \ocH_K} \lambda \Vert H \Vert_\infty \sqrt{C\cH(\bbP_h|\bbQ)} - \frac{1}{2b} \cH(\bbP_h|\bbQ) \\
        \leq &\lambda \bbE^\bbQ [(H \cdot X)_1] +a + \lambda 2 \Vert H \Vert_\infty \bbE^\bbQ |V^\bbQ|_{1\text{-}\var}+ \sup_{\zeta  \geq 0} \lambda \Vert H \Vert_\infty\sqrt{C\zeta} - \frac{1}{2b} \zeta \\
        = & \lambda \bbE^\bbQ [(H \cdot X)_1] +a + \lambda 2 \Vert H \Vert_\infty \bbE^\bbQ |V^\bbQ|_{1\text{-}\var}+ \frac{\lambda^2\Vert H \Vert_\infty^2 C b}{2}.
    \end{align*}
\end{proof}

\subsection{Non-anticipative functionals}
As a final application of the semimartingale AQTCI, consider a functional $F : [0,1] \times \cW^n \rightarrow \bbR$, which is assumed to be non-anticipative in the sense that for all $t \in [0,1]$ and $\omega \in \cW^n$, \begin{equation}
    \label{eqn:nonanticipative}
    F(t,\omega) = F(t,\omega_{\cdot \wedge t}).
\end{equation} In \cite{Dupire2019}, the notions of horizontal and vertical differentiation were introduced in order to establish a functional It\^o formula. To be more specific, we momentarily consider the space of c\`adl\`ag functions $\cD^n := \cD([0,1], \bbR^n)$ and some functional $F : [0,1] \times \cD^n \rightarrow \bbR$, for which we make the following assumption: 
\begin{asmp}
    \label{asmp:nonanticipativefnls}
    \begin{enumerate}
    \mbox{}
        \item (Non-anticipativity) For each $\omega \in \cD^n$ and $t \in [0,1]$, \eqref{eqn:nonanticipative} is satisfied. 
        \item (Horizontal differentiability) For all $(t,\omega) \in [0,1) \times \cD^n$, the map $h \mapsto u(h) := F(t+h, \omega_{\cdot \wedge t})$ is right-differentiable at $0$ with derivative $DF(t,\omega) := u'(0)$. 
        \item (Vertical differentiability) For all $(t,\omega) \in [0,1] \times \cD^n$, the map $e \mapsto v(e) := F(t, \omega_{\cdot \wedge t} + e \mathbbm{1}_{[t,1]})$ is twice differentiable at $0$ with first and second derivatives $\nabla F(t,\omega) := v'(0)$ and $\nabla^2 F(t,\omega) := v''(0)$ respectively. 
        \item (Continuity) $F, DF, \nabla F$ and $ \nabla^2 F$ are continuous w.r.t. the stopped path metric, i.e. for $G = F, DF, \nabla F, \nabla^2 F$, for all $t \in [0,1], \omega \in \cD^n$ and $\epsilon > 0$, there exists $\delta > 0$ such that for all $s \in [0,1]$ and $\omega' \in \cD^n$, \[
        |t - s| + |\omega_{\cdot \wedge t} - \omega'_{\cdot \wedge s}|_\infty \leq \delta \implies |G(t,\omega) - G(s,\omega')| \leq \epsilon.
        \] In particular, the map $F(\omega) := (t \mapsto F(t,\omega))$ maps $\cW^n$ into $\cW^1$ whenever $F(0,0) = 0$. 
        \item (Boundedness) $DF, \nabla F$ and $\nabla^2 F$ are boundedness preserving, i.e. for $G = DF, \nabla F, \nabla^2 F$, any compact $\scrK \subseteq \bbR^n$ and $t \in [0,1)$, there exists $C_{\scrK,t} > 0$ such that for all $s \in [0,t]$ and $\omega \in \cD^n$, \[
        \omega([0,s]) \subseteq \scrK \implies |G(s,\omega)| \leq C_{\scrK,t}.
        \]
            
        \item (Regularity) There exists some $K > 0$ such that $DF,\nabla F$ and $\nabla^2 F$ satisfy \[
        \sup_{(t,\omega) \in [0,1] \times \cW^n} |DF(t,\omega)| + |\nabla F(t,\omega)| + |\nabla^2 F(t,\omega)| \leq K 
        \] and for all $\omega, \omega' \in \cW^n$ and $t \in [0,1]$, \[
        \max\{ |DF(t,\omega) - DF(t,\omega')|,|\nabla F(t,\omega) - \nabla F(t,\omega')|, |\nabla^2 F(t,\omega) - \nabla^2 F(t,\omega')| \} \leq K|\omega_{\cdot \wedge t} - \omega'_{\cdot \wedge t}|_\infty.
        \]
    \end{enumerate}
\end{asmp}

Given any semimartingale $X$ and $F$ satisfying Assumption \ref{asmp:nonanticipativefnls} $1.$-$5.$, the functional It\^{o} formula \citep[Theorem $4.1$]{contfournie2013functionalito} yields that $(F(t,X))_{t \in [0,1]}$ is also a semimartingale with decomposition \[
F(t,X) = F(0,X) + \int_0^t DF(s,X) ds + \frac{1}{2} \int_0^t \tr{(\nabla^2 F(s,X) d[X]_s)} + \int_0^t \nabla F(s,X) dX_s. 
\] Using this, we can derive transport inequalities for the law of $F$ under $\bbQ \in \cS\cM_2(\cW^n)$ when $\bbQ$ satisfies the semimartingale AQTCI and some integrability conditions on its characteristics.

\begin{prop}
\label{prop:W1TCIfornonanticipativefnls}
    Fix $p \in [1,2]$. Suppose $\bbQ \in \cS\cM_2(\cW^n)$ satisfies the semimartingale AQTCI with constant $C > 0$. Furthermore, if $p \in (1,2]$, we assume additionally that the drift satisfies $\Vert |V^\bbQ|_{1-\var} \Vert_{L^{\frac{2p}{2-p}}(\bbQ)} < \infty$ and the quadratic variation satisfies \begin{enumerate}
        \item $\bbE^\bbQ \exp(\varsigma(\tr[M^\bbQ]_1)^{\frac{p}{2-p}}) < \infty$ for some $\varsigma > 0$ when $p \in (1,2)$,
        \item $\Vert \tr[M^\bbQ]_1 \Vert_{L^{\infty}(\bbQ)} < \infty$ when $p = 2$.
    \end{enumerate}  Then, for any $F$ satisfying $F(0,0) = 0$ and Assumption \ref{asmp:nonanticipativefnls},  $ F_{\#}\bbQ$ satisfies the following transport inequality: for all $\bbP' \in \cP(\cW^1)$, \[
    \cW_{|\cdot|_\infty,p}(\bbP',  F_{\#}\bbQ) \leq \sqrt{C}C'_p (1 + \cH(\bbP'|F_{\#}\bbQ)^{\frac{2-p}{2p}}) \sqrt{\cH(\bbP'| F_{\#}\bbQ)},
    \] where \begin{equation*}
        C'_p := K\sqrt{5}(2 + \Vert \tr[M^\bbQ]_1 \Vert_{L^{\frac{2p}{2-p}}(\bbQ)} + (1 + b_p)\Vert \sqrt{\tr[M^\bbQ]_1}\Vert_{L^{\frac{2p}{2-p}}(\bbQ)} + \Vert |V^\bbQ|_{\onevar} \Vert_{L^{\frac{2p}{2-p}} (\bbQ)} \\
        + A_p),
    \end{equation*} $b_p$ is the BDG inequality constant and \begin{enumerate}
        \item $A_1 = \Vert \sqrt{\tr[M^\bbQ]_1}\Vert_{L^{2}(\bbQ)} + \sqrt{C}$,
        \item $A_p = \frac{1}{\varsigma^{\frac{2-p}{2p}}}( 1 + (\log \bbE^\bbQ \exp(\varsigma (\tr[M^\bbQ]_1)^{\frac{p}{2-p}}))^{\frac{2-p}{2p}}))$ for $p \in (1,2)$,
        \item $A_2 = \Vert \sqrt{\tr[M^\bbQ]_1}\Vert_{L^{\infty}(\bbQ)}$.
    \end{enumerate} 
\end{prop}

\begin{rmk}
    \label{rmk:NAFLipschitzcounterexample}
    If $F : \cW^n \rightarrow \cW^1$ is Lipschitz, then since $\cW_{|\cdot|_\infty,2} \leq \cA \cW_{|\cdot|_\infty,2} \leq \sqrt{8 } \cA \cW_{sm}$, it is straightforward to establish the classical QTCI for $F_{\#}\bbQ$ by \cite[Lemma $2.1$]{djelloutguillinwu2004transportationci}. However, Assumption \ref{asmp:nonanticipativefnls} does not imply that $F$ is Lipschitz since the notion of vertical derivative only controls sensitivity to perturbations at the end-point of paths, instead of uniform perturbations of the entire path trajectory (see Example \ref{exam:NAFnotlipschitzbutregular} in Appendix \ref{appendix:misc}).
\end{rmk}

\begin{rmk}
    \label{rmk:NAFfiltration}
    The preceding proposition only establishes a classical transport inequality since in general, pushing $\pi \in \Pi_{bc}(\bbP,\bbQ)$ forward by $\left( F,F \right) $ need not preserve bicausality.
    However, if $\cF^F = \cF^X$, where $\cF^F_t := \sigma(\{F(s,X) : s \leq t\})$, then $\pi \in \Pi_{bc}(\bbP,\bbQ) \implies (F,F)_{\#}\pi \in \Pi_{bc}(F_{\#}\bbP, F_{\#}\bbQ)$ and we can conclude an adapted transport inequality as the proof will show.
\end{rmk}

\begin{proof}
    Fix $\bbP \ll \bbQ$ such that $\cH(\bbP|\bbQ) < \infty$ and $\pi \in \Pi_{bc}(\bbP, \bbQ)$. By the functional It\^o formula, the processes $(F(t,X))_{t \in [0,1]}$ and $(F(t,Y))_{t \in [0,1]}$ are semimartingales admitting decompositions \[
    F(t,X) = F(0,0) + \int_0^t DF(s,X) ds + \frac{1}{2} \int_0^t \tr(\nabla^2 F(s,X) d[X]_s) + \int_0^t \nabla F(s,X) dX_s
    \] with a similar expression for $F(\cdot,Y)$. We compute \begin{align*}
        &\Vert |F(\cdot,X) - F(\cdot,Y)|_\infty \Vert_{L^p(\pi)} \\
        \leq & \Vert |\int_0^\cdot DF(s,X) - DF(s,Y) ds|_\infty \Vert_{L^p(\pi)} + \frac{1}{2} \Vert |\int_0^\cdot \tr(\nabla^2 F(s,X) d[X]_s) - \int_0^\cdot \tr(\nabla^2 F(s,Y) d[Y]_s)|_\infty \Vert_{L^p(\pi)} \\
        &+ \Vert |\int_0^\cdot \nabla F(s,X) dX_s - \int_0^\cdot \nabla F(s,Y) dY_s|_\infty \Vert_{L^p(\pi)} =: \text{I} + \text{II} + \text{III}.
    \end{align*} We estimate each term separately. Using $1 \leq p \leq 2$, the BDG inequality, monotonicity of $L^p$-norms and Lipschitz property of $DF$, we obtain \begin{align*}
        \text{I} = \Vert |\int_0^\cdot DF(s,X) - DF(s,Y) ds|_\infty \Vert_{L^p(\pi)} \leq &K \Vert |X - Y|_\infty \Vert_{L^p(\pi)} \\
        \leq &K  \Vert |X - Y|_\infty \Vert_{L^2(\pi)} \\
        \leq &K(\Vert |V^\bbP - V^\bbQ|_{\onevar} \Vert_{L^2(\pi)} + 2 \sqrt{\Vert \tr[M^\bbP - M^\bbQ]_1 \Vert_{L^{1}(\pi)}}) \\
        \leq &K\sqrt{5} \frD(\pi),
    \end{align*} where $\frD(\pi) := \sqrt{\bbE^\pi (|V^\bbP - V^\bbQ|_{\onevar}^2 + \tr[M^\bbP - M^\bbQ]_1)}.$ Similarly, we estimate using H\"older's inequality, the boundedness and Lipschitz property of $\nabla^2 F$ that \begin{align*}
        &2\text{II}\\
        = & \Vert |\int_0^\cdot \tr(\nabla^2 F(s,X) d[X]_s) - \int_0^\cdot \tr(\nabla^2 F(s,Y) d[Y]_s)|_\infty \Vert_{L^p(\pi)}\\
        \leq &\Vert |\int_0^\cdot \tr((\nabla^2 F(s,Y) - \nabla^2 F(s,X)) d[Y]_s)|_\infty \Vert_{L^p(\pi)} + \Vert |\int_0^\cdot \tr(\nabla^2 F(s,X) d([X]- [Y])_s)|_\infty \Vert_{L^p(\pi)} \\
        \leq &\Vert \Vert |\nabla^2 F(\cdot,X) - \nabla^2 F(\cdot, Y)|_{\op} \Vert_{L^\infty(\Leb)} \tr[Y]_1 \Vert_{L^p(\pi)}\\
        &+ \Vert \Vert |\nabla^2 F(\cdot,X)|_{\op} \Vert_{L^\infty(\Leb)} \sqrt{\tr[X-Y]_1}(\sqrt{\tr[X]_1} + \sqrt{\tr[Y]_1}) \Vert_{L^p(\pi)} \\
        \leq &K (\Vert |X - Y|_\infty \tr[M^\bbQ]_1 \Vert_{L^p(\pi)} + \Vert \sqrt{\tr[M^\bbQ-M^\bbP]_1}(\sqrt{\tr[M^\bbP]_1} + \sqrt{\tr[M^\bbQ]_1}) \Vert_{L^p(\pi)} ) \\
        \leq &K \Vert |X - Y|_\infty \Vert_{L^{2}(\pi)} \Vert \tr[M^\bbQ]_1 \Vert_{L^{\frac{2p}{2-p}}(\bbQ)} \\
        &+ K\Vert \sqrt{\tr[M^\bbP-M^\bbQ]_1} \Vert_{L^2(\pi)} (\Vert \sqrt{\tr[M^\bbQ]_1}\Vert_{L^{\frac{2p}{2-p}}(\bbQ)} + \Vert \sqrt{\tr[M^\bbP]_1} \Vert_{L^{\frac{2p}{2-p}}(\bbP)}) \\
        \leq &K (\Vert \tr[M^\bbQ]_1 \Vert_{L^{\frac{2p}{2-p}}(\bbQ)} \sqrt{5} + \Vert \sqrt{\tr[M^\bbQ]_1}\Vert_{L^{\frac{2p}{2-p}}(\bbQ)} + \Vert \sqrt{\tr[M^\bbP]_1} \Vert_{L^{\frac{2p}{2-p}}(\bbP)} ) \frD(\pi).
    \end{align*} Note that since $\bbQ$ satisfies the semimartingale AQTCI, Corollary \ref{cor:QTCIforQV} implies that the $L^p$-norms of $\tr[M^\bbQ]_1$ and $ \sqrt{\tr[M^\bbQ]_1}$ are always finite and if \begin{enumerate}
        \item $p = 1$, then \[
        \Vert \sqrt{\tr[M^\bbP]_1} \Vert_{L^{2}(\bbP)} \leq \Vert \sqrt{\tr[M^\bbP - M^\bbQ]_1} \Vert_{L^{2}(\pi)} + \Vert \sqrt{\tr[M^\bbQ]_1} \Vert_{L^{2}(\bbQ)} \leq \frD(\pi) + \Vert \sqrt{\tr[M^\bbQ]_1} \Vert_{L^{2}(\bbQ)}
        \] and hence, taking infimum over $\pi \in \Pi_{bc}(\bbP, \bbQ)$, \begin{equation} 
            \label{eqn:p=1transfer}
             \Vert \sqrt{\tr[M^\bbP]_1} \Vert_{L^{2}(\bbP)} \leq \sqrt{C\cH(\bbP|\bbQ)} + \Vert \sqrt{\tr[M^\bbQ]_1} \Vert_{L^{2}(\bbQ)}
        \end{equation}
        \item $p \in (1,2)$, then since $\bbP$-a.s. $[M^\bbP] = [M^\bbQ]$ and  \cite[Proposition $B.1$]{gozlanleonard2010transport} gives \[
            \bbE^\bbP (\tr[M^\bbP]_1)^{\frac{p}{2-p}} = \bbE^\bbP (\tr[M^\bbQ]_1)^{\frac{p}{2-p}} \leq \frac{1}{\varsigma} (\cH(\bbP|\bbQ) + \log \bbE^\bbQ \exp(\varsigma (\tr[M^\bbQ]_1)^{\frac{p}{2-p}})) < \infty,
            \] we obtain \begin{equation}
            \label{eqn:expintimpliesmomenttransfer}
            \Vert \sqrt{\tr[M^\bbP]_1} \Vert_{L^{\frac{2p}{2-p}}(\bbP)} \leq \frac{1}{\varsigma^{\frac{2-p}{2p}}}( \cH(\bbP|\bbQ)^{\frac{2-p}{2p}} + (\log \bbE^\bbQ \exp(\varsigma (\tr[M^\bbQ]_1)^{\frac{p}{2-p}}))^{\frac{2-p}{2p}})) < \infty.
    \end{equation}
        \item $p = 2$, then since $\bbP \ll \bbQ$, we have \begin{equation}
            \label{eqn:p=2transfer}
            \Vert \sqrt{\tr[M^\bbP]_1} \Vert_{L^{\infty}(\bbP)} \leq \Vert \sqrt{\tr[M^\bbQ]_1} \Vert_{L^{\infty}(\bbQ)}.
        \end{equation}
    \end{enumerate}  Finally, we estimate \begin{align*}
        \text{III} = &\Vert |\int_0^\cdot \nabla F(s,X) dX_s - \int_0^\cdot \nabla F(s,Y) dY_s|_\infty \Vert_{L^p(\pi)} \\
        \leq &\Vert |\int_0^\cdot \nabla F(s,X) d(X-Y)_s|_\infty \Vert_{L^p(\pi)} + \Vert |\int_0^\cdot (\nabla F(s,X) - \nabla F(s,Y) )dY_s|_\infty \Vert_{L^p(\pi)} \\
        \leq &\Vert \Vert \nabla F(\cdot, X) \Vert_{L^\infty(\Leb)} |V^\bbQ - V^\bbP|_{\onevar} \Vert_{L^2(\pi)} + 2 \sqrt{ \Vert \Vert \nabla F(\cdot, X) \Vert_{L^\infty(\Leb)}^2 \tr[M^\bbQ - M^\bbP]_1 \Vert_{L^{1}(\pi)}} \\
        &+ \Vert \Vert \nabla F(\cdot, X) - \nabla F(\cdot,Y) \Vert_{L^\infty(\Leb)} |V^\bbQ|_{\onevar} \Vert_{L^p(\pi)} \\
        &+ b_p \sqrt{ \Vert \Vert \nabla F(\cdot, X) - \nabla F(\cdot,Y) \Vert_{L^\infty(\Leb)}^2 \tr[M^\bbQ]_1 \Vert_{L^{p/2}(\pi)}} \\
        \leq &K (\sqrt{5} \frD(\pi) +  \Vert |X - Y|_\infty |V^\bbQ|_{\onevar} \Vert_{L^p(\pi)} + b_p \Vert |X - Y|_\infty \sqrt{\tr[M^\bbQ]_1} \Vert_{L^{p}(\pi)}) \\
        \leq & K( \sqrt{5} \frD(\pi) + (\Vert |V^\bbQ|_{\onevar} \Vert_{L^{\frac{2p}{2-p}} (\bbQ)} + b_p \Vert \sqrt{\tr[M^\bbQ]_1}\Vert_{L^{\frac{2p}{2-p}}(\bbQ)} )\Vert |X-Y|_\infty \Vert_{L^2(\pi)} ) \\
        \leq & K\sqrt{5}(1 + \Vert |V^\bbQ|_{\onevar} \Vert_{L^{\frac{2p}{2-p}} (\bbQ)} + b_p \Vert \sqrt{\tr[M^\bbQ]_1}\Vert_{L^{\frac{2p}{2-p}}(\bbQ)} )\frD(\pi).
    \end{align*} Hence, for $p \in [1,2]$, we obtain a final estimate of \[
    \cW_{|\cdot|_\infty, p}(F_{\#}\bbP, F_{\#}\bbQ) \leq \Vert |F(\cdot,X) - F(\cdot,Y)|_\infty \Vert_{L^p(\pi)} \leq C'_p (1 + \cH(\bbP|\bbQ)^{\frac{2-p}{2p}}) \frD(\pi).
    \] Taking infimum over $\pi \in \Pi_{bc}(\bbP, \bbQ)$ appropriately, we obtain \[
    \cW_{|\cdot|_\infty, p}(F_{\#}\bbP, F_{\#}\bbQ) \leq C'_p(1 + \cH(\bbP|\bbQ)^{\frac{2-p}{2p}}) \cA\cW_{sm}(\bbP,\bbQ) \leq C'_p \sqrt{C}(1 + \cH(\bbP|\bbQ)^{\frac{2-p}{2p}}) \sqrt{\cH(\bbP|\bbQ)}.
    \] By \cite[Lemma $2.1$]{djelloutguillinwu2004transportationci}, $\cH(\bbP' | F_{\#}\bbQ) = \inf\{ \cH(\bbP|\bbQ) : F_{\#}\bbP = \bbP' \}$ for any $\bbP' \in \cP(\cW^1)$, which, together with the preceding estimates, yields the desired.
\end{proof}

\section{Adapted QTCI beyond diffusions}
\label{section:AQTCIbeyonddiffusions}

\subsection{Relative entropy as a causal transport problem}


We first work in an abstract setting, fixing some $\bbQ \in \cP(\cW^n)$ such that $\cT_{bcb}(\bbW^n, \bbQ)$ is non-empty. Let $T \in \cT_{bcb}(\bbW^n, \bbQ)$ with $\bbW^n$-a.s. inverse $S$. \begin{lem}
    \label{lem:imageofbicausalcpls}
    For any $\bbP \in \cP(\cW^n)$ satisfying $\bbP \ll \bbW^n$, we have \[
    \Pi_{bc}(\bbQ, T_{\#}\bbP) = \{ (T,T)_{\#}\pi : \pi \in \Pi_{bc}(\bbW^n, \bbP) \}.
    \]
\end{lem} 

\begin{proof}
    We first note that if $\bbP \ll \bbW^n$, then $T_{\#}\bbP \ll T_{\#}\bbW^n = \bbQ$. Furthermore, since $T \in \cT_{bcb}(\bbW^n, \bbQ)$, then by definition,  \[
    S \circ T = \identity \; \bbW^n\text{-a.s. and } T \circ S = \identity \; \bbQ\text{-a.s.}
    \] and these identities also hold $\bbP$-a.s. and $\bbT_{\#}\bbP$-a.s. respectively. We further note that $(\identity,T)_{\#}\bbW^n \in \Pi_c(\bbW^n, \bbQ)$ implies that \[
    A \in \cF^X_t \implies T^{-1}(A) \in \prescript{\bbW^n}{}{\cF}^{X}_t \subseteq \prescript{\bbP}{}{\cF}^{X}_t 
    \] so that $(\identity,T)_{\#}\bbP \in \Pi_c(\bbP, T_{\#}\bbP)$. Similarly, $(\identity,S)_{\#}\bbQ \in \Pi_c(\bbQ, \bbW^n) \implies (\identity,S)_{\#}(T_{\#}\bbP) \in \Pi_c(T_{\#}\bbP, \bbP).$ Together, we obtain that $(\identity,T)_{\#}\bbW^n \in \Pi_{bc}(\bbW^n, \bbQ) \implies (\identity,T)_{\#}\bbP \in \Pi_{bc}(\bbP, T_{\#}\bbP).$ Hence, if $\pi \in \Pi_{bc}(\bbW^n, \bbP)$, then since bicausality is stable under gluing \citep[Lemma $3.2$]{backhoff2020adaptedwd}, we have $(T,T)_{\#}\pi \in \Pi_{bc}(\bbQ, T_{\#}\bbP)$. On the other hand, given $\pi \in \Pi_{bc}(\bbQ, T_{\#}\bbP)$, noting that $(\identity, S)_{\#}\bbQ \in \Pi_{bc}(\bbQ, \bbW^n)$ and $(\identity, S)_{\#}(T_{\#}\bbP) \in \Pi_{bc}(T_{\#}\bbP, \bbP)$, we again have by gluing that $(S,S)_{\#}\pi \in \Pi_{bc}(\bbW^n, \bbP)$. 
\end{proof}

From the above lemma, it holds that whenever $\bbP \ll \bbW^n$, \begin{equation}
    \label{eqn:bicausaltptproblemequality}
    \inf_{\pi \in \Pi_{bc}(\bbW^n, \bbP)}  \bbE^\pi c_{\mathfrak{h}} = \inf_{\pi \in \Pi_{bc}(\bbQ, T_{\#}\bbP)}  \bbE^\pi c^T_{\mathfrak{h}}
\end{equation} for $c^T_{\mathfrak{h}} := c_{\mathfrak{h}} \circ (S,S).$


\begin{prop}
    \label{prop:relentequality}
    For any $\bbP \in \cP(\cW^n)$ such that $\bbP \ll \bbQ$,  we have $S_{\#}\bbP \ll \bbW^n$ and \begin{equation}
        \label{eqn:relendequality}
        \inf_{\pi \in \Pi_{bc}(\bbQ, \bbP)} \bbE^\pi c^T_{\mathfrak{h}} = \inf_{\pi \in \Pi_{bc}(S_{\#}\bbP, \bbW^n)} \bbE^\pi c_{\mathfrak{h}} = 2\cH(S_{\#}\bbP|\bbW^n) = 2 \cH(\bbP|\bbQ).
    \end{equation}
\end{prop}

\begin{rmk}
    \label{rmk:necessityofabscty}
    We note that it is not possible in general to drop the assumption that $\bbP \ll \bbQ$ and still obtain equality in \eqref{eqn:relendequality} since $T \circ S = \identity$ only holds $\bbQ$-a.s.. In particular, if $\bbP$ is singular w.r.t. $\bbQ$, it is possible that $\cH(S_{\#}\bbP|\bbW^n) = 0$ but $\cH(\bbP|\bbQ) = \infty$.
\end{rmk}

\begin{proof}
    That $\bbP \ll \bbQ \implies S_{\#}\bbP \ll \bbW^n$ is straightforward. The first equality is exactly \eqref{eqn:bicausaltptproblemequality}, while the second equality in \eqref{eqn:relendequality} follows from \cite[Lemma $5$]{lassalle2018causaltp}. The third equality follows from the fact that relative entropy decreases under pushforwards by the data processing inequality so that \[
    \cH(\bbP|\bbQ) \geq \cH(S_{\#}\bbP|S_{\#}\bbQ)= \cH(S_{\#}\bbP|\bbW^n) \geq \cH(T_{\#}S_{\#}\bbP|\bbQ) = \cH(\bbP|\bbQ).
    \] That $T_{\#}S_{\#}\bbP = (T \circ S)_{\#}\bbP = \bbP$ follows from $\bbP \ll \bbQ$ and $T \circ S = \identity$ $\bbQ$-a.s.. 
\end{proof}

The preceding proposition can be understood as an abstract generalization of \eqref{eqn:lassalletptinequality}:
\begin{equation}
    \label{eqn:abstractrelativeentropy}
    2 \cH(\bbP | \bbQ) = \inf_{\pi \in \Pi_{bc}(\bbQ, \bbP)} \bbE^\pi c_{\mathfrak{h}}^T \text{ for any } \bbP \ll \bbQ.
\end{equation} In fact, an optimal coupling for the problem on the RHS of \eqref{eqn:abstractrelativeentropy} can be made explicit using the notion of a Girsanov shift defined in \cite{lassalle2018causaltp}. 

\begin{cor}
    \label{cor:optimalbicausalcpling}
     For any $\bbP \in \cP(\cW^n)$ such that $\bbP \ll \bbQ$, there exists some $\bbR^n$-valued $\cF^X$-adapted process $v = (v_t)_{t \in [0,1]}$ such that the Girsanov shift $V^{T} : \cW^n \rightarrow \cW^n$ given by \begin{equation}
    \label{eqn:girsanovshift}
        V^{T}(X) := X - \int_0^\cdot v(s,X) ds
\end{equation} satisfies $V^{T} \in \cT_{bc}(S_{\#}\bbP, \bbW^n)$ and $\pi^* := (T, T \circ V^{T})_{\#}(S_{\#}\bbP) \in \Pi_{bc}(\bbP, \bbQ)$  satisfies \[
     2 \cH(\bbP|\bbQ) = \inf_{\pi \in \Pi_{bc}(\bbP, \bbQ)} \bbE^\pi c_{\mathfrak{h}}^T = \bbE^{\pi^*} c_{\mathfrak{h}}^T = \bbE^{S_{\#}\bbP} c_{\mathfrak{h}}(V^T,\identity) = \bbE^{S_{\#}\bbP} \int_0^1 |v(s,X)|^2 ds.
     \]
\end{cor}

\begin{proof}
     The first equality follows from Proposition \ref{prop:relentequality} and the third equality is straightforward. Since $S_{\#} \bbP \ll \bbW^n$, then \cite[Proposition $6$, Lemma $5$]{lassalle2018causaltp} gives the existence of the Girsanov shift $V^T \in \cT_{bc}(S_{\#}\bbP, \bbW^n)$ and \[
     \inf_{\pi \in \Pi_{bc}(\bbP, \bbQ)} \bbE^\pi c_{\mathfrak{h}}^T = \inf_{\pi \in \Pi_{bc}(S_{\#}\bbP, \bbW^n)} \bbE^\pi c_{\mathfrak{h}} = \bbE^{S_{\#}\bbP} c_{\mathfrak{h}}(V_T,\identity) = 2\cH(S_{\#}\bbP|\bbW^n) = \bbE^{S_{\#}\bbP} \int_0^1 |v(s,X)|^2 ds.
     \] Furthermore, since $T \in \cT_{bcb}(\bbW^n, \bbQ)$, we easily conclude that $(T \circ V^T, T)_{\#}(S_{\#}\bbP) \in \Pi_{bc}(\bbQ, \bbP)$ as required.
\end{proof}




\begin{exam}
    To make things concrete, we collect some examples of $\bbQ \in \cP(\cW^n)$ such that $\cT_{bcb}(\bbW^n, \bbQ) \neq \emptyset$. Then, $\bbQ$ can be
\begin{enumerate}
    \item the law of a stochastic differential equation \eqref{eqn:sde} with Brownian noise, where $b \in \scrA^{n,1}, \sigma \in \scrA^{n,n}$ are such that $\det(\sigma) \neq 0$ and there exists a unique solution to \eqref{eqn:sde} adapted to $\cF^B$. 
    \item the law of a scalar SDE with fractional Brownian noise \begin{equation}
        \label{eqn:fsde}
        X_t = \int_0^t b(X_{s}) ds + \int_0^t \sigma(X_s) dB^H_s, X_0 = 0,
    \end{equation} where $B^H_t := \int_0^t K_H(t,s) dB_s$, $K_H$ is the Molchan-Golosov kernel for $H \in (0,1)$, $b,\sigma \in \scrM^{1,1}$ are independent of time, $\sigma > 0$ and there exists a unique solution to \eqref{eqn:fsde} adapted to $\cF^B$ in the sense of Doss-Sussmann. In this case, we denote the law of a solution to \eqref{eqn:fsde} as $\varrho^{\sigma,b} \in \cP(\cW^1).$
    \item the law of a scalar stochastic Volterra equation (SVE) \begin{equation}
        \label{eqn:sve}
        X_t = \int_0^t G(t - s) b(s,X_s) ds + \int_0^t G(t - s) \sigma(s,X_s) dB_s, X_0 = 0,
    \end{equation} where $b, \sigma \in \scrM^{1,1}, \sigma > 0$, $(s,t) \mapsto G(t-s)$ is a Volterra kernel with resolvent of the first kind and there exists a unique solution to \eqref{eqn:sve} adapted to $\cF^B$. In this case, we denote the law of a solution to \eqref{eqn:sve} as $\kappa^{\sigma,b}_{G} \in \cP(\cW^1)$.
\end{enumerate} In all of the above examples, the key property is the equality of filtrations of the solution and driving noise: $\cF^X = \cF^B$. This implies that the solution map $T$ sending a Brownian motion path to a solution path lies in $\cT_{bcb}(\bbW^n, \bbQ)$. We refer readers to Appendix \ref{appendix:bcbmaps} for the full details.
\end{exam} In what follows, we will use Corollary \ref{cor:optimalbicausalcpling} with the above examples to conclude the AQTCI for the solutions of \eqref{eqn:fsde} and \eqref{eqn:sve}.

\subsection{Scalar SDEs driven by fractional Brownian noise}

\begin{thm}
    \label{thm:adaptedQTCIfSDEs}\mbox{}
   \begin{enumerate}
       \item Let $b \in \scrA^{1,1}$ satisfy \eqref{eqn:lipschitzsdecoeff} and $\int_0^1 |b(t,0)| dt < \infty$. Then, the law $\varrho^{1,b} \in \cP(\cW^1)$ of the unique solution to \begin{equation}
        \label{eqn:additivede}
            X_t = \int_0^t b(s,X) ds + B^H_t = \int_0^t b(s,X) ds + \int_0^t K_H(t,s) dB_s, \quad X_0 = 0
    \end{equation} satisfies an AQTCI w.r.t. $d = |\cdot|_\infty$ with constant $4 e^{2K^2}$.
    \item Let $b \in C_b^1(\bbR), \sigma \in C^2_b(\bbR)$ and $\delta < \sigma < \frac{1}{\delta}$ for some $\delta > 0$. Then, the law $\varrho^{\sigma,b}$ of the solution to \eqref{eqn:fsde} satisfies an AQTCI w.r.t. $d = |\cdot|_\infty$ with constant $4 \Vert \sigma \Vert_{\infty}^2 e^{2 \tilde{K}^2}$, where $\tilde{K} := (1 + \frac{\Vert b \Vert_\infty}{\delta}) \max \{ \Vert b' \Vert_\infty, \Vert \sigma' \Vert_\infty \} $. 
   \end{enumerate} 
\end{thm}

\begin{proof}

We fix $\bbP \ll \varrho^{\sigma,b}$ such that $\cH(\bbP | \varrho^{\sigma,b}) < \infty$. 

Assuming $1.$, by Proposition \ref{prop:additivegennoise}, the solution map $\frF^{1,b}$ lies in $\cT_{bcb}(\bbW^1, \varrho^{1,b})$. Let $\tilde{\bbP} := (\frF^{1,b})^{-1}_{\#}\bbP.$ By Corollary \ref{cor:optimalbicausalcpling}, there exists some real-valued $\cF^X$-adapted process $v = (v_t)_{t \in [0,1]}$ such that the Girsanov shift $V := V^{\frF^{1,b}}$ defined by \eqref{eqn:girsanovshift} satisfies $V \in \cT_{bc}(\tilde{\bbP}, \bbW^1)$, $(\frF^{1,b} \circ V, \frF^{1,b})_{\#}\tilde{\bbP} \in \Pi_{bc}(\varrho^{1,b},\bbP)$ and \[
2\cH(\bbP | \varrho^{1,b}) = \bbE^{\tilde{\bbP}} \int_0^1 |v(s,X)|^2 ds.
\] Since $\tilde{\bbP} \ll \bbW^1$, $Z := \frF^{1,b}(X)$ satisfies $\bbW^1$-a.s. and hence, $\tilde{\bbP}$-a.s. \[
Z_t = \int_0^t b(s,Z) ds + \int_0^t K_H(t,s) dX_s.
\] Furthermore, since $V_{\#}\tilde{\bbP} = \bbW^1$, i.e. $V$ is a Brownian motion, $\tilde{Z} := \frF^{1,b} \circ V$ satisfies $\tilde{\bbP}$-a.s. \[
\tZ_t = \int_0^t b(s, \tZ) ds + \int_0^t K_H(t,s) dV_s.
\] Noting that $X - V = \int_0^\cdot v(s,X) ds$, it is straightforward to verify that $\tilde{\bbP}$-a.s. \[
\int_0^t K_H(t,s) dX_s = \int_0^t K_H(t,s) dV_s + \int_0^t K_H(t,s)v(s,X) ds,
\] for fixed $t \in [0,1]$ and then, for all $t \in [0,1]$ by continuity of sample paths of fBMs. Hence, we compute using H\"older's inequality that \begin{align*}
    \sup_{s \in [0,t]} |Z_s - \tZ_s|^2      \leq & 2\int_0^t |b(s,Z) - b(s, \tZ)|^2 ds + 2\sup_{s \in [0,t]} |\int_0^s K_H(s,u)v(u,X) du|^2 \\
    \leq  & 2K^2 \int_0^t \sup_{u \in [0,s]} |Z_u - \tZ_u|^2 ds  + 2 (\sup_{s \in [0,t]}\int_0^s K_H(s,u)^2 du ) \int_0^1 |v(s,X)|^2 ds.
\end{align*} Since \[
t^{2H} = \bbE (B^H_t)^2 = \bbE (\int_0^t K_H(t,s) dB_s)^2 = \int_0^t K_H(t,s)^2 ds,
\] we obtain \[
  \bbE^{\tilde{\bbP}} \sup_{s \in [0,t]} |Z_s - \tZ_s|^2 \leq 2 K^2 \int_0^t \bbE^{\tilde{\bbP}} \sup_{u \in [0,s]} |Z_u - \tZ_u|^2 ds +2 \bbE^{\tilde{\bbP}} \int_0^1 |v(s,X)|^2 ds.
\] By Gronwall's inequality, we obtain \[
\bbE^{\tilde{\bbP}} |Z - \tZ|_\infty^2 \leq 4 \cH(\bbP |\varrho^{1, b}) \exp(2K^2).
\] Furthermore, since $\Law(\tZ, Z) = (\frF^{1,b} \circ V, \frF^{1,b})_{\#} \tilde{\bbP} \in \Pi_{bc}(\varrho^{1,b},\bbP)$, we conclude \[
\cA\cW^2_{|\cdot|_\infty}(\bbP, \varrho^{1,b}) \leq 4 e^{2K^2} \cH(\bbP |\varrho^{1, b}).
\]

Assuming $2.$, we first define $h \in C^1(\bbR)$ to be the solution to the ODE \[
\frac{dh}{d\beta} = \sigma(h(\beta)), \quad h(0) = 0,
\] and note that since $\delta < \sigma \leq \Vert \sigma \Vert_{\infty}$, then $h$ is bi-Lipschitz, i.e. for all $\alpha, \beta \in \bbR$, \[
\delta |\alpha - \beta| \leq |h(\alpha) - h(\beta)| \leq \Vert \sigma \Vert_{\infty} |\alpha - \beta|.
\] In particular, the map $\mathbf{h} : \cW^1 \rightarrow \cW^1$ given by $\mathbf{h}(\omega) := (t \mapsto h(\omega_t))$ lies in $\cT_{bcb}(\bbP, \mathbf{h}_{\#}\bbP)$ for any $\bbP \in \cP(\cW^1)$ and satisfies \[
\delta |\omega - \omega'|_\infty \leq |\mathbf{h}(\omega) - \mathbf{h}(\omega')|_\infty \leq \Vert \sigma \Vert_{\infty} |\omega - \omega'|_\infty.
\]
Defining $\tb(x) := \frac{b(h(x))}{\sigma(h(x))}$, Proposition \ref{prop:dosssussmanbcbmap} gives $\mathbf{h}_{\#}\varrho^{1,\tb} = \varrho^{\sigma,b}$. Hence, for any $\bbP' \in \cP(\cW^1)$ and $\pi \in \Pi_{bc}(\bbP', \varrho^{1,\tb})$, we have by gluing that $(\mathbf{h},\mathbf{h})_{\#}\pi \in \Pi_{bc}(\mathbf{h}_{\#}\bbP', \varrho^{\sigma,b})$ and \[
\cA\cW^2_{|\cdot|_\infty}(\mathbf{h}_{\#}\bbP', \varrho^{\sigma,b}) \leq \bbE^{\pi} |\mathbf{h}(X) - \mathbf{h}(Y)|^2_\infty \leq \Vert \sigma \Vert^2_{\infty} \bbE^\pi |X - Y|^2_\infty.
\] Taking infimum over $\pi$ yields \[
\cA\cW^2_{|\cdot|_\infty}(\mathbf{h}_{\#}\bbP', \varrho^{\sigma,b}) \leq \Vert \sigma \Vert^2_{\infty} \cA\cW^2_{|\cdot|_\infty}(\bbP', \varrho^{1,\tb}).
\] Furthermore, since $\tb$ is Lipschitz with Lipschitz constant $\tilde{K} := (1 + \frac{\Vert b \Vert_\infty}{\delta})\max\{ \Vert b' \Vert_\infty, \Vert \sigma' \Vert_\infty\}$, the previous part gives \[
\cA\cW^2_{|\cdot|_\infty}(\bbP', \varrho^{1, \tb}) \leq 4 e^{2 \tilde{K}^2} \cH(\bbP'|\varrho^{1,\tb}).
\] Combining the preceding bound and the following identity \[
\cH(\bbP|\varrho^{\sigma,b}) = \inf\{ \cH(\bbP' | \varrho^{1,\tb}) : \mathbf{h}_{\#}\bbP' = \bbP \}
\] used in the proof of Corollary \ref{cor:QTCIforQV} yields the desired.
\end{proof}

\subsection{Scalar SVEs}

In order to consider SVEs with unique $\cF^B$-adapted solutions, we impose the following conditions on $\sigma,b \in \scrM^{1,1}$ and the Volterra kernel $G( \cdot - \cdot) : [0,1]^2 \rightarrow \bbR_{\geq 0}$: \begin{asmp}
    \label{asmp:svecoeff}
    $b, \sigma \in \scrM^{1,1}$ satisfy \eqref{eqn:lipschitzsdecoeff} and $\sup_{t \in [0,1]} |b(t,0)| + |\sigma(t,0)| < \infty$. Furthermore, there exists some $\beta > 1$, $\kappa < \infty$ and $\theta \in (0,1]$ such that \[
    \sup_{t \in [0,1]} \int_0^t G(t - s)^{2\beta} ds < \infty \text{ and } \sup_{t \in [0,1]} (\int_0^t (G((t + \delta) \wedge 1 - s) - G(t - s))^2 ds)^{1/2} \leq \kappa \delta^\theta
    \] for all $\delta \in (0,1).$
    \end{asmp} Under the above assumption, \eqref{eqn:sve} admits a unique $\cF^B$-adapted solution $X$ such that $\bbE |X|_\infty^p < \infty$ for any $p > 0$ by \cite[Theorem $1.1$]{jourdainpages2025convex} or \cite[Theorem $3.1$]{zhang2010volterra}.

\begin{thm}
    \label{thm:adaptedQTCIforSVEs}
    Suppose $b,\sigma,G$ satisfy Assumption \ref{asmp:svecoeff}, $\sigma > 0$ is bounded and $G$ admits a resolvent of the first kind: there exists a finite signed measure $\mu$ on $[0,1]$ such that for all $t \in (0,1]$, \[
        \int_{[0,t]} G(t - s) \mu(ds) = 1.
        \]
    Then, the law $\kappa^{\sigma,b}_G$ of the solution to \eqref{eqn:sve} satisfies an AQTCI w.r.t. $d = \sqrt{\int_0^1 |\cdot|^2 ds}$ with constant $C := 4D \exp(\frac{2^{1+ q/\beta}A^{1 + q/\beta} B^q}{q \beta^{q/\beta}}) > 0$, where $q = \frac{\beta}{\beta - 1}$, $A := 6K^2, B := (\sup_{t \in [0,1]} \int_0^t G(t - u)^{2\beta} du)^{1/\beta},D:= 3 \Vert \sigma \Vert_{\infty}^2  B$.
\end{thm}

\begin{proof}
    We fix $\bbP \ll \kappa_G^{\sigma,b}$ such that $\cH(\bbP | \kappa_G^{\sigma,b}) < \infty$. 

By Proposition \ref{prop:svebcbmap}, the solution map $F^{\sigma,b}_G$ lies in $\cT_{bcb}(\bbW^1, \kappa_G^{\sigma,b})$. Let $\tilde{\bbP} := (F^{\sigma,b}_G)^{-1}_{\#}\bbP.$ By Corollary \ref{cor:optimalbicausalcpling}, there exists some real-valued $\cF^X$-adapted process $v = (v_t)_{t \in [0,1]}$ such that the Girsanov shift $V := V^{F^{\sigma,b}_G}$ defined by \eqref{eqn:girsanovshift} satisfies $V \in \cT_{bc}(\tilde{\bbP}, \bbW^1)$, $(F^{\sigma,b}_G \circ V, F^{\sigma,b}_G)_{\#}\tilde{\bbP} \in \Pi_{bc}(\kappa^{\sigma,b}_G,\bbP)$ and \[
2\cH(\bbP | \kappa^{\sigma,b}_G) = \bbE^{\tilde{\bbP}} \int_0^1 |v(s,X)|^2 ds.
\] Since $\tilde{\bbP} \ll \bbW^1$, $Z := F^{\sigma,b}_G(X)$ satisfies $\bbW^1$-a.s. and hence, $\tilde{\bbP}$-a.s. \[
Z_t = \int_0^t G(t - s) b(s,Z_s) ds + \int_0^t G(t - s) \sigma(s,Z_s) dX_s.
\] Furthermore, since $V_{\#}\tilde{\bbP} = \bbW^1$, i.e. $V$ is a Brownian motion, then $\tilde{Z} := F^{\sigma,b}_G \circ V$ satisfies $\tilde{\bbP}$-a.s. \[
\tZ_t = \int_0^t G(t - s) b(s, \tZ_s) ds + \int_0^t G(t -s) \sigma(s,\tZ_s) dV_s.
\] Since $X - V = \int_0^\cdot v(s,X) ds$, the argument in \cite[Appendix $C$, Property $2$]{jourdainpages2025convex} with minor modifications, gives $\sup_{t \in [0,1]} \bbE^{\tilde{\bbP}} |Z_t|^2 + |\tZ_t|^2 < \infty$. We then use H\"older's inequality, Cauchy-Schwarz inequality and It\^o's isometry to obtain \begin{align*}
    &\bbE^{\tilde{\bbP}} |Z_s - \tZ_s|^2 \\
    \leq & 3(\bbE^{\tilde{\bbP}} | \int_0^s G(s- u)(b(u,Z_u) - b(u, \tZ_u)) du|^2 + \bbE^{\tilde{\bbP}} | \int_0^s G(s- u)( \sigma(u,Z_u) - \sigma(u,\tZ_u)) dV_u|^2)\\
    &+ 3\bbE^{\tilde{\bbP}} |\int_0^s G(s- u) v(u,X) \sigma(u,Z_u) du|^2 \\
    \leq &6 K^2 \bbE^{\tilde{\bbP}}\int_0^s G(s- u)^2 |Z_u - \tZ_u|^2  du + 3 \Vert \sigma \Vert_{\infty}^2 \bbE^{\tilde{\bbP}} (\int_0^s G(s- u) |v(u,X)| du)^2 \\
    \leq &6 K^2 \bbE^{\tilde{\bbP}}\int_0^s G(s- u)^2 |Z_u - \tZ_u|^2  du + 3 \Vert \sigma \Vert_{\infty}^2  (\int_0^s G(s-u)^2 du) \bbE^{\tilde{\bbP}} (\int_0^s v(u,X)^2 du) \\
    \leq &6 K^2 \bbE^{\tilde{\bbP}}\int_0^s G(s- u)^2 |Z_u - \tZ_u|^2  du + 6 \Vert \sigma \Vert_{\infty}^2  (\sup_{t \in [0,1]} \int_0^t G(t - u)^{2\beta} du)^{1/\beta} \cH(\bbP|\kappa^{\sigma,b}_G).
\end{align*} Temporarily denoting $f(s) := \bbE^{\tilde{\bbP}} |Z_s - \tZ_s|^2$, the above gives \[
f(s) \leq A \int_0^s G(s- u)^2 f(u) du + 2D\cH(\bbP|\kappa^{\sigma,b}_G).
\] Let $q := \frac{\beta}{\beta - 1}$ be the H\"older conjugate of $\beta$. Then, writing $f^*(s) := \sup_{u \in [0,s]} f(u)$, we have $f^*(1) < \infty$ and H\"older's inequality yields \[
\int_0^s G(s - u)^2 f(u) du \leq (\int_0^s G(s - u)^{2\beta} du)^{1/\beta} (\int_0^s f(u)^q du)^{1/q} 
    \leq B (\int_0^s f(u) du)^{1/q} f^*(s)^{\frac{q-1}{q}},
\] where we have used that $f(u)^q \leq f(u) f^*(s)^{q-1}$. Fix now some $\epsilon > 0$. Applying Young's inequality for the conjugate pair $(\frac{q}{q-1},q) = (\beta,q)$ gives \[
f^*(s)^{\frac{q-1}{q}} B (\int_0^s f(u) du)^{1/q} \leq \frac{\epsilon^\beta}{\beta} f^*(s) + \frac{B^q}{q\epsilon^q} \int_0^s f(u) du.
\] Hence, we obtain \[
f(s) \leq \frac{A \epsilon^\beta}{\beta} f^*(s) + \frac{A B^q}{q \epsilon^q} \int_0^s f(u) du + 2D\cH(\bbP|\kappa^{\sigma,b}_G).
\] and for any $0 \leq s \leq t \leq 1$, \[
f(s) \leq \frac{A \epsilon^\beta}{\beta} f^*(t) + \frac{A B^q}{q \epsilon^q} \int_0^t f^*(u) du + 2D\cH(\bbP|\kappa^{\sigma,b}_G).
\] Taking supremum over $s \in [0,t]$ yields the inequality \[
f^*(t) \leq \frac{A \epsilon^\beta}{\beta} f^*(t) + \frac{A B^q}{q \epsilon^q} \int_0^t f^*(u) du + 2D\cH(\bbP|\kappa^{\sigma,b}_G).
\] By choosing $\epsilon := (\frac{\beta}{2A})^{1/\beta}$, we obtain $1 - \frac{A \epsilon^\beta}{\beta} = \frac{1}{2}$ and \[
f^*(t) \leq \frac{2^{1+q/\beta}A^{1 + q/\beta} B^q}{q \beta^{q/\beta}} \int_0^t f^*(u) du + 4D\cH(\bbP|\kappa^{\sigma,b}_G).
\] Gronwall's inequality now yields \[
f^*(1) = \sup_{s \in [0,1]} \bbE^{\tilde{\bbP}} |Z_s - \tZ_s|^2 \leq 4D \exp(\frac{2^{1+q/\beta}A^{1 + q/\beta} B^q}{q \beta^{q/\beta}}) \cH(\bbP|\kappa^{\sigma,b}_G).
\] Since $\bbE^{\tilde{\bbP}} \int_0^1 |Z_s - \tZ_s|^2 ds \leq f^*(1)$ and $\Law( \tZ,Z) = (F^{\sigma,b}_G \circ V, F^{\sigma,b}_G)_{\#}\tilde{\bbP} \in \Pi_{bc}(\kappa^{\sigma,b}_G,\bbP)$, we conclude \[
\cA\cW_{\sqrt{\int_0^1 | \cdot|^2 ds}, 2}^2(\bbP, \kappa^{\sigma,b}_G) \leq \bbE^{\tilde{\bbP}} \int_0^1 |Z_s - \tZ_s|^2 ds \leq 4D \exp(\frac{2^{1+q/\beta}A^{1 + q/\beta} B^q}{q \beta^{q/\beta}}) \cH(\bbP|\kappa^{\sigma,b}_G)
\] as required.
\end{proof}

\subsubsection*{Acknowledgments} The authors would like to thank Dr. Yifan Jiang  for the numerous fruitful discussions, helpful comments and for the computations used in Proposition \ref{prop:svebcbmap}. Erhan Bayraktar is supported in part by National Science Foundation Grant DMS-2602036 and in part by the Susan M.
Smith Chair.


\appendix
\section{}
\label{section:appendix}

\subsection{Martingale representation property under finite relative entropy}
\label{appendix:mrpfiniterelent}

Let $\bbQ := \mu^{\sigma,b} \in \cS \cM(\cW^n)$ be the law of an SDE with a unique (in law) weak solution and $\bbP \in \cP(\cW^n)$ be such that $\cH(\bbP|\bbQ) < \infty$. 

We first recall Girsanov's theorem established in \cite{leonard2012girsanov} specialized to our context:
\begin{lem}[{\cite[Theorem $1\& 2$]{leonard2012girsanov}}]
    \label{thm:cleonardgirsanov}
    For any $\bbP \in \cP(\cW^n)$ such that $\cH(\bbP|\bbQ) < \infty$, there exists an $\bbR^n$-valued $\cF^X$-adapted process $v : [0,1] \times \cW^n \rightarrow \bbR^n$ defined $\bbP$-a.s. such that \[
    \bbE^\bbP \int_0^1 |\sigma^*(s,X) v(s,X)|^2 ds = 2 \cH(\bbP| \mu^{\sigma,b}) 
    \] and under $\bbP$, the canonical process $X$ has $\cF^X$-semimartingale decomposition $X = V^\bbP + M^\bbP$, where $V^\bbP$ (resp. $M^\bbP$) is a finite-variation process (resp. continuous local martingale) satisfying \begin{equation}
        \label{eqn:semimartdecompunderfiniterelent}
        V^\bbP_t = \int_0^t b(s,X) + \sigma \sigma^*(s,X) v(s,X) ds \text{ and } [M^\bbP]_t = \int_0^t \sigma(s,X) \sigma^*(s,X) ds.
    \end{equation} Furthermore, since $\bbQ$ is the law of an SDE with a unique weak solution, then $\bbP$ is also unique in the following sense: if $\bbP' \in \cP(\cW^n)$ is such that $\cH(\bbP'|\bbQ) <\infty$ and under $\bbP'$, the canonical process $X$ has the same semimartingale decomposition as $\bbP$, then $\bbP = \bbP'$.
\end{lem}

This theorem and the following lemmata lead to the MRP for $\bbP$ in Proposition \ref{prop:mrpunderfiniterelent}.
For any $\bbP \in \cP(\cW^n)$ such that $\cH(\bbP|\bbQ) < \infty$, the above gives that $\bbP \in \cS \cM(\cW^n)$ and we can define the set of stochastic integrals \[
    \Gamma^\bbP := \{ N \in L^2(\bbP) : N = \bbE^\bbP(N) + (\alpha \cdot M^{\bbP})_1, \alpha \in K \},
    \] where $K$ is the set of all $\bbR^n$-valued $\cF^X$-progressively measurable processes $\alpha : [0,1] \times \cW^n \rightarrow \bbR^n$ such that \[
    \bbE^\bbP \int_0^1 |\sigma^*(s,X) \alpha(s,X)|^2 ds < \infty.
    \]

\begin{lem}
    \label{lem:orthogonalityreduction}
    Denoting $(\Gamma^\bbP)^\perp$ as the orthogonal complement of $\Gamma^\bbP$ in $L^2(\bbP)$, we have \[
    (\Gamma^\bbP)^\perp \cap L^\infty(\bbP) = \{ 0 \} \implies (\Gamma^\bbP)^\perp \cap L^2(\bbP) = \{ 0 \}.
    \]
\end{lem}

\begin{proof}
    As preparation, we note that since $\bbQ$ is assumed to be the law of an SDE with weak uniqueness, it is trivially an extremal solution to its martingale problem and \cite[Theorem $2.3$]{stroockyor1980extremal} gives that $\cF^X_0 = \cap_{\epsilon > 0} \sigma(X_s : s \leq \epsilon)$ is trivial modulo $\bbQ$-null sets. In particular, any $\cF^X_0$-measurable random variable must be $\bbQ$-a.s. and hence, $\bbP$-a.s. constant.

    We prove the contrapositive. Assume that $0 \neq M_1 \in (\Gamma^\bbP)^\perp \cap L^2(\bbP)$. Since $\Gamma^\bbP$ contains the constants, we conclude $M_0 = \bbE^\bbP M_1 = 0$. There exists a decomposition of the martingale $M := (M_t)_{t \in [0,1]}$, where $M_t := \bbE^\bbP(M_1 | \cF^X_t)$, into continuous and purely discontinuous martingale parts, i.e. $M = M^c + M^d$ such that $M^c_0 = M^d_0 = 0$. In particular, for any $N \in \Gamma^\bbP$, we compute that \[
    \bbE^\bbP (N M_1) = 0 \implies \bbE^\bbP (N M^c_1) = 0,
    \] since, using that $\bbE^\bbP M^d_1 = 0$ and $M^d$ is purely discontinuous, we obtain $\bbE^\bbP(N M^d_1) = 0$. Since $N$ was arbitrary in $\Gamma^\bbP$, $M^c_1 \in (\Gamma^\bbP)^\perp$. By introducing stopping time $\tau_\epsilon := \inf \{t \in [0,1] : |M^c_t| \geq \epsilon \}$, we obtain that if $M^c \neq 0$, then $0 \neq M^c_{\tau_\epsilon} \in (\Gamma^\bbP)^\perp \cap L^\infty(\bbP)$ for some $\epsilon > 0$ small enough. Hence, we now assume that $M^c = 0$ and $M = M^d \neq 0.$ In this case, it suffices to exhibit one non-zero purely discontinuous, bounded martingale $\tilde{M}$, so that $\tilde{M}_1$ is bounded and automatically orthogonal to $\Gamma^\bbP.$ This follows immediately from \cite[Theorem $5.125$]{vandervaartstochasticcalculus}, but we provide an alternative argument here. By \cite[Theorem $6.22.2$]{hewangyan1992semimart}, there exists a sequence of stopping times $(T_n)_{n \geq 0}$, each either predictable or totally inaccessible, exhausting the jumps of $M$ such that $M^d = \sum_n M^n$, where $M^n = A^n - \tilde{A^n}$, $A^n := \Delta M_{T_n} \mathbbm{1}_{[T_n, \infty)}$ and $\tilde{A^n}$ its compensator. Fix any purely discontinuous $M^n \neq 0$ and let $\xi^{n,m} :=  \max\{ - m, \min\{ \Delta M_{T_n}, m\}\}$. Note that since $M^n \neq 0$, then $\Delta M_{T_n} \neq 0$ and for $m$ large enough, we must have $\xi^{n,m} \neq 0$. For any such $m$, define $A^{n,m}:= \xi^{n,m} \mathbbm{1}_{[T_n,\infty)}$ with compensator $\tilde{A^{n,m}}$ so that by \cite[Theorem $6.21.1$]{hewangyan1992semimart}, $M^{n,m} := A^{n,m} - \tilde{A^{n,m}}$ is purely discontinuous. 

    For $T_n$ totally inaccessible, $A^{n,m}$ is quasi-left-continuous \cite[Definition $4.22$]{hewangyan1992semimart} and hence, \cite[Corollary $5.28.3$]{hewangyan1992semimart} gives that $\tilde{A^{n,m}}$ is continuous. Then, introducing another stopping time $\tau'_\epsilon := \inf\{ t \in [0,1] : | M^{n,m}_t| \geq \epsilon \}$, we have, for $\epsilon
     > 0$ small enough, that $M^{n,m}_{\cdot \wedge \tau'_\epsilon}$ is a non-zero bounded, purely discontinuous martingale since the space of purely discontinuous martingales is stable \cite[Corollary $6.17$]{hewangyan1992semimart} and \[
    |M^{n,m}_{\tau'_\epsilon}| \leq |M^{n,m}_{\tau'_\epsilon -}| + |\Delta M^{n,m}_{\tau'_\epsilon}| \leq \epsilon + |\Delta A^{n,m}| \leq \epsilon + m.
    \] 

    For $T_n$ predictable, $\tilde{A^{n,m}} = \bbE^\bbP(\xi^{n,m}|\cF_{T_n -}) \mathbbm{1}_{[T_n, \infty)}$ by \cite[Theorem $5.29.2$]{hewangyan1992semimart}. Hence, $|M^{n,m}| \leq 2m$ and $M^{n,m}$ is a purely discontinuous, bounded martingale. It suffices now to verify that $M^{n,m} \neq 0$. If $M^{n,m} = 0$ for all sufficiently large $m$, then we must have $\xi^{n,m} = \bbE^\bbP(\xi^{n,m}|\cF_{T_n -})$ and hence, taking limits as $m \rightarrow \infty$, we conclude that $\Delta M_{T_n}$ is $\cF_{T_n -}$-measurable. \cite[Theorem $4.41$]{hewangyan1992semimart} gives $\Delta M_{T_n} = \bbE^\bbP (\Delta M_{T_n} | \cF_{T_n - }) = 0$, which is a contradiction. Hence, there exists some $m$ large enough such that $M^{n,m}$ is a non-zero purely discontinuous, bounded martingale. 
\end{proof}

\begin{lem}
    \label{lem:densityofstoint}
    For any $\bbP \in \cP(\cW^n)$ such that $\cH(\bbP|\bbQ) < \infty$, $\Gamma^{\bbP}$ is dense in $L^2(\bbP)$.
\end{lem}

\begin{proof}
    The proof method follows closely that of \cite[Theorem $1$]{ustunel2019mrp}. By Lemma \ref{lem:orthogonalityreduction}, it suffices to show that $(\Gamma^\bbP)^\perp \cap L^\infty(\bbP) = \{ 0 \}$. Suppose for contradiction that there exists some $M_1 \in L^\infty(\bbP) \backslash \{ 0 \}$ orthogonal to $\Gamma^\bbP$. In particular,  $\bbE^\bbP M_1 = 0$ since $\Gamma^\bbP$ contains the constants and there exists some $\epsilon > 0$ such that $D_1 := 1+ \epsilon M_1$ is strictly positive, $\bbE^\bbP D_1 = 1$ and $D_1$ is orthogonal to $\Gamma^\bbP_0$, where \[
    \Gamma^{\bbP}_0 := \{ N \in L^2(\bbP) : N = (\alpha \cdot M^{\bbP})_1, \alpha \in K \}.\] We claim that $D_1 = 1$ $\bbP$-a.s.. In this case, we must have $M_1 = 0$ $\bbP$-a.s., which is clearly a contradiction. 
    
    Fix any smooth function $f : \bbR^n \rightarrow \bbR$ with bounded derivatives and $Y = (Y_t)_{t \in [0,1]}$, where $Y_t := f(X_t) - f(X_0) - \int_0^t Lf(s,X) ds$ and \[
    Lf(t,X) = \frac{1}{2} \sum_{i,j} (\sigma \sigma^*)^{i,j}(t,X) \partial_{i,j} f(X_t) + \sum_{i} (b^i(t,X) + (\sigma \sigma^*(t,X)v(t,X))^i)\partial_i f(X_t).
    \] Under $\bbP$, by using It\^o's formula and the semimartingale decomposition of $X$ in \eqref{eqn:semimartdecompunderfiniterelent}, $Y = \nabla f (X_\cdot) \cdot M^\bbP$ is a continuous local $\cF^X$-martingale. Consider now $\bbP' \in \cP(\cW^n)$ defined by $d\bbP' := D_1 d\bbP.$ We claim that $Y$ is again a continuous local $\cF^X$-martingale under $\bbP'$. Indeed, up to the usual stopping techniques, we can assume that $Y$ is bounded. Fixing $0 \leq s < t \leq 1$ and some $\cF^X_s$-measurable bounded function $G$, we compute \[
    \bbE^{\bbP'}[G(Y_t - Y_s)] = \bbE^{\bbP'}[ \int_0^1 G \mathbbm{1}_{[s,t)}(u) dY_u] = \bbE^\bbP[ D_1 \sum_{i = 1}^n \int_0^1 G \mathbbm{1}_{[s,t)}(u) \partial_i f(X_u) d(M^\bbP)^i_u ] = 0,
    \] since $(u,\omega) \mapsto G(\omega) \mathbbm{1}_{[s,t)}(u) \nabla f(X_u(\omega)) \in K$ and $D_1$ is orthogonal to $\Gamma^\bbP_0$. As in the proof of \cite[Proposition $2.1$]{ikedawatanabe1989stochasticde}, by varying $f$ over all smooth functions with bounded derivatives, one can conclude that $X$ is again a semimartingale under $\bbP'$ with the same semimartingale decomposition as $X$ under $\bbP$. Finally, we note that $\frac{d\bbP'}{d\bbQ} = D_1 \frac{d\bbP}{d\bbQ}$ and compute \[
    \cH(\bbP'|\bbQ) = \bbE^{\bbP}[D_1\log(D_1)] + \bbE^\bbP[D_1 \log(\frac{d\bbP}{d\bbQ})].
    \] Since $D_1$ is bounded, the first term on the RHS is finite. To see that the second term is finite, note that $x \mapsto x \log(x)$ is bounded on $(0, 1]$ so that \[
    \bbE^\bbP \log(\frac{d\bbP}{d\bbQ}) \mathbbm{1}_{\{\frac{d\bbP}{d\bbQ} \geq 1\}} = \cH(\bbP|\bbQ) - \bbE^\bbP \log(\frac{d\bbP}{d\bbQ}) \mathbbm{1}_{\{\frac{d\bbP}{d\bbQ} < 1\}} < \infty. 
    \]  Hence, using the decomposition \[
    \bbE^\bbP (D_1 \log(\frac{d\bbP}{d\bbQ})) = \bbE^\bbP D_1\log(\frac{d\bbP}{d\bbQ}) \mathbbm{1}_{\{\frac{d\bbP}{d\bbQ} \geq 1\}} + \bbE^\bbQ D_1 \frac{d\bbP}{d\bbQ} \log(\frac{d\bbP}{d\bbQ}) \mathbbm{1}_{\{\frac{d\bbP}{d\bbQ} < 1\}}
    \] and that $D_1$ is non-negative and bounded, we can conclude $\bbE^\bbP (D_1 \log(\frac{d\bbP}{d\bbQ})) < \infty$. In particular, by the uniqueness statement in Lemma \ref{thm:cleonardgirsanov}, we have $\bbP' = \bbP$, that is, $D_1 = 1$ $\bbP$-a.s..
\end{proof}

\begin{prop}
    \label{prop:mrpunderfiniterelent}
    On $(\cW^n, \cF^X_1, \cF^X, \bbP)$, we have that for any $F \in L^2(\bbP)$, there exists an $(\cF^X_t)_{t \in [0,1]}$-adapted process $\xi : [0,1] \times \cW^n \rightarrow \bbR^n$ such that \[
    \bbE^\bbP \left ( \int_0^1 |\xi(s,X)|^2 ds\right ) < \infty
    \] and \[
    F = \bbE^\bbP F +  (\xi \cdot L)_1 \quad \bbP\text{-}a.s.,
    \] where $L_\cdot = \int_0^\cdot \sigma^\dag(s,X) dM^\bbP_s$ and $A^\dag$ denotes the Moore-Penrose pseudoinverse of matrix $A$. In particular, under $\bbP$, any square integrable $\cF^X$-martingale $(M_t)_{t \in [0,1]}$ has representation \[
    M_t = M_0 + (\xi \cdot L)_t \quad \forall t \in [0,1] \quad \bbP\text{-}a.s..
    \]
\end{prop}

\begin{proof}
Fix $F \in L^2(\bbP)$ and assume, w.l.o.g., that $\bbE^\bbP F = 0$. Then, by Lemma \ref{lem:densityofstoint}, there exists a sequence $(\alpha^n)_{n \geq 0} \subseteq K$ such that \[
\Vert (\alpha^n \cdot M^\bbP)_1 - F \Vert_{L^2(\bbP)} \rightarrow 0 \text{ as } n \rightarrow \infty.
\] Note that by It\^o's isometry, we have \[
\Vert ( (\alpha^n - \alpha^m) \cdot M^\bbP)_1 \Vert_{L^2(\bbP)} = \Vert \sigma^* (\alpha^n - \alpha^m)\Vert_{L^2(\bbP \otimes \Leb)}, 
\] so that the sequence $(\sigma^* \alpha^n)_{n \geq 0}$ is a Cauchy sequence of $\cF^X$-adapted $L^2(\bbP \otimes \Leb)$-integrable random processes. Hence, there exists some $\xi$ as in the proposition statement satisfying \[
\Vert \sigma^* \alpha^n - \xi \Vert_{L^2(\bbP \otimes \Leb)} \rightarrow 0 \text{ as } n \rightarrow \infty.
\] In fact, since $\sigma^\dag \sigma = ( \sigma^\dag \sigma)^* = \sigma^* (\sigma^\dag)^* = \sigma^* (\sigma^*)^\dag$ is an orthogonal projection onto the orthogonal complement of the kernel of $\sigma$ (see \cite{golubvanloan1996matrix}), we have $|\sigma^\dag \sigma x | \leq |x|$ for any $x \in \bbR^n$ and hence, we have that \[
\Vert \sigma^* \alpha^n - \sigma^* (\sigma^\dag)^* \xi \Vert_{L^2(\bbP \otimes \Leb)} = \Vert \sigma^* (\sigma^\dag)^* \sigma^* \alpha^n - \sigma^* (\sigma^\dag)^* \xi \Vert_{L^2(\bbP \otimes \Leb)} \rightarrow 0 \text{ as } n \rightarrow \infty.
\] Here, we have used the property of the Moore-Penrose pseudoinverse that $\sigma^* (\sigma^\dag)^* \sigma^* = \sigma^*$. We compute by the triangle inequality \[
\Vert F - (\xi \cdot L)_1 \Vert_{L^2(\bbP)} \leq \limsup_n \Vert F - (\alpha^n \cdot M^\bbP)_1 \Vert_{L^2(\bbP)} + \limsup_n \Vert (\alpha^n \cdot M^\bbP)_1 - (\xi \cdot L)_1 \Vert_{L^2(\bbP)}.
\] The first limit superior on the RHS is $0$. Since $(\xi \cdot L)_1 = ((\sigma^\dag)^* \xi \cdot M^{\bbP})_1$, by It\^o's isometry, we have \[
\Vert (\alpha^n \cdot M^\bbP)_1 - (\xi \cdot L)_1 \Vert_{L^2(\bbP)} = \Vert \sigma^* (\alpha^n - (\sigma^\dag)^* \xi) \Vert_{L^2(\bbP \otimes \Leb)},
\] from which we conclude that the second limit superior is also $0$. In particular, $F = (\xi \cdot L)_1$ $\bbP$-a.s.. The second assertion is straightforward.
\end{proof}

\subsection{Bicausal bijective maps and solutions to stochastic equations}
\label{appendix:bcbmaps}

We begin with an abstract lemma:

\begin{lem}
    \label{lem:bcbmapsequiv}
    Fix $\bbP,\bbQ \in \cP(\cW^n)$. Suppose on some abstract probability space $(\Omega, \cG, \bbP^*)$, there exist two stochastic processes $Z = (Z_t)_{t\in [0,1]}$ and $\tZ = (\tZ_t)_{t\in [0,1]}$ with continuous sample paths $\bbP^*$-a.s. such that $\Law(Z) = \bbP, \Law(\tZ) = \bbQ$. Then, $\tZ = T(Z)$ $\bbP^*$-a.s. for some $T \in \cT_{bcb}(\bbP,\bbQ)$ if and only if $\cF^Z \equiv \cF^{\tZ}$ $\bbP^*$-a.s.. 
\end{lem}

\begin{proof}
    Since $T$ is causal, we obtain that $\cF^{\tZ} \subseteq \cF^Z$ $\bbP^*$-a.s. and since $S:= T^{-1}$ is also causal satisfying $S(\tZ) = S(T(Z)) = Z$ $\bbP^*$-a.s., we obtain $\cF^Z \subseteq \cF^{\tZ}$ $\bbP^*$-a.s.. The converse follows from Doob's factorization lemma.
\end{proof}

We now proceed to verify under various conditions that solutions of stochastic equations generate the same filtration as their driving noise. For the rest of the section, we fix a Brownian motion $B = (B_t)_{t \in [0,1]}$ on some abstract probability space $(\Omega, \cG, \bbP^*)$.

\subsubsection*{Brownian noise}

\begin{prop}
    \label{prop:sdebcbmap}
    Suppose $\sigma \in \scrA^{n,n}, b \in \scrA^{n,1}$ are such that $\det(\sigma) \neq 0$ and there exists a pathwise unique, strong solution to \eqref{eqn:sde}. Let $X= (X_t)_{t \in [0,1]}$ be any solution of \eqref{eqn:sde} driven by Brownian motion $B = (B_t)_{t \in [0,1]}$. Then, $\cF^X \equiv \cF^B$ $\bbP^*$-a.s. and the It\^o solution map $F^{\sigma,b} $ lies in $ \cT_{bcb}(\bbW^n, \mu^{\sigma,b}).$
\end{prop}

\begin{proof}
    By definition of a pathwise unique, strong solution, we have that $F^{\sigma,b}(B) = X$ and $\cF^X \subseteq \cF^B$ $\bbP^*$-a.s.. To see the converse, we note that \[
    B_t = \int_0^t \sigma^{-1}(s,X_{\cdot \wedge s}) dM_s,
    \] where $M_t = X_t - \int_0^t b(s,X_{\cdot \wedge s}) ds$, $t \in [0,1]$, from which it is clear that $\cF^B \subseteq \cF^X$ $\bbP^*$-a.s..
\end{proof}

\subsubsection*{Fractional Brownian noise}
We restrict ourselves to the scalar case ($n = 1$) and let $B^H = (B^H_t)_{t \in [0,1]}$ be a fractional Brownian motion (fBM) with Hurst parameter $H \in (0,1)$. Consider the Molchan-Golosov representation of fBM: \[
B^H_t = \int_0^t K_H(t,s) dB_s, \quad t \in [0,1],
\] where $K_H \in L^2([0,1]^2)$ is the Molchan-Golosov kernel given by \[
K_H(t,s) := c_H (t - s)^{H-1/2}F(H-1/2,1/2 - H, H + 1/2, 1 - t/s),
\] 
where $c_H > 0$ is some normalizing constant and $F$ is the hypergeometric function.  In this case, $\cF^{B^H} \equiv \cF^B$ $\bbP^*$-a.s.. We refer readers to \cite{decreusefondustunel1999} for more details. 
We first investigate the case of additive noise, whose solution theory is straightforward:

\begin{prop}
    \label{prop:additivegennoise}
    For any $b \in \scrA^{1,1}$ satisfying \eqref{eqn:lipschitzsdecoeff} and $\int_0^1 |b(t,0)| dt < \infty$, there exists a unique solution $X = (X_t)_{t \in [0,1]}$ to \eqref{eqn:additivede} such that $\cF^X \equiv \cF^B$ $\bbP^*$-a.s.. In particular, the solution map $\mathfrak{F}^{1,b} : \cW^1 \rightarrow \cW^1$ sending $B$ to $X$ lies in $\cT_{bcb}(\bbW^1, \varrho^{1,b})$.
\end{prop}

\begin{proof}
    The existence of a process $X = (X_t)_{t \in [0,1]}$ satisfying \eqref{eqn:additivede} follows directly from Picard's iteration, from which $\cF^X \subseteq \cF^{B^H} \equiv \cF^B$ $\bbP^*$-a.s. follows. The converse $\cF^{B^H} \equiv \cF^B \subseteq \cF^X$ $\bbP^*$-a.s. follows easily since \[
    B^H_t = X_t - \int_0^t b(s,X_{\cdot \wedge s}) ds.
    \]
\end{proof}

To handle multiplicative noise, we adopt Doss-Sussmann's approach \citep{doss1977liens,sussmann1978}. Fix $b ,\sigma : \bbR \rightarrow \bbR$ and $g \in \cW^1$. 

\begin{defn}
    We say that $x \in \cW^1$ is a solution to \begin{equation}
    \label{eqn:dosssussmande}
    dx_t = b(x_t) dt + \sigma(x_t) dg_t, \quad x_0 = 0
\end{equation} if there exists a continuous map $\Gamma^{\sigma,b} : \cW^1 \rightarrow \cW^1$ such that: \begin{enumerate}
    \item for each differentiable $f \in \cW^1$, the continuous function $\Gamma^{\sigma,b}(f)$ is a solution to the ordinary differential equation (ODE) \begin{equation}
        \label{eqn:dosssussmanode}
        \frac{dx}{dt} = b(x_t) + \sigma(x_t) \frac{df}{dt}, x_0 = 0,
    \end{equation}
    \item $\Gamma^{\sigma,b}(g) = x$.
\end{enumerate}
\end{defn}

When $b$ is locally Lipschitz with linear growth and $\sigma \in C^1(\bbR)$ with locally Lipschitz and bounded derivative, there exists a unique solution by \cite[Theorem $2$]{sussmann1978}. 

In the scalar setting, $\Gamma^{\sigma,b}$ has an explicit form: Let $h : \bbR \times \bbR \rightarrow \bbR$ be such that for each $\alpha \in \bbR$, $\beta \mapsto h(\alpha,\beta)$ is the solution to the ODE given by \[
\frac{dx}{dt} = \sigma(x_t), \quad x_0 = \alpha.
\] In particular, we have $\frac{\partial h}{\partial \beta}(\alpha,\beta) = \sigma(h(\alpha,\beta))$ and $\frac{\partial h}{\partial \alpha}(\alpha,\beta) = \exp(\int_0^\beta \sigma'(h(\alpha,u)) du)$. We further fix $f(x,y) :=\exp(-\int_0^y \sigma'(h(x,s)) ds) b(h(x,y))$ and let $\gamma \in \cW^1$ be the solution to the ODE given by \[
\frac{dx}{dt} = f(x_t,g_t), \quad x_0 = 0.
\] Then, we have \[
\Gamma^{\sigma,b}(g) = h(\gamma_{\cdot},g_{\cdot}).
\]

\begin{prop}
    \label{prop:dosssussmanbcbmap}
    Suppose $\sigma \in C^2_b(\bbR)$, $b \in C^1_b(\bbR)$ and $\sigma > \delta$ for some $\delta > 0$. Then, setting $X := \Gamma^{\sigma,b}(B^H)$ and $ Y := \Gamma^{1,\tb}(B^H)$ for \[
    \tb(x) := \frac{b(h(0,x))}{\sigma(h(0,x))},
    \] we have $X_t = h(0, Y_t)$. In particular, $\cF^X \equiv \cF^{B}$ $\bbP^*$-a.s. and $\Gamma^{\sigma,b} \circ \cK_H \in \cT_{bcb}(\bbW^1, \varrho^{\sigma,b}),$ where $\cK_H : \cW^1 \rightarrow \cW^1$ is the map satisfying $\bbP^*$-a.s. \[
    \cK_H(B) = \int_0^\cdot K_H(\cdot,s)dB_s.
    \]
\end{prop}

\begin{proof}
    We first note that for $f \in \cW^1$ differentiable, if $x := \Gamma^{\sigma, b}(f)$ and $y := \Gamma^{1,\tb}(f)$, then standard calculus gives $ x_t = h(0,y_t)$. Set $B^{H,\epsilon}_t := \frac{1}{2\epsilon} \int_0^t B^H_{(s + \epsilon) \wedge 1} - B^H_{(s - \epsilon) \vee 0} ds$ so that $B^{H,\epsilon} \in C^1([0,1])$ a.s. and $|B^{H,\epsilon} - B^H|_\infty \rightarrow 0$ a.s.. Then, for $X^{\epsilon} := \Gamma^{\sigma,b}(B^{H,\epsilon})$ and $Y^\epsilon := \Gamma^{1,\tb}(B^{H, \epsilon})$, \[
    X^\epsilon_t = h(0, Y^\epsilon_t) \stackrel{\epsilon \rightarrow 0}{\implies} X_t = h(0, Y_t),
    \] by continuity of the maps $\Gamma^{\sigma,b}$ and $\Gamma^{1, \tb}$. Since $h(0,\cdot)$ is invertible with inverse $h^{-1}(0,x) = \int_0^x \frac{1}{\sigma(s)} ds$, it is clear that $\cF^X \equiv \cF^Y$ $\bbP^*$-a.s.. By Proposition \ref{prop:additivegennoise}, $\cF^Y \equiv \cF^B$ $\bbP^*$-a.s., from which we conclude $\cF^X \equiv \cF^B$ $\bbP^*$-a.s..
\end{proof}

\begin{rmk}
    \label{rmk:dosssussmansolutions}
    A priori, it is unclear how to interpret the ``integral" $\int_0^t \sigma(x_t) dg_t := x_t - \int_0^t b(x_s) ds$ for arbitrary $g \in \cW^1$. However, in the case that $g$ is given by Brownian motion $B = (B_t)_{t \in [0,1]}$, then $X:=\Gamma^{\sigma,b}(B)$ coincides with the solution to the \emph{Stratonovich} SDE: \[
    X_t = \int_0^t b(X_s) ds + \int_0^t \sigma(X_s) \circ dB_s, \quad X_0 = 0,
    \] so that $X_t - \int_0^t b(X_s) ds = \int_0^t \sigma(X_s) \circ dB_s$ is the Stratonovich integral. More generally, setting $X := \Gamma^{\sigma,b}(B^H)$, we have \begin{enumerate}
        \item for $H \in (1/2,1)$, $X$ is a solution in the sense studied in \cite{nualartrascanu2002differentialeq} and $X_t - \int_0^t b(X_s) ds = \int_0^t \sigma(X_s)dB^H_s$ is a Riemann-Stieltjes integral, and 
        \item for $H \in (1/4,1/2)$, $X$ is a solution in the sense of \cite[Definition $5$]{alosleonnualart2001stochastic} and $X_t - \int_0^t b(X_s) ds = \int_0^t \sigma(X_s) \circ dB^H_s$ is the Stratonovich integral defined in \cite[Definition $1$]{alosleonnualart2001stochastic}.
    \end{enumerate}
\end{rmk}


\subsubsection*{SVEs}

    \begin{prop}
        \label{prop:svebcbmap}
        Under the assumptions of Theorem \ref{thm:adaptedQTCIforSVEs}, the unique solution $X = (X_t)_{t \in [0,1]}$ to \eqref{eqn:sve} driven by Brownian motion $B$ satisfies $\cF^X \equiv \cF^B$ $\bbP^*$-a.s.. In particular, the solution map $F^{\sigma,b}_G$ sending $B$ to $X$ lies in $\cT_{bcb}(\bbW^1, \kappa_G^{\sigma,b})$.
    \end{prop}

    \begin{rmk}
        \label{rmk:volterrakernel}
        A concrete example of $G$ satisfying Assumption \ref{asmp:svecoeff} and admitting a resolvent of first kind is $G(t - s) := c (t - s)^{\alpha - 1}$ for $c> 0, \alpha \in (1/2,1)$. Indeed, one can take $\beta \in (1, \frac{1}{2(1-\alpha)}), \theta = \alpha - 1/2$ as shown in \cite[Practitioner's corner]{jourdainpages2025convex} and $\mu(ds) := \tilde{c} s^{-\alpha} ds$ for some appropriate $\tilde{c} > 0$ by \cite{abijaberlarssonpulido2019affine}. We refer to \cite{abijaberlarssonpulido2019affine} for more examples of kernels admitting resolvents of first kind.
    \end{rmk}

    \begin{proof}
        Since \cite{jourdainpages2025convex} gives that $X$ is $\cF^B$-adapted, we trivially have $\cF^X \subseteq \cF^B$ $\bbP^*$-a.s.. Furthermore, we calculate \begin{align*}
            \int_{[0,t]} X_{t - s} \mu(ds) = & \int_{[0,t]} \int_0^{t - s} G( t - s -u) b(u,X_u) du \mu(ds) + \int_{[0,t]} \int_0^{t - s} G(t - s - u) \sigma(u,X_u) dB_u \mu(ds) \\
            = &\int_0^t \int_{[0,t - u]} G( t - u - s) \mu(ds) b(u,X_u) du + \int_0^t \int_{[0,t - u]} G( t - u - s) \mu(ds) \sigma(u,X_u) dB_u \\
            = &\int_0^t b(u,X_u) du + \int_0^t \sigma(u,X_u) dB_u.
        \end{align*} Here, the exchange of integrals in the second equality is justified by the stochastic Fubini theorem proved in \cite{veraar2012stochastic} because $b$ and $\sigma$ have linear growth and $\bbE |X|_\infty^p < \infty$, so that almost surely \[
        \int_{[0,t]} (\int_{[0,t-s]} G(t - s -u) |b(u,X_u)| du + (\int_0^{t - s} G(t - s - u)^2 \sigma(u,X_u)^2 du)^{1/2}) |\mu|(ds) < \infty.  
        \] In particular, \[
        B_t = \int_0^t \frac{1}{\sigma(u,X_u)} dM_u, 
        \] where $M_t := \int_{[0,t]} X_{t - s} \mu(ds) - \int_0^t b(s,X_s) ds$, from which it is clear that $\cF^B \subseteq \cF^X$ $\bbP^*$-a.s..
    \end{proof}

We conclude with a final proposition regarding the density of bicausal bijective maps among bicausal couplings. 
Related density results were established by \cite{pratelli2007equality} for classical transport, by \cite{beiglbockpammerschrott2025denseness} for causal transport in discrete time, and by \cite{beiglbocklacker2020denseness} in continuous time. 
In particular, the following extends \cite[Propositions $4.4$ and $4.5$]{cont2024causaltransportpathspace}: 

\begin{prop}
    \label{prop:densityofbcbmaps}
    Suppose $\bbQ_i, i = 1,2$ are such that $\cT_{bcb}(\bbW^n, \bbQ_i) \neq \emptyset$ for $i = 1,2$. Then, $\cT_{bcb}(\bbQ_1,\bbQ_2)$ is dense among $\Pi_{bc}(\bbQ_1,\bbQ_2)$ for the topology of weak convergence. In particular, for any continuous, bounded cost $c : \cW^n \times \cW^n \rightarrow \bbR$, we have \[
    \inf_{\pi \in \cT_{bcb}(\bbQ_1,\bbQ_2)} \bbE^\pi c = \inf_{\pi \in \cT_{bc}(\bbQ_1,\bbQ_2)} \bbE^\pi c = \inf_{\pi \in \Pi_{bc}(\bbQ_1,\bbQ_2)} \bbE^\pi c.
    \]
\end{prop} 

\begin{proof}
    Fix $T_i \in \cT_{bcb}(\bbW^n, \bbQ_i)$ for $i = 1,2$.
    By \cite[Proposition $4.3$]{cont2024causaltransportpathspace}, $\cT_{bcb}(\bbW^n, \bbW^n)$ is dense in $\Pi_{bc}(\bbW^n,\bbW^n).$ Hence, by \cite[Lemma $4.1$]{cont2024causaltransportpathspace}, we obtain that \[
    (T_1,T_2)_{\#}\Pi_{bc}(\bbW^n,\bbW^n) := \{(T_1,T_2)_{\#}\pi : \pi \in \Pi_{bc}(\bbW^n,\bbW^n) \} \subseteq \Pi_{bc}(\bbQ_1, \bbQ_2) \text{ is closed } 
    \] and \[
    \{ (T_1,T_2)_{\#}\pi : \pi \in \cT_{bcb}(\bbW^n,\bbW^n) \} \text{ is dense in } (T_1,T_2)_{\#}\Pi_{bc}(\bbW^n,\bbW^n).
    \] Hence, it is sufficient to verify that $(T_1,T_2)_{\#}\Pi_{bc}(\bbW^n,\bbW^n) = \Pi_{bc}(\bbQ_1, \bbQ_2)$ and $\{ (T_1,T_2)_{\#}\pi : \pi \in \cT_{bcb}(\bbW^n,\bbW^n) \} \subseteq \cT_{bcb}(\bbQ_1, \bbQ_2).$ The former is easily verified since $(T_1^{-1},T_2^{-1})_{\#}\Pi_{bc}(\bbQ_1, \bbQ_2) = \Pi_{bc}(\bbW^n, \bbW^n)$ and the latter follows from the fact that if $\pi \in \cT_{bcb}(\bbW^n,\bbW^n)$, then there exists some bicausal bijective map $\hat{T}$ such that $\hat{T} \in \cT_{bcb}(\bbW^n, \bbW^n)$ and \[
    (T_1, T_2)_{\#}(\identity,\hat{T})_{\#}\bbW^n = (T_1, T_2 \circ \hat{T})_{\#}\bbW^n = (\identity, T_2 \circ \hat{T} \circ T_1^{-1})_{\#} \bbQ_1.
    \] One verifies in a straightforward manner that if $T' := T_2 \circ \hat{T} \circ T_1^{-1}$ and $S' = T_1 \circ \hat{T}^{-1} \circ T_2^{-1}$, then \[
    S' \circ T' = \identity \quad \bbQ_1 \text{-a.s. and }    T' \circ S' = \identity \quad \bbQ_2 \text{-a.s.}.
    \] In other words, $(T_1, T_2)_{\#} \pi \in \cT_{bcb}(\bbQ_1, \bbQ_2)$ as required. The final assertion is straightforward to verify.
\end{proof}

As an application of the preceding proposition, let $\bbQ_i = \varrho^{1,b_i}$ be the law of the solution to the SDE \eqref{eqn:additivede} driven by fBM for $b_i \in \scrA^{1,1}$ satisfying \eqref{eqn:lipschitzsdecoeff} and $\int_0^1 |b_i(t,0)| dt< \infty$. Then, Proposition \ref{prop:additivegennoise} gives $\cT_{bcb}(\bbW^1, \bbQ_i) \neq \emptyset$ and hence, the values of the bicausal Monge and Kantorovich problems between $\bbQ_1$ and $\bbQ_2$ for bounded, continuous cost agree by Proposition \ref{prop:densityofbcbmaps}. One can derive similar statements using Proposition \ref{prop:dosssussmanbcbmap} or \ref{prop:svebcbmap}.

\subsection{Mean-variance portfolio selection stability}
\label{appendix:mvpselect}

For any $\bbP \in \cS\cM_2(\cW^n)$ and $k,\gamma \geq 0$, we define the mean-variance portfolio selection problem as \[
MV(\bbP) := \sup_{H \in \ocH_k} \bbE^\bbP (H\cdot X)_1 - \gamma \text{Var}^\bbP((H\cdot X)_1),
\] where $\text{Var}^\bbP(Z) := \bbE^\bbP (Z - \bbE^\bbP Z)^2.$ 

\begin{prop}
    \label{prop:mvprobstab}
    For any $\bbP, \bbQ \in \cS\cM_2(\cW^n)$, it holds that \[
    |MV(\bbP) - MV(\bbQ)| \leq (\sqrt{2}k + 4 \gamma k^2(\cA\cW_{sm}(\bbP, \delta_0)+\cA\cW_{sm}(\bbQ, \delta_0)))\cA\cW_{sm}(\bbP, \bbQ).
    \]
\end{prop}

\begin{proof}
    Fix $H \in \ocH_k$ and $\pi \in \Pi_{bc}(\bbP, \bbQ)$. By \cite[Lemma $3.7$]{backhoff2020adaptedwd}, there exists $G \in \ocH_k$ such that $\pi$-a.s. $G_t(Y) = \bbE^\pi[ H_t(X)| \cF^Y_t]$ for every $t \in [0,1]$ and \[
    (G(Y) \cdot Y)_1 = \bbE^\pi [ (H(X) \cdot Y)_1 | \cF^Y_1].
    \] Since $MV(\bbQ)$ is a supremum, we have \[
MV(\bbQ) \geq \bbE^\pi  (G(Y) \cdot Y)_1 - \gamma \text{Var}^\pi( (G(Y) \cdot Y)_1).
\] Furthermore, by properties of the conditional expectation, \[
\bbE^\bbQ  (G(Y) \cdot Y)_1 = \bbE^\pi (H(X) \cdot Y)_1 \text{ and } \text{Var}^\bbQ( (G(Y) \cdot Y)_1) \leq \text{Var}^\pi( (H(X) \cdot Y)_1),
\] from which we obtain \[
MV(\bbQ) \geq \bbE^\pi (H(X) \cdot Y)_1 - \gamma \text{Var}^\pi( (H(X) \cdot Y)_1).
\] 
The differences in mean and variance satisfy the following bounds: \[
|\bbE^\pi (H(X) \cdot Y)_1 - \bbE^\pi (H(X) \cdot X)_1| \leq \Vert (H(X) \cdot (Y - X))_1 \Vert_{L^2(\pi)} 
\] and \begin{align*}
     &|\text{Var}^\pi((H(X) \cdot Y)_1) - \text{Var}^\pi((H(X) \cdot X)_1)| \\
    \leq &
    2 (\Vert (H(X) \cdot X)_1 \Vert_{L^2(\pi)} + \Vert (H(X) \cdot Y)_1 \Vert_{L^2(\pi)}) \Vert (H(X) \cdot (Y - X))_1 \Vert_{L^2(\pi)}. 
\end{align*} Hence, taking infimum over $\pi \in \Pi_{bc}(\bbP, \bbQ)$ appropriately, we obtain \begin{multline*}
    MV(\bbQ) \\
    \geq \bbE^\bbP(H(X)\cdot X)_1 - \gamma \text{Var}^\bbP((H(X)\cdot X)_1) - (k \sqrt{2} + 4k^2 \gamma (\cA\cW_{sm}(\bbP, \delta_0) + \cA\cW_{sm}(\bbQ, \delta_0)))\cA\cW_{sm}(\bbP, \bbQ) 
\end{multline*} so that taking supremum over $H \in \ocH_k$ on the RHS yields \[
MV(\bbQ) \geq MV(\bbP) - (k \sqrt{2} + 4k^2 \gamma (\cA\cW_{sm}(\bbP, \delta_0) + \cA\cW_{sm}(\bbQ, \delta_0)))\cA\cW_{sm}(\bbP, \bbQ).
\] Reversing the roles of $\bbP$ and $\bbQ$ gives the desired.
\end{proof}


\subsection{Optional projection stability}
\label{appendix:optprojstab}

Optional projections also play a role in the stability of portfolio optimization: \citet[Theorem~3.2]{bayraktardolinskyidolinsky2020extended} show that the optional projection of a subsequential limit of asymptotically optimal strategies onto the investor's filtration is optimal, under their hypotheses. Here we study quantitative stability of optional projections with respect to changes in the underlying law.

Given any probability measure $\bbP \in \cP(\cW^n)$, the optional projection of a measurable process $Z : [0,1] \times \cW^n \rightarrow \bbR$ is an $\cF^X$-optional process $\optproj{\bbP}{Z} : [0,1] \times \cW^n \rightarrow \bbR$ such that for any $\cF^X$-stopping time $\tau$, \[
	\optproj{\bbP}{Z}_\tau = \bbE^\bbP (Z_\tau | \cF^X_\tau) \quad \bbP\text{-a.s.}.
\] 
If $\bbE^{\bbP}|Z|_\infty < \infty,$
then $\optproj{\bbP}{Z}$ exists by \cite[Theorem $5.1$]{hewangyan1992semimart}.
When necessary, given $\pi \in \cP(\cW^n \times \cW^n)$ and measurable process $\tZ : [0,1] \times \cW^n \times \cW^n \rightarrow \bbR$, the optional projection of $\tZ$ w.r.t. the filtration $\cF^{X,Y}$ jointly generated by marginal processes $X$ and $Y$ is denoted $\jointoptproj{\pi}{\tZ}$. Finally, we endow $\cM_b([0,1],\bbR^n)$ with the $\sigma$-algebra generated by the evaluation maps $\sigma(\{\omega \mapsto \omega(t) : t \in [0,1]\})$.

 \begin{prop}
    \label{prop:optionalprojAWbounds}
    Suppose $\bbP,\bbQ \in \cP_{|\cdot|_\infty,2}(\cW^n)$, $F : \cW^n \rightarrow \cW^n$ and measurable $ h : \cM_b([0,1], \bbR^n) \rightarrow \bbR$ satisfying \[
    |F(\omega) - F(\omega')|_\infty \leq K |\omega - \omega'|_\infty \text{ and }  |h(\overline{\omega}) - h(\overline{\omega}')| \leq K \sup_{t \in [0,1]} |\overline{\omega}_t - \overline{\omega}'_t|
    \] for all $\omega, \omega' \in \cW^n$ and $\overline{\omega}, \overline{\omega}' \in \cM_b([0,1], \bbR^n)$. Then, \[
    |\bbE^\bbP h(\optproj{\bbP}{F}) - \bbE^\bbQ h(\optproj{\bbQ}{F})| \leq 2K^2 \cA\cW_{|\cdot|_\infty, 2}(\bbP,\bbQ).
    \]
\end{prop}

\begin{rmk}
\label{rmk:optionalprojectionexistence}
    By $\optproj{\bbP}{F}$, we mean the optional projection of $F = (F_t)_{t \in [0,1]}$ given by $F_t(\omega) := F(\omega)(t)$ for all $t \in [0,1]$. 
    In particular, since $F$ is Lipschitz and $\bbP \in \cP_{|\cdot|_\infty, 2}(\cW^n)$, $\optproj{\bbP}{F}$ exists. 
\end{rmk}

\begin{proof}
    Fix $\pi \in \Pi_{bc}(\bbP,\bbQ)$. By Lemma \ref{lem:optionalprojcausalcondition} below, we have \[
    \optproj{\bbP}{F} \circ X = \jointoptproj{\pi}{F(X)} \text{ and } \optproj{\bbQ}{F} \circ Y = \jointoptproj{\pi}{F(Y)}.
    \] In particular, Lemma \ref{lem:optionalprojectionbound} gives \[
    \bbE^\bbP |\optproj{\bbP}{F}|_\infty^2 + \bbE^\bbQ |\optproj{\bbQ}{F}|_\infty^2 = \bbE |\jointoptproj{\pi}{F(X)}|^2_\infty + \bbE |\jointoptproj{\pi}{F(Y)}|^2_\infty \leq 4 (\bbE^\bbP |F|^2_\infty + \bbE^\bbQ |F|^2_\infty) < \infty
    \] and $h(\optproj{\bbP}{F})$ and $h(\optproj{\bbQ}{F})$ are well-defined. We compute using linearity of the optional projection and Lemma \ref{lem:optionalprojectionbound} below that \begin{align*}
        |\bbE^\bbP h(\optproj{\bbP}{F}) - \bbE^\bbQ h(\optproj{\bbQ}{F})| \leq &\bbE^\pi | h(\optproj{\bbP}{F} \circ X) - h(\optproj{\bbQ}{F} \circ Y)| \\
        \leq &K \bbE^\pi |\optproj{\bbP}{F}\circ X - \optproj{\bbQ}{F} \circ Y|_\infty\\
        = & K \bbE^\pi |\jointoptproj{\pi}{[F(X) - F(Y)]}|_\infty\\
        \leq &K \sqrt{\bbE^\pi |\jointoptproj{\pi}{[F(X) - F(Y)]}|^2_\infty} \\
        \leq & 2K \sqrt{\bbE^\pi |F(X) - F(Y)|^2_\infty} \\
        \leq &2K^2 \sqrt{\bbE^\pi |X - Y|^2_\infty}.
    \end{align*} Taking infimum over $\pi$ on the RHS yields the desired.
\end{proof}

\begin{lem}
    \label{lem:optionalprojectionbound}
    Fix $\pi \in \cP(\cW^n \times \cW^n)$ and any stochastic process $Z = (Z_t)_{t \in [0,1]}$ defined on $\cW^n \times \cW^n$ with continuous sample paths $\pi$-a.s. such that $\bbE^\pi |Z|^2_\infty < + \infty$. Then, 
    \[
    \bbE^\pi |\jointoptproj{\pi}{Z}|^2_\infty \leq 4\bbE^\pi |Z|^2_\infty.
    \]
\end{lem}

\begin{proof}
    We first denote $M := |Z|_\infty$ and 
    note that by definition, for any $\cF^{X,Y}$-stopping time $\tau$, we have $\pi$-a.s. that \[
    |\jointoptproj{\pi}{Z}_\tau| = |\bbE^\pi (Z_\tau | \cF^{X,Y}_\tau)| \leq \bbE^\pi(|Z_\tau||\cF^{X,Y}_\tau) \leq \bbE^\pi (M | \cF^{X,Y}_\tau) = \optproj{\pi}{M}_\tau. 
    \] This implies that $\pi$-a.s. $|\jointoptproj{\pi}{Z}| \leq \jointoptproj{\pi}{M}$ by \cite[Theorem $4.10$]{hewangyan1992semimart} (more precisely, the set $\{(\omega, \omega') \in \cW^n \times \cW^n : |\jointoptproj{\pi}{Z}_t(\omega,\omega')| > \jointoptproj{\pi}{M}_t(\omega,\omega') \text{ for some } t \in [0,1] \}$ is contained in a $\pi$-null set). This implies that $\pi$-a.s. $|\jointoptproj{\pi}{Z}|^2_\infty \leq |\jointoptproj{\pi}{M}|^2_\infty$ and it suffices now to show that $\bbE^\pi |\jointoptproj{\pi}{M}|^2_\infty \leq 4\bbE M^2.$ But this is Doob's $L^2$-martingale inequality since $\jointoptproj{\pi}{M}$ is the c\`adl\`ag version of the martingale  $(\bbE^\pi (M |\cF^{X,Y}_t))_{t \in [0,1]}$.
\end{proof}

\begin{lem}
    \label{lem:optionalprojcausalcondition}
    Let $\pi \in \Pi_c(\bbP, \bbQ)$ and $Z = (Z_t)_{t \in [0,1]}$ be a measurable process on $\cW^n$ with continuous sample paths $\bbP$-a.s. such that $\bbE^\bbP |Z|_\infty < \infty$. Then, \[
    \jointoptproj{\pi}{(Z \circ X)} = \optproj{\bbP}{Z} \circ X \quad \pi\text{-a.s.,} 
    \] where $X$ is the first marginal process on $\cW^n \times \cW^n$. More precisely, $\{ (\omega, \omega') \in \cW^n \times \cW^n :  \jointoptproj{\pi}{(Z \circ X)}_t \neq \optproj{\bbP}{Z}_t \circ X \text{ for some } t \in [0,1] \}$ is contained in a set of null $\pi$-measure.
\end{lem}

\begin{rmk}
    \label{rmk:cffiltrationenlargement}
    \cite{acciaio2020causalot} establish a related result showing that under causal couplings, the optional projection and dual optional projection (w.r.t. $\cF^X$) of c\`adl\`ag processes with integrable variation are actually indistinguishable. Here, we show that optional projections of continuous processes w.r.t. $\cF^X$ and $\cF^{X,Y}$ are indistinguishable.
\end{rmk}

\begin{proof}
    By definition of the optional projection, for any $\cF^{X,Y}$-stopping time $\tau$, we have \[
    \bbE^\pi((Z \circ X)_\tau|\cF^{X,Y}_\tau) = \jointoptproj{\pi}{(Z\circ X)}_\tau \quad \pi\text{-a.s..}
    \] Furthermore, since both $\jointoptproj{\pi}{(Z \circ X)}$ and $\optproj{\bbP}{Z} \circ X$ are optional (w.r.t. $\cF^{X,Y}$) processes, by \cite[Corollary $4.13$]{hewangyan1992semimart}, it suffices to show that \[
    \bbE^\pi((Z \circ X)_\tau) = \bbE^\pi (\jointoptproj{\pi}{(Z \circ X)}_\tau) = \bbE^\pi(\optproj{\bbP}{Z}_\tau \circ X).
    \] In order to show the preceding, we first consider a general measurable function $L : [0,1] \times \cW^n \rightarrow \bbR$ such that $\bbE^\pi |L(\tau,X)| < \infty$. Denoting $\lambda_\omega := \Law^\pi(\tau | X = \omega)$, we have $\Law^\pi(\tau,X| X = \omega) = \lambda_\omega \otimes \delta_\omega$. The quantile of $\lambda_\omega \in \cP([0,1])$ is given by \[
    Q_\omega(u) := \inf\{t \in [0,1] : \lambda_\omega([0,t]) \geq u\} = \inf\{ t \in [0,1] : \pi (\tau \leq t | X =\omega) \geq u \}, \quad u \in (0,1).
    \] We compute \begin{multline}
    \label{eqn:quantilerepresentation}
        \bbE^\pi L(\tau,X) = \int_{\cW^n} \int_{[0,1]} L(s,\omega) \lambda_\omega(ds) \bbP(d\omega) = \\
        \int_{\cW^n} \int_{[0,1]} L(Q_\omega(u),\omega) du \bbP(d\omega) = \int_0^1 \bbE^\bbP L(Q_X(u),X) du
    \end{multline} where the middle equality holds by inverse transform sampling. In particular, if we choose $L(t,\omega) = Z_t(\omega)$, we obtain \begin{equation}
        \label{eqn:quantilerepjointoptproj}
        \bbE^\pi ((Z \circ X)_\tau) = \int_0^1 \bbE^\bbP Z_{Q_X(u)} du.
    \end{equation}
    
    By causality of $\pi$, for fixed $u \in (0,1)$, $\omega \mapsto Q_\omega(u)$ is actually an $\prescript{\bbP}{}{\cF}^X$-stopping time. Indeed, by definition of the quantile, we have \[
    Q_\omega(u) \leq t \iff \lambda_\omega([0,t]) \geq u
    \] and hence, \[
    \{Q_\omega(u) \leq t\} = \{ \pi (\tau \leq t | X =\omega) \geq u \}.
    \]  We further have that since $\pi \in \Pi_c(\bbP,\bbQ)$, $\cF^X$ is conditionally independent of $\cF^{X,Y}_t$ given $\cF^X_t$ and since $\{\tau \leq t \} \in \cF^{X,Y}_t$, \[
    \pi(\tau \leq t | \cF^X) = \pi(\tau \leq t | \cF^X_t) \implies \omega \mapsto \pi(\tau \leq t | X = \omega) \text{ is } \prescript{\bbP}{}{\cF}^X_t\text{-measurable.}
    \] This implies that $\{Q_\omega(u) \leq t\} \in \prescript{\bbP}{}{\cF}^X_t$, which implies that $\omega \mapsto Q_\omega(u)$ is an $\prescript{\bbP}{}{\cF}^X$-stopping time. In particular, \[
    \bbE^\bbP (\optproj{\bbP}{Z}_{Q(u)}) = \bbE^\bbP (Z_{Q(u)}).
    \] Hence, choosing $L(t,\omega) := \optproj{\bbP}{Z}_t(\omega)$ in \eqref{eqn:quantilerepresentation}, we obtain \begin{equation}
        \label{eqn:quantilerepoptproj}
        \bbE^\pi(\optproj{\bbP}{Z}_\tau \circ X) = \int_0^1 \bbE^\bbP (\optproj{\bbP}{Z}_{Q_X(u)}) du = \int_0^1 \bbE^\bbP (Z_{Q_X(u)}) du.
    \end{equation} Comparing the RHSs of \eqref{eqn:quantilerepjointoptproj} and \eqref{eqn:quantilerepoptproj} yields the desired equality.
\end{proof}

\subsection{Miscellaneous}
\label{appendix:misc}

\begin{lem}
    \label{lem:infsupprob}
    Let $g : [0,\infty) \rightarrow [0,\infty)$ be a strictly increasing, strictly concave surjective function such that $g$ is smooth on $(0,\infty)$ and \begin{equation}
        \label{eqn:concavefnlimbehavior}
        \lim_{\zeta \downarrow 0} \frac{g(\zeta)}{\zeta} = \infty \text{ and } 
    \lim_{\zeta \rightarrow \infty} \frac{g(\zeta)}{\zeta} = 0.
    \end{equation} Then, for any $r > 0$, \[
    \inf_{\lambda > 0} \sup_{\zeta > 0}  \lambda  [g(\zeta) - r] - \zeta =  \sup_{\zeta > 0} \inf_{\lambda > 0}  \lambda  [g(\zeta) - r] - \zeta = -g^{-1}(r).
    \]  
\end{lem}

\begin{proof}
    Set $F(\lambda, \zeta) := \lambda [g(\zeta) - r] - \zeta$. We first show that \[
     \inf_{\lambda > 0} \sup_{\zeta > 0}  \lambda  [g(\zeta) - r] - \zeta =   \sup_{\zeta > 0} \inf_{\lambda > 0} \lambda  [g(\zeta) - r] - \zeta.
    \] By strict monotonicity and surjectivity of $g$, $g^{-1}$ exists and $g^{-1}(r) \in (0,\infty)$. For fixed $\lambda > 0$, $\zeta \mapsto F_\lambda(\zeta) := F(\lambda, \zeta)$ is strictly concave so its unique maximum $\zeta^*_\lambda \in (0,\infty)$ has characterization   \[
    \zeta = \zeta^*_\lambda \iff g'(\zeta) = \frac{1}{\lambda} \iff \lambda = \frac{1}{g'(\zeta)}.
    \] Hence, \[
    \sup_{\zeta > 0} \inf_{\lambda > 0} F(\lambda, \zeta) \leq \inf_{\lambda > 0} \sup_{\zeta > 0} F(\lambda,\zeta) \leq \sup_{\zeta > 0} F(\frac{1}{g'(g^{-1}(r))}, \zeta) = F(\frac{1}{g'(g^{-1}(r))}, g^{-1}(r)) = - g^{-1}(r).
    \] On the other hand, it is straightforward to check that \[
    \inf_{\lambda > 0} F(\lambda, \zeta) = \begin{cases}
        - \zeta &\text{ if } g(\zeta) \geq r \\
        -\infty &\text{ if } g(\zeta) < r.
    \end{cases}
    \] Since $g(\zeta) \geq r \iff \zeta \geq g^{-1}(r)$, we have \[
    \sup_{\zeta > 0} \inf_{\lambda > 0} F(\lambda, \zeta) = \sup_{\zeta \geq g^{-1}(r)} -\zeta = -g^{-1}(r),
    \] from which we deduce the desired identity.
\end{proof}

\begin{exam}
    \label{exam:NAFnotlipschitzbutregular}
    Let $(t, \omega) \in [0,1] \times \cD^1$ and define $F(t,\omega) := f(a(\max\{\omega_t^*, 0\})(\max\{\omega_t^*, 0\} - \omega_t))$, where $\omega^*_t := \sup_{s \leq t} \omega_s$, $f, a \in C^\infty_b(\bbR)$ such that \[
     a(0) = f'(0) = f''(0) = 0, \sup_{z} |z|(|f''(z)| + |f'''(z)|) < \infty, a'(0) \neq 0 \text{ and } f(a'(0)) \neq f(0) = 0.
    \]  Then, $F$ satisfies Assumption \ref{asmp:nonanticipativefnls} but is not Lipschitz w.r.t. $|\cdot|_\infty$. We first verify Assumption \ref{asmp:nonanticipativefnls}. It is clear that $F$ is non-anticipative. Since $(t,\omega) \mapsto \max\{ \omega^*_t,0 \},\omega_t$ is continuous w.r.t. the stopped path metric, it is straightforward to verify $F$, as the composition of continuous functions, is also continuous w.r.t. the stopped path metric. Furthermore, the horizontal derivative $DF$ is $0$. To see that $F$ is twice vertically differentiable, we need to show that $e \mapsto v(e)$ is twice differentiable at $0$ where  \[
    v(e) := f(a(\max\{\omega_{t-}^*, \omega_t + e, 0\})(\max\{\omega_{t-}^*, \omega_t + e, 0\} - \omega_t - e))
    \] and $\omega_{t-}^* := \sup_{s < t} \omega_s$. It is clear that if $\omega_t < \max\{\omega^*_{t-},0\}$, then for $e$ small enough, we have \[
    v(e) = f(a(\max\{\omega_{t-}^*, 0\})(\max\{\omega_{t-}^*, 0\} - \omega_t - e))
    \] and \begin{gather*}
        v'(0) = -f'(a(\max\{\omega_{t}^*, 0\})(\max\{\omega_{t}^*, 0\} - \omega_t))a(\max\{\omega_{t}^*, 0\}) \\
        v''(0) = f''(a(\max\{\omega_{t}^*, 0\})(\max\{\omega_{t}^*, 0\} - \omega_t))a(\max\{\omega_{t}^*, 0\})^2.
    \end{gather*} If $\omega_t > \max \{ \omega^*_{t-}, 0 \}$, then $v(e)$ is constant for $e$ small enough and $v'(0) = v''(0) = 0$. If $\omega_t = \max\{\omega^*_{t-},0\}$, we compute the left and right derivatives of $v$ separately. For $e > 0$, we have \[
    v(e) = f(a(\omega_t + e)(0)) = f(0) \implies v'_+(0) = 0.
    \] Otherwise, for $e < 0$, we have \[
     v(e) = f(a(\max\{\omega_{t-}^*, 0\})(\max\{\omega_{t-}^*, 0\} - \omega_t - e)) \implies v'_-(0) = -f'(0)a(\omega_t) = 0.
    \] Since the left and right derivatives agree, we conclude that $v'(0) = 0$ and \[
    \nabla F(t,\omega) = - f'(a(\max\{\omega_{t}^*, 0\})(\max\{\omega_{t}^*, 0\} - \omega_t))a(\max\{\omega_{t}^*, 0\}).
    \] Similar computations yield \[
    \nabla^2 F(t,\omega) = f''(a(\max\{\omega_{t}^*, 0\})(\max\{\omega_{t}^*, 0\} - \omega_t))a(\max\{\omega_{t}^*, 0\})^2.
    \] Hence, by continuity of $f'',f'$ and $a$, we see that $DF, \nabla F$ and $\nabla^2 F$ are boundedness preserving and continuous w.r.t. the stopped path metric. To see that the vertical derivatives are Lipschitz (boundedness is clear), it is straightforward to verify that the partial derivatives of the smooth functions $(x,y) \mapsto - f'(a(x)y)a(x), f''(a(x)y)a(x)^2$ are bounded by the assumption on $f$. Hence, $F$ satisfies Assumption \ref{asmp:nonanticipativefnls}. Nonetheless, $F$ is not even uniformly continuous on $\cW^1$. 
    Indeed, we take $\omega_n,\omega'_n \in \cW^1$ to be the piecewise linear functions such that $\omega_n(1/2) = 0, \omega_n(1) = -n$ and $\omega'_n(1/2) = 1/n, \omega'_n(1) = -n$. We have that $|\omega_n - \omega'_n|_\infty \leq \frac{1}{n}$, but \[
    F(1,\omega_n) = f(0) \text{ and } F(1, \omega'_n) = f(a(n^{-1})(n^{-1} + n)) \stackrel{n \rightarrow \infty}{\rightarrow} f(a'(0)) \neq f(0).
    \] 
\end{exam}

\bibliography{mybiblio}




\end{document}